\documentclass[10pt,journal,compsoc]{IEEEtran}

\ifCLASSOPTIONcompsoc
  \usepackage[nocompress]{cite}
\else
  \usepackage{cite}
\fi
\ifCLASSINFOpdf
\else
\fi
\usepackage{amsmath}
\usepackage{amssymb}
\usepackage{mathtools}
\usepackage{amsthm}

\usepackage{amsfonts}
\usepackage{multirow}
\usepackage{multicol}

\usepackage{booktabs}
\usepackage{color}
\usepackage{algorithm}
\usepackage{algorithmic}

\usepackage{tcolorbox}
\usepackage{graphicx}
\usepackage{subfig}

\newtheorem{theorem}{Theorem}
\newtheorem{lemma}{Lemma}

\newtheorem{assumption}{Assumption}
\newtheorem{remark}{Remark}

\begin{document}
%
\title{BiAdam: Fast Adaptive Bilevel Optimization Methods}
%
%
%
%

\author{Feihu~Huang,~Junyi Li,~Shangqian~Gao 
\IEEEcompsocitemizethanks{\IEEEcompsocthanksitem Feihu Huang is with College of Computer Science and Technology, Nanjing University of Aeronautics and Astronautics, Nanjing, China; and also with MIIT Key Laboratory of Pattern Analysis and Machine Intelligence, Nanjing, China.
E-mail: huangfeihu2018@gmail.com \\
Junyi Li is with Department of Electrical and Computer Engineering, University of Pittsburgh, Pittsburgh, USA.
E-mail: junyili.ai@gmail.com \\
Shangqian Gao
is with Department of Electrical and Computer Engineering, University of Pittsburgh, Pittsburgh, USA.
E-mail: shg84@pitt.edu
}}

%
%

\markboth{Journal of \LaTeX\ Class Files,~Vol.~14, No.~8, August~2015}%
{Shell \MakeLowercase{\textit{et al.}}: Bare Advanced Demo of IEEEtran.cls for IEEE Computer Society Journals}
%



\IEEEtitleabstractindextext{%
\begin{abstract}
Bilevel optimization recently has attracted increased interest in machine
learning due to its many applications such as hyper-parameter optimization and meta learning.
Although many bilevel methods recently have been proposed,
these methods do not consider using adaptive learning rates. It is well known that adaptive learning rates can accelerate
optimization algorithms.
To fill this gap, in the paper, we propose a novel fast adaptive bilevel framework to solve stochastic bilevel optimization problems that
the outer problem is possibly nonconvex and the inner problem is strongly convex.
Our framework uses unified adaptive matrices including many types of adaptive learning rates, and
can flexibly use the momentum and variance reduced techniques.
In particular, we provide a useful convergence analysis framework for the bilevel optimization.
Specifically, we propose a fast single-loop adaptive bilevel optimization (BiAdam) algorithm,
which achieves a sample complexity of $\tilde{O}(\epsilon^{-4})$ for finding an $\epsilon$-stationary solution.
Meanwhile, we propose an accelerated version of BiAdam algorithm (VR-BiAdam), which reaches
the best known sample complexity of $\tilde{O}(\epsilon^{-3})$.
To the best of our knowledge, we first
study the adaptive bilevel optimization methods with adaptive learning rates.
Experimental results on data hyper-cleaning and hyper-representation learning tasks demonstrate the efficiency of our algorithms.
\end{abstract}

}

\maketitle

\IEEEdisplaynontitleabstractindextext

%
\IEEEpeerreviewmaketitle

\section{Introduction}
\IEEEPARstart{B}{ilevel} optimization is known as a class of popular hierarchical optimization, which has been applied to a wide range of machine learning problems
such as hyperparameter optimization~\cite{shaban2019truncated}, meta-learning~\cite{ji2021bilevel,liu2021investigating} and policy optimization~\cite{hong2020two}.
In the paper, we consider solving the following stochastic bilevel optimization problem, defined as
\begin{align}
 \min_{x \in \mathcal{X}} & \ F(x):=\mathbb{E}_{\xi\sim \mathcal{D}}\Big[f\big(x,y^*(x);\xi\big)\Big]  & \mbox{(Outer)} \label{eq:1} \\
 \mbox{s.t.} & \ y^*(x) \in \arg\min_{y\in \mathcal{Y}} \ \mathbb{E}_{\zeta\sim \mathcal{H}}\Big[g(x,y;\zeta)\Big],   &  \mbox{(Inner)}
\end{align}
where $F(x)=f(x,y^*(x))=\mathbb{E}_{\xi}\big[f(x,y^*(x);\xi)\big]$ is a differentiable and possibly nonconvex function,
and $g(x,y)=\mathbb{E}_{\zeta}\big[g(x,y;\zeta)\big]$ is a differentiable and strongly convex
function in variable $y$, and $\xi$ and $\zeta$ are random variables follow unknown distributions $\mathcal{D}$ and
$\mathcal{H}$, respectively.
Here $\mathcal{X}\subseteq \mathbb{R}^d$ and $\mathcal{Y}\subseteq \mathbb{R}^p$ are convex closed sets.
Problem (\ref{eq:1}) involves many machine learning problems with a hierarchical structure, which include hyper-parameter
optimization~\cite{franceschi2018bilevel}, meta-learning~\cite{franceschi2018bilevel}, policy optimization~\cite{hong2020two} and neural network architecture search \cite{liu2018darts}.
Here, we specifically provide two popular machine learning applications
that can be formulated as Problem \eqref{eq:1}.

\begin{table*}
  \centering
  \caption{ Sample complexity of the representative bilevel optimization methods for finding an $\epsilon$-stationary point of the bilevel problem (\ref{eq:1}), i.e., $\mathbb{E}\|\nabla F(x)\| \leq \epsilon$ or its equivalent variants.
  \textbf{BSize} denotes mini-batch size; \textbf{ALR} denotes adaptive learning rate. $C(x,y)$ denotes the constraint sets in $x$ and $y$, where \textbf{Y} denotes the fact that there exists a convex constraint on variable, otherwise is \textbf{N}. \textbf{DD} denotes dimension dependence in the gradient estimators, and $p$ denotes the dimension of variable $y$.
  \textbf{1} denotes Lipschitz continuous of $\nabla_xf(x,y;\xi)$, $\nabla_yf(x,y;\xi)$, $\nabla_yg(x,y;\zeta)$, $\nabla^2_{xy}g(x,y;\zeta)$ and $\nabla^2_{yy}g(x,y;\zeta)$ for all $\xi,\zeta$;
  \textbf{2} denotes Lipschitz continuous of $\nabla_xf(x,y)$, $\nabla_yf(x,y)$, $\nabla_yg(x,y)$, $\nabla^2_{xy}g(x,y)$ and $\nabla^2_{yy}g(x,y)$; \textbf{3} denotes bounded stochastic partial derivatives $\nabla_yf(x,y;\xi)$ and $\nabla^2_{xy}g(x,y;\zeta)$; \textbf{4} denotes bounded stochastic partial derivatives $\nabla_xf(x,y;\xi)$, and $\nabla^2_{yy}g(x,y;\zeta)$; \textbf{5} denotes the bounded true partial derivatives $\nabla_y f(x,y)$ and $\nabla^2_{xy} g(x,y)$; \textbf{6} denotes Lipschitz continuous of function $f(x,y;\xi)$; \textbf{7} denotes
  $g(x,y;\zeta)$ is $L_g$-smooth and $\mu$-strongly convex function \emph{w.r.t.} $y$ for all $\zeta$; \textbf{8} denotes
  $g(x,y)$ is  $L_g$-smooth and $\mu$-strongly convex function \emph{w.r.t.} $y$.
   }
  \label{tab:1}
  \resizebox{1.0\textwidth}{!}{
  \begin{tabular}{c|c|c|c|c|c|c|c|c}
  \hline
 \textbf{Algorithm} & \textbf{Reference} &  \textbf{Complexity}  & \textbf{BSize} & \textbf{Loop} & $C(x,y)$ & \textbf{DD} & \textbf{ ALR} & \textbf{Conditions}\\ \hline
  BSA & \cite{ghadimi2018approximation} & $O(\epsilon^{-6})$  & $\tilde{O}(1)$ & Double & Y, N & $p^2$ & & \textbf{2},\ \textbf{5},\ \textbf{7}  \\ \hline
  TTSA & \cite{hong2020two} & $\tilde{O}(\epsilon^{-5})$  & $\tilde{O}(1)$ &  Single & Y, N & $p^2$ &  & \textbf{1},\ \textbf{3},\ \textbf{7}  \\ \hline
  stocBiO & \cite{ji2021bilevel} & $O(\epsilon^{-4})$ & $\tilde{O}(\epsilon^{-2})$ & Double & N, N & $p^2$ &  & \textbf{1},\ \textbf{6},\ \textbf{7} \\ \hline
  STABLE & \cite{chen2022single} & $\tilde{O}(\epsilon^{-4})$  & $O(1)$ & Single & N, N & $p^3$ & & \textbf{1},\ \textbf{3},\ \textbf{4},\ \textbf{8} \\ \hline
  SMB & \cite{guo2021stochastic} & $\tilde{O}(\epsilon^{-4})$  & $\tilde{O}(1)$ & Single & N, Y & $p^2$ &  & \textbf{2},\ \textbf{5},\ \textbf{7} \\ \hline
  SUSTAIN & \cite{khanduri2021near} & $\tilde{O}(\epsilon^{-3})$ & $\tilde{O}(1)$  & Single  & N, N & $p^2$ & & \textbf{1},\ \textbf{3},\ \textbf{7}\\ \hline
  SVRB & \cite{guo2021randomized} & $\tilde{O}(\epsilon^{-3})$  & $O(1)$ & Single & N, N & $p^3$ &  & \textbf{1},\ \textbf{5},\ \textbf{8} \\ \hline
  MRBO & \cite{yang2021provably} & $\tilde{O}(\epsilon^{-3})$ & $\tilde{O}(1)$  & Single  & N, N & $p^2$ & & \textbf{1},\ \textbf{6},\ \textbf{7}\\ \hline
  VRBO & \cite{yang2021provably} & $\tilde{O}(\epsilon^{-3})$  & $\tilde{O}(\epsilon^{-2})$ & Double & N, N & $p^2$ &  & \textbf{1},\ \textbf{6},\ \textbf{7} \\ \hline
  BiAdam & Ours & $\tilde{O}(\epsilon^{-4})$ & $\tilde{O}(1)$ & Single & Y/N, Y/N & $p^2$ & $\surd$ & \textbf{2},\ \textbf{5},\ \textbf{7} \\ \hline
  VR-BiAdam & Ours & $\tilde{O}(\epsilon^{-3})$ & $\tilde{O}(1)$ & Single & Y/N, Y/N & $p^2$ & $\surd$ & \textbf{1},\ \textbf{5},\ \textbf{7} \\ \hline
  \end{tabular}
  }
\end{table*}

\textbf{Hyper-parameter optimization}. The goal of hyper-parameter optimization is to search for the optimal hyper-parameters $x\in \mathcal{X} \subseteq \mathbb{R}^d$ (e.g., regularization coefficient, learning rate, neural network architecture). These optimal hyper-parameters
used in training a model $\theta \in \Theta \subseteq \mathbb{R}^d$ on the training set aims to make the learned model achieve the
lowest risk on the validation set. Specifically, consider hyper-parameter as the regularization coefficient as in \cite{franceschi2018bilevel},
we solve the following bilevel optimization problem:
\begin{align}
 & \min_{x\in \mathcal{X}} \  \mathbb{E}_{\xi\sim \mathcal{D}_{val}}\Big[\ell\big(\theta^*(x);\xi\big)\Big]  \nonumber   \\
 & \mbox{s.t.} \ \theta^*(x) \in \arg\min_{\theta \in \Theta} \bigg\{ \mathbb{E}_{\xi \sim \mathcal{D}_{tra} }\Big[\ell(\theta;\xi)\Big]
  + \lambda\sum_{i=1}^d x_i\theta_i^2 \bigg\}, \nonumber
\end{align}
where $\lambda>0$ is a tuning parameter, and $\ell(\theta;\xi)$ is the loss function on sample $\xi$, and $\mathcal{D}_{tra}$ and
$\mathcal{D}_{val}$ denote the training and validation datasets, respectively.

\textbf{ Meta-Learning}. Model-agnostic meta learning (MAML) is an effective learning paradigm, which is to find a good model
to achieve the best performance for individual tasks by using more experiences.
Consider the few-shot meta-learning problem with $m$ tasks $\{\mathcal{T}_i\}_{i=1}^m$, each task $\mathcal{T}_i$ has
training and test datasets $\mathcal{D}_{tr}^i$ and $\mathcal{D}_{te}^i$.
As in \cite{ji2021bilevel,guo2021randomized}, the MAML can be formulated as the following bilevel optimization problem:
\begin{align}
 \min_{\theta \in \Theta} & \ \frac{1}{m}\sum_{i=1}^m \frac{1}{|\mathcal{D}^i_{te}|}\sum_{\xi\in \mathcal{D}^i_{te}} \mathcal{L}(\theta,\theta^{i*};\xi)  \nonumber \\
 \mbox{s.t.} & \ \theta^{i*} \in \arg\min_{\theta^i\in \mathbb{R}^d} \bigg\{ \frac{1}{|\mathcal{D}^i_{tr}|} \sum_{\xi\in \mathcal{D}^i_{tr}} \mathcal{L}(\theta,\theta^{i};\xi)  + \frac{\lambda}{2}\|\theta-\theta^i\|^2 \bigg\}   \label{eq:6},
\end{align}
where $\theta^i$ is the model parameter of the $i$-th task for all $i\in [m]$, and $\theta$ is the shared model parameter.
Here $\mathcal{L}(\cdot)$ denotes loss function, and $\lambda\geq 0$ is a tuning parameter. Given a sufficiently large $\lambda$, the  inner problem \eqref{eq:6} is strongly-convex.

Since bilevel optimization has been widely applied in machine learning,
some works recently have been begun to study the bilevel optimization. For example,  \cite{ghadimi2018approximation,ji2021bilevel}
proposed a class of double-loop methods to solve the problem (\ref{eq:1}). However,
to obtain an accurate estimate, the BSA in \cite{ghadimi2018approximation} needs to
solve the inner problem to a high accuracy, and the stocBiO in \cite{ji2021bilevel} requires large batch-sizes in solving the
inner problem. \cite{hong2020two} proposed a class of single-loop methods to solve the bilevel problems.
Subsequently, \cite{khanduri2021near,guo2021randomized,yang2021provably,chen2022single} presented some accelerated single-loop methods by using the momentum-based variance reduced technique of STORM \cite{cutkosky2019momentum}.
Although these methods can effectively solve the bilevel problems, they
do not consider using the adaptive learning rates and only consider the bilevel problems under unconstrained setting.
Since using generally different learning rates for the inner and outer problems
to ensure the convergence of bilevel optimization problems, we will consider using different adaptive learning rates
for the inner and outer problems with convergence guarantee. Clearly, this can not follow the exiting adaptive methods
for single-level problems. Thus, there exists a natural question:
\begin{center}
\begin{tcolorbox}
\textbf{ How to design the effective optimization methods with adaptive learning rates for the bilevel problems ? }
\end{tcolorbox}
\end{center}

In the paper, we provide an affirmative answer to this question and propose
a class of fast single-loop adaptive bilevel optimization methods based on unified adaptive matrices,
which including many types of adaptive learning rates.
Moreover, our framework can flexibly use the momentum and variance reduced techniques.
Our main \textbf{contributions} are summarized as follows:
\begin{itemize}
\setlength{\itemsep}{0pt}
\item[1)] We propose a fast single-loop adaptive bilevel optimization algorithm (BiAdam) based on the basic momentum technique, which achieves a sample complexity of $\tilde{O}(\epsilon^{-4})$ for finding an $\epsilon$-stationary solution.
\item[2)] We further propose a single-loop accelerated version of BiAdam algorithm (VR-BiAdam) by using the momentum-based
variance reduced technique, which reaches the best known sample complexity of $\tilde{O}(\epsilon^{-3})$.
\item[3)] Moreover, we provide a useful convergence analysis framework for both the constrained and unconstrained bilevel programming under some mild conditions (Please see Table \ref{tab:1}).
\item[4)] The  experimental results on hyper-parameter learning demonstrate the efficiency of the proposed algorithms.
\end{itemize}

\section{Related Works}
In this section, we overview the existing bilevel optimization methods and
adaptive methods for single-level optimization,  respectively.

\subsection{ Bilevel Optimization Methods }
Bilevel optimization has shown successes in many machine learning problems with hierarchical structures such as policy optimization~\cite{hong2020two},
model-agnostic meta-learning~\cite{liu2021investigating} and adversarial training~\cite{zhang2021revisiting}.
Thus, its researches have become active in the machine learning community, and
some bilevel optimization methods recently have been proposed.
For example, one class of successful methods \cite{colson2007overview,kunapuli2008classification} are to reformulate the bilevel problem as a single-level problem by replacing
the inner problem by its optimality conditions. Another class of successful methods \cite{ghadimi2018approximation,hong2020two,ji2021bilevel,
chen2022single,khanduri2021near,guo2021randomized,chen2021closing,liu2021towards,liu2022general,li2021fully} for bilevel optimization
are to iteratively approximate the (stochastic) gradient of the outer problem either
in forward or backward. Specifically, \cite{liu2022general} proposed a general gradient-based descent aggregation framework for
bilevel optimization. Moreover, the non-asymptotic analysis of these gradient-based
bilevel optimization methods has been recently studied.
For example, \cite{ghadimi2018approximation} first studied the sample complexity of $O(\epsilon^{-6})$ of the proposed double-loop algorithm
for the bilevel problem \eqref{eq:1} (Please see Table \ref{tab:1}). Subsequently, \cite{ji2021bilevel} proposed an accelerated double-loop algorithm that
reaches the sample complexity of $O(\epsilon^{-4})$ relying on large batches.
At the same time, \cite{hong2020two} studied a single-loop algorithm that
reaches the sample complexity of $O(\epsilon^{-5})$ without relying on large batches.
Moreover, \cite{khanduri2021near,guo2021randomized,yang2021provably} proposed a class of accelerated single-loop methods for the bilevel problem \eqref{eq:1}
by using momentum-based variance reduced technique,
which achieve the best known sample complexity of $O(\epsilon^{-3})$. More recently,
\cite{ji2022lower} studied the lower bound of bilevel optimization methods.
Meanwhile, \cite{chen2022fast} studied the nonconvex and nonsmooth bilevel Optimization.

\subsection{ Adaptive Gradient Methods }
Adaptive gradient methods recently have been shown great successes in current machine learning problems such as training Deep Neural Networks (DNNs).
Recently, thus many adaptive gradient methods \cite{duchi2011adaptive,kingma2014adam,loshchilov2018decoupled,zhuang2020adabelief} have been developed and studied. For example,
Adagrad \cite{duchi2011adaptive} is the first
adaptive gradient method that shows good performances under the sparse gradient setting.  One variant of Adagrad, i.e.,
Adam \cite{kingma2014adam} is a very popular adaptive gradient method and
basically is a default method of choice for training DNNs.
Subsequently, some variants of Adam algorithm \cite{reddi2019convergence,chen2019convergence}
have been developed and studied, and especially they have  convergence guarantee under the nonconvex setting.
At the same time, some adaptive gradient methods \cite{loshchilov2018decoupled,chen2018closing,zhuang2020adabelief} have been presented to improve the generalization performance of Adam algorithm.
The norm version of AdaGrad (i.e., AdaGrad-Norm) \cite{ward2019adagrad} has been presented
to accelerate AdaGrad without sacrificing generalization.
Moreover, some accelerated adaptive gradient methods such as STORM \cite{cutkosky2019momentum} and SUPER-ADAM \cite{huang2021super} have been proposed by using variance-reduced technique.
Meanwhile, \cite{huang2021super,guo2021novel} studied the convergence analysis framework for adaptive gradient methods.

\section{Preliminaries}
\subsection{Notations}
$\mathcal{U}\{1,2,\cdots,K\}$ denotes a uniform distribution over a discrete set $\{1,2,\cdots,K\}$.
$\|\cdot\|$ denotes the $\ell_2$ norm for vectors and spectral norm for matrices.
$\langle x,y\rangle$ denotes the inner product of two vectors $x$ and $y$. For vectors $x$ and $y$, $x^r \ (r>0)$ denotes the element-wise
power operation, $x/y$ denotes the element-wise division and $\max(x,y)$ denotes the element-wise maximum. $I_{d}$ denotes a $d$-dimensional identity matrix. $A \succ 0$ denotes that the matrix $A$ is positive definite.
Given function $f(x,y)$, $f(x,\cdot)$ denotes  function \emph{w.r.t.} the second variable with fixing $x$,
and $f(\cdot,y)$ denotes function \emph{w.r.t.} the first variable
with fixing $y$.
$a=O(b)$ denotes that $a\leq C b$ for some constant $C>0$. The notation $\tilde{O}(\cdot)$ hides logarithmic terms.
Given a convex closed set $\mathcal{X}$, we define a projection operation to $\mathcal{X}$ as
$\mathcal{P}_{\mathcal{X}}(z) = \arg\min_{x\in \mathcal{X}}\frac{1}{2}\|x-z\|^2$.

\subsection{ Some Mild Assumptions}
In this subsection, we give some mild assumptions on the problem (\ref{eq:1}).
\begin{assumption} \label{ass:1}
For any $x$ and $\zeta$, $g(x,y;\zeta)$ is $L_g$-smooth and $\mu$-strongly convex function, i.e., $L_g I_p\succeq \nabla^2_{yy}g(x,y;\zeta) \succeq \mu I_p$.
\end{assumption}
\begin{assumption} \label{ass:2}
For functions $f(x,y)$ and $g(x,y)$ for all $x\in \mathcal{X}$ and $y\in \mathcal{Y}$, we assume the following conditions hold:
$\nabla_x f(x,y)$ and $\nabla_y f(x,y)$ are $L_{f}$-Lipschitz continuous,
$\nabla_y g(x,y)$ is $L_{g}$-Lipschitz continuous, $\nabla^2_{xy}g(x,y)$ is $L_{gxy}$-Lipschitz continuous,
$\nabla^2_{yy}g(x,y)$ is $L_{gyy}$-Lipschitz continuous. For example, for all $x,x_1,x_2 \in \mathcal{X}$ and $y,y_1,y_2\in \mathcal{Y}$, we have
\begin{align}
 & \|\nabla_x f(x_1,y)-\nabla_x f(x_2,y)\| \leq L_f\|x_1-x_2\|, \nonumber \\
 & \|\nabla_x f(x,y_1)-\nabla_x f(x,y_2)\| \leq L_f\|y_1-y_2\|. \nonumber
\end{align}
\end{assumption}
\begin{assumption} \label{ass:3}
For functions $f(x,y;\xi)$ and $g(x,y;\zeta)$ for all $x\in \mathcal{X}$, $y\in \mathcal{Y}$, $\xi$ and $\zeta$, we assume the following conditions hold:
$\nabla_x f(x,y;\xi)$ and $\nabla_y f(x,y;\xi)$ are $L_{f}$-Lipschitz continuous,
$\nabla_y g(x,y;\zeta)$ is $L_{g}$-Lipschitz continuous, $\nabla^2_{xy}g(x,y;\zeta)$ is $L_{gxy}$-Lipschitz continuous,
$\nabla^2_{yy}g(x,y;\zeta)$ is $L_{gyy}$-Lipschitz continuous. For example, for all $x,x_1,x_2 \in \mathcal{X}$ and $y,y_1,y_2\in \mathcal{Y}$, we have
\begin{align}
 & \|\nabla_x f(x_1,y;\xi)-\nabla_x f(x_2,y;\xi)\| \leq L_f\|x_1-x_2\|, \nonumber \\
 &\|\nabla_x f(x,y_1;\xi)-\nabla_x f(x,y_2;\xi)\| \leq L_f\|y_1-y_2\|. \nonumber
\end{align}
\end{assumption}
\begin{assumption} \label{ass:4}
The partial derivatives $\nabla_y f(x,y)$ and $\nabla^2_{xy} g(x,y)$ are bounded, i.e., $\|\nabla_y f(x,y)\|^2 \leq C^2_{fy}$
and $\|\nabla^2_{xy} g(x,y)\|^2 \leq C^2_{gxy}$.
\end{assumption}
\begin{assumption} \label{ass:5}
Stochastic functions $f(x,y;\xi)$ and $g(x,y;\zeta)$ have unbiased stochastic partial derivatives
with bounded variance, e.g.,
\begin{align}
& \mathbb{E}[\nabla_x f(x,y;\xi)] = \nabla_x f(x,y), \nonumber \\
& \mathbb{E}\|\nabla_x f(x,y;\xi) - \nabla_x f(x,y) \|^2 \leq \sigma^2. \nonumber
\end{align}
The same assumptions hold for $\nabla_y f(x,y;\xi)$, $\nabla_y g(x,y;\zeta)$, $\nabla^2_{xy} g(x,y;\zeta)$ and $\nabla^2_{yy} g(x,y;\zeta)$.
\end{assumption}

Assumptions \ref{ass:1}-\ref{ass:5} are commonly used in stochastic bilevel optimization problems \cite{ghadimi2018approximation,hong2020two,ji2021bilevel,chen2022single,khanduri2021near}.
Note that Assumption \ref{ass:3} is clearly stricter than Assumption \ref{ass:2}. For example,
given Assumption \ref{ass:3}, we have $\|\nabla_x f(x_1,y)-\nabla_x f(x_2,y)\|=\|\mathbb{E}[\nabla_x f(x_1,y;\xi)-\nabla_x f(x_2,y;\xi)]\|
\leq \mathbb{E}\|\nabla_x f(x_1,y;\xi)-\nabla_x f(x_2,y;\xi)\| \leq L_f\|x_1-x_2\|\|$ for any $x,y,\xi$.
At the same time, based on Assumptions \ref{ass:4}-\ref{ass:5}, we also have
$\|\nabla_y f(x,y;\xi)\|^2 = \|\nabla_y f(x,y;\xi)-\nabla_y f(x,y)-\nabla_y f(x,y)\|^2 \leq 2\|\nabla_y f(x,y;\xi)-\nabla_y f(x,y)\|^2 + 2\|\nabla_y f(x,y)\|^2 \leq 2\sigma^2 + 2C^2_{fy}$ and $\|\nabla^2_{xy} g(x,y;\zeta)\|^2 \leq 2\sigma^2 + 2C^2_{gxy}$. Thus we argue that under Assumption \ref{ass:5}, the bounded $\nabla_y f(x,y)$ and $\nabla^2_{xy} g(x,y)$ are not milder than the bounded $\nabla_y f(x,y;\xi)$ and $\nabla^2_{xy} g(x,y;\zeta)$ for all $\xi$ and $\zeta$.

\subsection{ Bilevel Optimization }
In this subsection, we review the basic first-order method for solving the problem (\ref{eq:1}).
Naturally, we give the following iteration to update the variables $x,y$: at the $t$-th step
\begin{align}
 & y_{t+1} = \mathcal{P}_{\mathcal{Y}}\Big(y_t - \lambda \nabla_y g(x_t,y_t)\Big), \nonumber \\
 & x_{t+1} = \mathcal{P}_{\mathcal{X}}\Big(x_t - \gamma \nabla_x f(x_t,y^*(x_t))\Big), \nonumber
\end{align}
where $\lambda>0$ and $\gamma>0$ denote the step sizes.
Clearly, if there does not exist a closed form solution of the inner problem in the problem (\ref{eq:1}), i.e., $y_{t+1}\neq y^*(x_t)$, we can not
easily obtain the gradient $\nabla F(x_t) = \nabla f(x_t,y^*(x_t))$. Thus, one of key points in solving the problem (\ref{eq:1}) is to estimate the gradient $\nabla F(x_t)$.

\begin{lemma} \label{lem:0}
(Lemma 2.1 in \cite{ghadimi2018approximation})
Under the above Assumption \ref{ass:2}, we have, for any $x\in \mathcal{X}$
\begin{align}
 \nabla F(x)  & = \nabla_x f(x,y^*(x)) - \nabla^2_{xy} g(x,y^*(x))\big[\nabla^2_{yy}g(x,y^*(x))\big]^{-1} \nonumber \\
 & \quad \cdot \nabla_y f(x,y^*(x)). \nonumber
\end{align}
\end{lemma}
From the above Lemma \ref{lem:0}, it is natural to use the following form to estimate $\nabla F(x)$, defined as,
\begin{align}
 \bar{\nabla} f(x,y) & = \nabla_xf(x,y) - \nabla^2_{xy}g(x,y) \big(\nabla^2_{yy}g(x,y)\big)^{-1}  \nonumber \\
 & \quad \cdot \nabla_yf(x,y), \quad \forall x\in \mathcal{X}, y\in \mathcal{Y} \nonumber
\end{align}
Note that although the inner problem of the problem (\ref{eq:1}) is a constrained optimization, we
assume that the optimal condition of the inner problem still is $\nabla_y g(x,y^*(x))=0$ and $y^*(x)\in \mathcal{Y}$.

\begin{lemma} \label{lem:1}
(Lemma 2.2 in \cite{ghadimi2018approximation})
Under the above Assumptions (\ref{ass:1}, \ref{ass:2}, \ref{ass:4}), for all $x,x_1,x_2\in \mathcal{X}$ and $y\in \mathcal{Y}$, we have
\begin{align}
 & \|\bar{\nabla}f(x,y)-\nabla F(x)\| \leq L_y\|y^*(x)-y\|, \nonumber \\
 & \|y^*(x_1)-y^*(x_2)\| \leq \kappa\|x_1-x_2\|, \nonumber \\
 & \|\nabla F(x_1) - \nabla F(x_2)\|\leq L\|x_1-x_2\|,  \nonumber
\end{align}
where $L_y=L_f + L_fC_{gxy}/\mu + C_{fy}\big(L_{gxy}/\mu + L_{gyy}C_{gxy}/\mu^2\big)$, $\kappa=C_{gxy}/\mu$, and
$L=L_f + (L_f+L_y)C_{gxy}/\mu + C_{fy}\big(L_{gxy}/\mu + L_{gyy}C_{gxy}/\mu^2\big)$.
\end{lemma}

For the stochastic bilevel optimization, \cite{yang2021provably,hong2020two} provided a stochastic estimator $\nabla F(x)$ as follows:
\begin{align} \label{eq:5}
 \hat{\nabla} f(x,y;S) & = \nabla_xf(x,y;\xi) - \nabla^2_{xy}g(x,y;\zeta)\vartheta\sum_{q=-1}^{Q-1}\prod_{i=Q-q}^Q \nonumber \\
 & \quad \cdot \Big( I_p - \vartheta \nabla^2_{yy}g(x,y;\zeta^i)\Big)\nabla_yf(x,y;\xi),
\end{align}
where $\vartheta>0$ and $Q\geq 1$.
Here $S=\big\{\xi,\zeta,\zeta^1,\cdots\zeta^Q\big\}$, where $\xi$ is drawn from distribution $\mathcal{D}$,
and $\{\zeta,\zeta^1,\cdots\zeta^Q\}$ are drawn from distribution $\mathcal{H}$.

\section{ Adaptive Bilevel Optimization Methods }
In this section, we propose a class of fast single-loop adaptive bilevel optimization methods to solve the problem \eqref{eq:1}.
Specifically, our methods adopt the universal adaptive learning rates as in \cite{huang2021super}.
Moreover, our methods can be flexibly incorporate the momentum and variance reduced techniques.

\subsection{ BiAdam Algorithm }
In this subsection, we propose a fast single-loop adaptive bilevel optimization method (BiAdam) based on
the basic momentum technique. Algorithm \ref{alg:1} shows
the algorithmic framework of our BiAdam algorithm.

\begin{algorithm}[tb]
\caption{ BiAdam Algorithm }
\label{alg:1}
\begin{algorithmic}[1] 
\STATE {\bfseries Input:} $T, K \in \mathbb{N}^+$, parameters $\{\gamma, \lambda, \eta_t, \alpha_t, \beta_t\}$
and initial input $x_1 \in \mathcal{X}$ and $y_1 \in \mathcal{Y}$; \\
\STATE {\bfseries initialize:} Draw $K+2$ independent samples $\bar{\xi}_1=\{\xi_1,\zeta^0_1,\zeta^1_1, \cdots, \zeta^{K-1}_1\}$ and $\zeta_1$,
and then compute $v_1 = \nabla_y g(x_1,y_1;\zeta_1)$, and $w_1 = \bar{\nabla} f(x_1,y_1;\bar{\xi}_1)$ generated from (\ref{eq:13});  \\
\FOR{$t = 1, 2, \ldots, T$}
\STATE Generate adaptive matrices $A_t \in \mathbb{R}^{d\times d}$, $B_t \in \mathbb{R}^{p\times p}$;
\textcolor{blue}{One example of $A_t$ and $B_t$ by using update rule ($a_0 = 0$, $b_0 = 0$, $ 0 < \tau < 1$, $\rho>0$.) } \\
\textcolor{blue}{  $ a_t = \tau a_{t-1} + (1 - \tau)\big(\nabla_x f(x_t,y_t;\xi_t)\big)^2$, $A_t = \mbox{diag}(\sqrt{a_t} + \rho)$}; \\
\textcolor{blue}{ $ b_t = \tau b_{t-1} + (1 - \tau)||\nabla_y g(x_t,y_t;\zeta_t)||$, $B_t = (b_t + \rho)I_p$}; \\
\STATE $\tilde{x}_{t+1} = \arg\min_{x\in \mathcal{X}}\big\{ \langle w_t, x\rangle + \frac{1}{2\gamma}(x-x_t)^TA_t(x-x_t)\big\}$, and
 $x_{t+1} = x_t+\eta_t(\tilde{x}_{t+1}-x_t)$;
\STATE $\tilde{y}_{t+1} = \arg\min_{y\in \mathcal{Y}}\big\{ \langle v_t, y\rangle + \frac{1}{2\lambda}(y-y_t)^TB_t(y-y_t) \big\}$, and
 $y_{t+1} = y_t+\eta_t(\tilde{y}_{t+1}-y_t)$;
\STATE Draw $K+2$ independent samples $\bar{\xi}_{t+1}=\{\xi_{t+1},\zeta^0_{t+1}, \cdots, \zeta^{K-1}_{t+1}\}$ and $\zeta_{t+1}$:
\STATE $v_{t+1} = \alpha_{t+1}\nabla_y g(x_{t+1},y_{t+1};\zeta_{t+1}) + (1-\alpha_{t+1})v_t $;
\STATE $w_{t+1} = \beta_{t+1}\bar{\nabla} f(x_{t+1},y_{t+1};\bar{\xi}_{t+1}) + (1-\beta_{t+1})w_t $;
\ENDFOR
\STATE {\bfseries Output:} Chosen uniformly random from $\{x_t\}_{t=1}^{T}$.
\end{algorithmic}
\end{algorithm}

At the line 4 of Algorithm~\ref{alg:1}, we generate the adaptive matrices $A_t$ and $B_t$
for updating variables $x$ and $y$, respectively.
In Algorithm~\ref{alg:1}, we give an example, where
 the matrix $A_t$ be generated as the Adam \cite{kingma2014adam}, and
 the matrix $B_t$ as a new version of AdaGrad-Norm \cite{ward2019adagrad}. In fact, we can also defined some other adaptive matrices $A_t$ and $B_t$, as Adabelief~\cite{zhuang2020adabelief} :
\begin{align}
  & a_t = \tau a_{t-1} + (1-\tau) (\nabla_x f(x_t,y_t;\xi_t)-w_t)^2, \ a_0=0, \nonumber \\
  & A_t=\mbox{diag}\big( \sqrt{a_t} + \rho \big), \ t\geq 1   \label{eq:10}  \\
  & b_t = \tau b_{t-1} + (1-\tau)\|\nabla_y g(x_t,y_t;\zeta_t)-v_t\|, \ b_0>0, \nonumber \\
  & B_t=(b_t + \rho )I_p, \  t\geq 1,  \label{eq:11}
\end{align}
where $\tau \in (0,1)$ and $\rho >0$.

At the lines 5-6 of Algorithm \ref{alg:1},
we use the generalized projection gradient iteration with Bregman distance \cite{censor1992proximal,huang2021super} to update the
variables $x$ and $y$, respectively. When $\mathcal{X}=\mathbb{R}^d$ and $\mathcal{Y}=\mathbb{R}^p$, i.e., unconstrained optimization problem (\ref{eq:1}),
we have
$x_{t+1}=x_t - \gamma \eta_t A_t^{-1}w_t$ and $y_{t+1}=y_t - \lambda \eta_t B_t^{-1}v_t$.

At the line 7 of Algorithm \ref{alg:1},
we draw $K+1$ independent samples $\bar{\xi} = \{\xi,\zeta^0,\zeta^1,\cdots, \zeta^{K-1}\}$ from distributions $\mathcal{D}$
and $\mathcal{H}$, then
we define a stochastic gradient estimator as in \cite{khanduri2021near}:
\begin{align} \label{eq:13}
 \bar{\nabla} f(x,y,\bar{\xi}) & = \nabla_xf(x,y;\xi) - \nabla^2_{xy}g(x,y;\zeta^0)  \\
 &\quad \cdot \bigg[ \frac{K}{L_g}\prod_{i=1}^k \big( I_p - \frac{1}{L_g}\nabla^2_{yy}g(x,y;\zeta^i)\big) \bigg] \nabla_yf(x,y;\xi), \nonumber
\end{align}
where $K\geq 1$ and $k\sim \mathcal{U}\{0,1, \cdots, K-1\}$ is a uniform random variable independent on $\bar{\xi}$. In fact, the estimator (\ref{eq:13})
is a specific case of the above estimator (\ref{eq:5}). In practice, thus, we can use a tuning parameter $\vartheta\in (0,\frac{1}{L_g}]$ instead of
$\frac{1}{L_g}$ in the estimator (\ref{eq:13}) as in \cite{yang2021provably}.
Here we use the term $\frac{K}{L_g}\prod_{i=1}^k \big( I_p - \frac{1}{L_g}\nabla^2_{yy}g(x,y;\zeta^i)\big)$ to
approximate the Hessian inverse, i.e., $\big(\nabla^2_{yy} g(x,y;\zeta)\big)^{-1}$.
Clearly, the above $\bar{\nabla} f(x,y,\bar{\xi})$ is a biased estimator in estimating $\bar{\nabla} f(x,y)$, i.e.
$\mathbb{E}_{\bar{\xi}}\big[\bar{\nabla}f(x,y;\bar{\xi})\big] \neq \bar{\nabla}f(x,y)$.
In the following, we give
Lemma \ref{lem:2},
which shows that the bias $R(x,y)=\bar{\nabla}f(x,y) - \mathbb{E}_{\bar{\xi}}\big[\bar{\nabla} f(x,y;\bar{\xi})\big]$
in the gradient estimator (\ref{eq:13}) decays exponentially fast with number $K$.

\begin{lemma} \label{lem:2}
(Lemma 2.1 in \cite{khanduri2021near} and Lemma 11 in \cite{hong2020two})
 Under the about Assumptions (\ref{ass:1}, \ref{ass:4}), for any $K\geq 1$, the gradient estimator in \eqref{eq:13} satisfies
 \begin{align}
  \|R(x,y)\| \leq \frac{C_{gxy}C_{fy}}{\mu}\big(1 - \frac{\mu}{L_g}\big)^K,
 \end{align}
 where $R(x,y)=\bar{\nabla}f(x,y) - \mathbb{E}_{\bar{\xi}}\big[\bar{\nabla} f(x,y;\bar{\xi})\big]$.
\end{lemma}
From Lemma \ref{lem:2}, choose $K=\frac{L_g}{\mu}\log(C_{gxy}C_{fy}T/\mu)$ in Algorithm \ref{alg:1}, we have
 $\|R(x,y)\|\leq \frac{1}{T}$ for all $t\geq 1$. Thus, this result guarantees convergence of our algorithms only requiring a small mini-batch samples.
For notational simplicity, let $R_t=R(x_t,y_t)$ for all $t\geq 1$.

\begin{lemma} \label{lem:3}
(Lemma 3.1 in \cite{khanduri2021near})
Under the above Assumptions (\ref{ass:1}, \ref{ass:3}, \ref{ass:4}), stochastic gradient estimate $\bar{\nabla}f(x,y;\bar{\xi})$ is Lipschitz continuous,
 such that for $x,x_1,x_2\in \mathcal{X}$ and $y,y_1,y_2\in \mathcal{Y}$,
 \begin{align}
   & \mathbb{E}_{\bar{\xi}}\|\bar{\nabla} f(x_1,y;\bar{\xi}) - \bar{\nabla} f(x_2,y;\bar{\xi})\|^2 \leq L^2_K\|x_1-x_2\|^2, \nonumber \\
   & \mathbb{E}_{\bar{\xi}}\|\bar{\nabla} f(x,y_1;\bar{\xi}) - \bar{\nabla} f(x,y_2;\bar{\xi})\|^2 \leq L^2_K\|y_1-y_2\|^2, \nonumber
 \end{align}
 where $L_K^2 = 2L^2_f + 6C^2_{gxy}L^2_f\frac{K}{2\mu L_g - \mu^2} + 6C^2_{fy}L^2_{gxy}\frac{K}{2\mu L_g - \mu^2} + 6C^2_{gxy}L^2_f\frac{K^3L^2_{gyy}}{(L_g-\mu)^2(2\mu L_g - \mu^2)}$.
\end{lemma}

\subsection{ VR-BiAdam Algorithm }
In this subsection, we propose an accelerated version of BiAdam method (VR-BiAdam) by using
the momentum-based variance reduced technique. Algorithm \ref{alg:2} demonstrates
the algorithmic framework of our VR-BiAdam algorithm.

\begin{algorithm}[tb]
\caption{ VR-BiAdam Algorithm }
\label{alg:2}
\begin{algorithmic}[1] 
\STATE {\bfseries Input:} $T, K \in \mathbb{N}^+$, parameters $\{\gamma, \lambda, \eta_t, \alpha_t, \beta_t\}$
and initial input $x_1 \in \mathcal{X}$ and $y_1 \in \mathcal{Y}$; \\
\STATE {\bfseries initialize:} Draw $K+2$ independent samples $\bar{\xi}_1=\{\xi_1,\zeta^0_1,\zeta^1_1, \cdots, \zeta^{K-1}_1\}$ and $\zeta_1$,
and then compute $v_1 = \nabla_y g(x_1,y_1;\zeta_1)$, and $w_1 = \bar{\nabla}f(x_1,y_1;\bar{\xi}_1)$ generated from (\ref{eq:13});  \\
\FOR{$t = 1, 2, \ldots, T$}
\STATE Generate adaptive matrices $A_t \in \mathbb{R}^{d\times d}$, $B_t \in \mathbb{R}^{p\times p}$;
\textcolor{blue}{One example of $A_t$ and $B_t$ by using update rule ($a_0 = 0$, $b_0 = 0$, $ 0 < \tau < 1$, $\rho>0$.) } \\
\textcolor{blue}{  $ a_t = \tau a_{t-1} + (1 - \tau)\big(\nabla_x f(x_t,y_t;\xi_t)\big)^2$, $A_t = \mbox{diag}(\sqrt{a_t} + \rho)$}; \\
\textcolor{blue}{ $ b_t = \tau b_{t-1} + (1 - \tau)||\nabla_y g(x_t,y_t;\zeta_t)||$, $B_t = (b_t + \rho)I_p$}; \\
\STATE $\tilde{x}_{t+1} = \arg\min_{x\in \mathcal{X}}\big\{ \langle w_t, x\rangle + \frac{1}{2\gamma}(x-x_t)^TA_t(x-x_t)\big\}$, and
 $x_{t+1} = x_t+\eta_t(\tilde{x}_{t+1}-x_t)$;
\STATE $\tilde{y}_{t+1} = \arg\min_{y\in \mathcal{Y}}\big\{ \langle v_t, y\rangle + \frac{1}{2\lambda}(y-y_t)^TB_t(y-y_t) \big\}$, and
 $y_{t+1} = y_t+\eta_t(\tilde{y}_{t+1}-y_t)$;
\STATE Draw $K+2$ independent samples $\bar{\xi}_{t+1}=\{\xi_{t+1},\zeta^0_{t+1}, \cdots, \zeta^{K-1}_{t+1}\}$ and $\zeta_{t+1}$;
\STATE $v_{t+1} = \nabla_y g(x_{t+1},y_{t+1};\zeta_{t+1}) + (1-\alpha_{t+1})\big[v_t - \nabla_y g(x_t,y_t;\zeta_{t+1})\big]$;
\STATE $w_{t+1} = \bar{\nabla} f(x_{t+1},y_{t+1};\bar{\xi}_{t+1}) + (1-\beta_{t+1})\big[w_t - \bar{\nabla} f(x_t,y_t;\bar{\xi}_{t+1})\big] $;
\ENDFOR
\STATE {\bfseries Output:} Chosen uniformly random from $\{x_t\}_{t=1}^{T}$.
\end{algorithmic}
\end{algorithm}

At the lines 8-9 of Algorithm \ref{alg:2}, we use the momentum-based variance reduced technique to
estimate the stochastic partial derivatives $v_t$ and $w_t$.
For example, the estimator of partial derivative $\bar{\nabla }f(x_{t+1},y_{t+1})$ is defined as
\begin{align}
 w_{t+1}
  &= \beta_{t+1}\bar{\nabla} f(x_{t+1},y_{t+1};\bar{\xi}_{t+1}) + (1-\beta_{t+1})  \big[w_t \nonumber \\
  & \quad + \bar{\nabla} f(x_{t+1},y_{t+1};\bar{\xi}_{t+1}) - \bar{\nabla} f(x_t,y_t;\bar{\xi}_{t+1})\big]. \nonumber
\end{align}
Compared with the estimator $w_{t+1}$ in Algorithm \ref{alg:1}, $w_{t+1}$ in Algorithm \ref{alg:2}
adds the term $\bar{\nabla} f(x_{t+1},y_{t+1};\bar{\xi}_{t+1}) - \bar{\nabla} f(x_t,y_t;\bar{\xi}_{t+1})$ to control the variances of estimator.

In fact, our algorithms use the the momentum techniques to update variables and estimate stochastic gradients simultaneously. For example, our BiAdam uses the momentum iteration (i.e., $x_{t+1}=x_t + \eta_t(\tilde{x}_{t+1}-x_t)$ at the line 5 of our Algorithm 1) to update variable $x$, and applies the basic momentum technique to estimate the stochastic gradient $w_t$
 (i.e., the line 9 of our Algorithm 1). Our VR-BiAdam also uses the momentum iteration (i.e.,
$x_{t+1}=x_t + \eta_t(\tilde{x}_{t+1}-x_t)$ at the line 5 of our Algorithm 2) to update variable
$x$, and applies the variance-reduced momentum technique of STROM to estimate the stochastic gradient $w_t$ (i.e., the line 9 of our Algorithm 2).

\begin{figure*}[!t]
  \centering
  \includegraphics[width=0.92\textwidth]{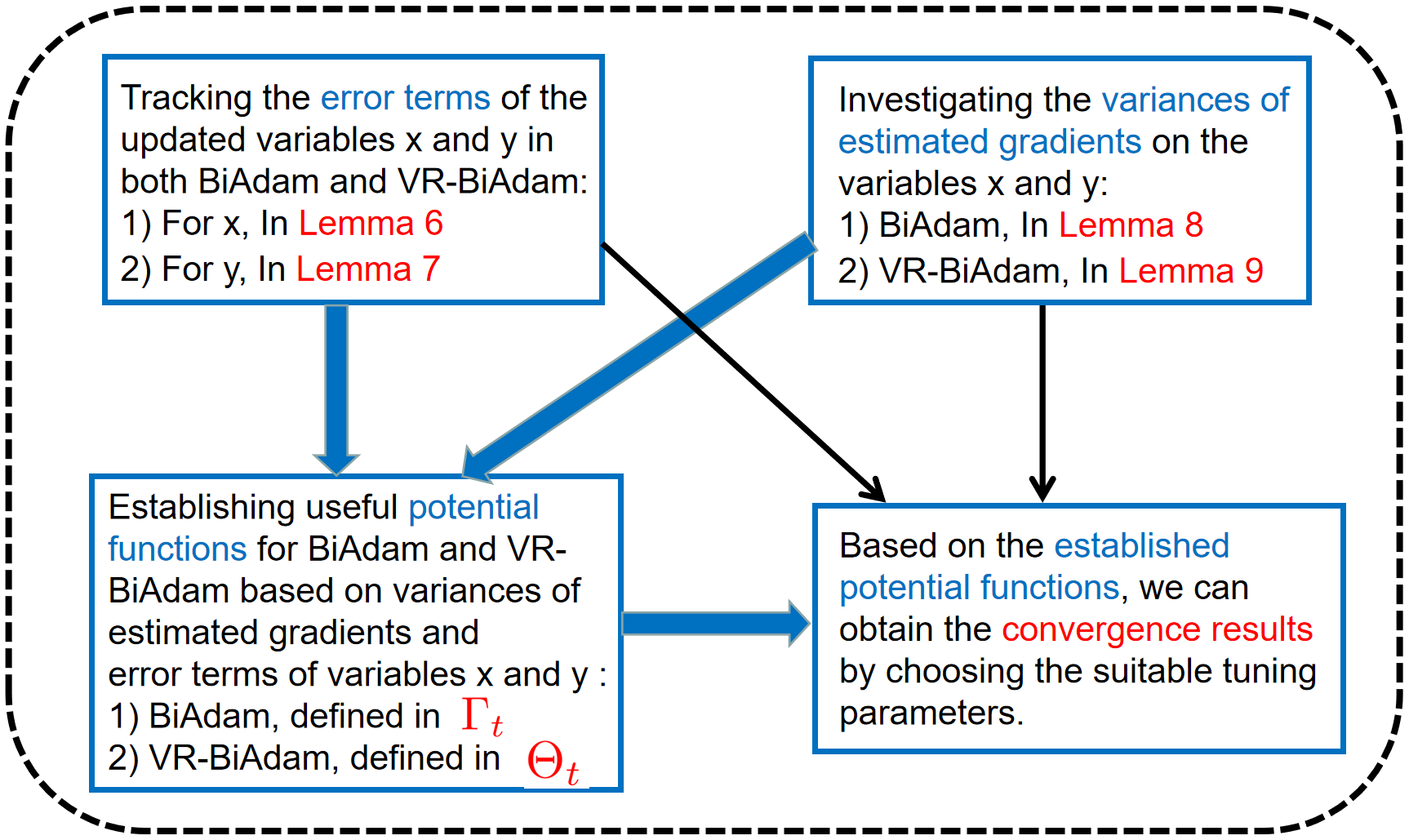}
  \caption{ The basic idea of our convergence analysis. }
  \label{fig:1}
\end{figure*}

\section{ Theoretical Analysis }
In this section, we study the convergence properties of our algorithms (BiAdam and VR-BiAdam) under some mild conditions. The basic idea of our convergence analysis is given in Fig.~\ref{fig:1}.
All proofs are provided in the Appendix A.

\subsection{ Additional Mild Assumptions }

\begin{assumption} \label{ass:6}
 The estimated stochastic partial derivative $\bar{\nabla} f(x,y;\bar{\xi})$ satisfies
 \begin{align}
 & \mathbb{E}_{\bar{\xi}}\big[\bar{\nabla} f(x,y;\bar{\xi})\big] = \bar{\nabla} f(x,y) + R(x,y), \nonumber \\
 & \mathbb{E}_{\bar{\xi}}\|\bar{\nabla}f(x,y;\bar{\xi}) - \bar{\nabla} f(x,y) - R(x,y)\|^2 \leq \sigma^2. \nonumber
 \end{align}
 The stochastic partial derivative $\nabla_y g(x,y;\zeta)$ satisfies
 \begin{align}
  & \mathbb{E}[\nabla_y g(x,y;\zeta)] = \nabla_y g(x,y), \nonumber \\
  & \mathbb{E}\|\nabla_y g(x,y;\zeta)-\nabla_y g(x,y)\|^2 \leq \sigma^2. \nonumber
 \end{align}
\end{assumption}
\begin{assumption} \label{ass:7}
In our algorithms, the adaptive matrices $A_t$ and $B_t$ for all $t\geq 1$ satisfy
$A_t \succeq \rho I_d \ (\rho>0)$ and $B_t = bI_p \ (b_u \geq b \geq b_l>0)$, respectively, where $\rho$, $b_u$ and $b_l$ are appropriate positive numbers.
\end{assumption}

Assumption \ref{ass:6} is commonly used in the stochastic bilevel optimization methods \cite{ji2021bilevel,yang2021provably,khanduri2021near}.
In the paper, we consider the general adaptive learning rates (including the coordinate-wise and global learning rates)
for variable $x$ and the global learning rate for variable $y$. Assumption \ref{ass:7} ensures that the adaptive matrices $A_t$ for all $t\geq 1$ are positive definite as in \cite{huang2021super}.
Assumption \ref{ass:7} also guarantees the global  adaptive matrices $B_t$ for all $t\geq 1$ are positive definite and bounded.
In fact, Assumption \ref{ass:7} is mild.
For example, in the problem $\min_{x\in \mathbb{R}^p} \mathbb{E}[f(x;\xi)]$, \cite{ward2019adagrad} apply a global adaptive learning rate
to the update form $x_t = x_{t-1} - \eta\frac{\nabla f(x_{t-1};\xi_{t-1})}{b_t}, \ b_t^2 = b_{t-1}^2 + \|\nabla f(x_{t-1};\xi_{t-1})\|^2,  \ b_0>0, \eta>0$ for all $t\geq 1$,
which is equivalent to the form $ x_t = x_{t-1} - \eta B_t^{-1}\nabla f(x_{t-1};\xi_{t-1})$ with $B_t = b_tI_p$ and
$b_t \geq  \cdots \geq b_0 >0$.
\cite{li2019convergence,cutkosky2019momentum} use a global adaptive learning rate
to the update form $x_{t+1} = x_t - \eta g_t/b_t$, where $g_t$ is stochastic gradient and
$b_t=\big(\omega+\sum_{i=1}^t\|\nabla f(x_i;\xi_i)\|^2\big)^{\alpha}/k$, $k>0$, $\omega>0$ and $\alpha\in (0,1)$,
which is equivalent to $x_{t+1} = x_t - \eta B_t^{-1}g_t$ with $B_t = b_tI_p$ and $b_t \geq \cdots \geq b_0=\frac{\omega^\alpha}{k}>0$.
At the same time, the problem $\min_{x\in \mathbb{R}^p}f(x)=\mathbb{E}[f(x;\xi)]$ approaches the stationary points, i.e.,
$\nabla f(x)=0$ or even $\nabla f(x;\xi)=0$ for all $\xi$. Thus, these global adaptive learning rates are generally bounded, i.e.,
$b_u \geq b_t \geq b_l>0$ for all $t\geq 1$.

\subsection{ Useful Convergence Metric and Lemmas }
In the subsection, we first define a useful convergence metric for our algorithms.
\begin{lemma} \label{lem:4}
Given gradient estimator $w_t$ is generated from Algorithms \ref{alg:1} or \ref{alg:2}, for all $t\geq 1$,
we have
 \begin{align}
 \|w_t-\nabla F(x_t)\|^2 \leq L^2_0\|y^*(x_t)-y_t\|^2 + 2\|w_t-\bar{\nabla} f(x_t,y_t)\|^2, \nonumber
\end{align}
where $L^2_0 = 8\big(L^2_f+ \frac{L^2_{gxy}C^2_{fy}}{\mu^2} + \frac{L^2_{gyy} C^2_{gxy}C^2_{fy}}{\mu^4} +
 \frac{L^2_fC^2_{gxy}}{\mu^2}\big)$.
\end{lemma}

For our Algorithms \ref{alg:1} and \ref{alg:2}, based on Lemma \ref{lem:4},
we provide a convergence metric $\mathbb{E}[\mathcal{M}_t]$,
defined as
\begin{align}
 \mathcal{M}_t & =  \frac{1}{\gamma}\|\tilde{x}_{t+1} - x_t\| + \frac{1}{\rho}\Big( \sqrt{2}\|w_t - \bar{\nabla}f(x_t,y_t)\|
 \nonumber \\
 & \quad + L_0\|y^*(x_t)-y_t\| \Big), \nonumber
\end{align}
where the first two terms of $\mathcal{M}_t$ measure the convergence of the iteration solutions $\{x_t\}_{t=1}^T$,
and the last term measures the convergence of the iteration solutions $\{y_t\}_{t=1}^T$.

Let $\phi_t(x)=\frac{1}{2}x^TA_t x$, we define a prox-function (a.k.a., Bregman distance) \cite{censor1981iterative,censor1992proximal,ghadimi2016mini} associated with $\phi_t(x)$ as follows:
\begin{align}
 V_t(x,x_t) &= \phi_t(x) - \big[ \phi_t(x_t) + \langle\nabla \phi_t(x_t), x-x_t\rangle\big] \nonumber \\
 & = \frac{1}{2}(x-x_t)^TA_t(x-x_t).
\end{align}
The line 5 of Algorithm \ref{alg:1} or \ref{alg:2} is equivalent to the following generalized projection problem:
\begin{align} \label{eq:18}
 \tilde{x}_{t+1} = \min_{x\in \mathcal{X}}\Big\{ \langle w_t, x\rangle + \frac{1}{\gamma}V_t(x,x_t) \Big\}.
\end{align}
As in \cite{ghadimi2016mini}, we define a generalized projected gradient $\mathcal{G}_{\mathcal{X}}(x_t,w_t,\gamma)=\frac{1}{\gamma}(x_t-\tilde{x}_{t+1})$. At the same time, we define a gradient mapping $\mathcal{G}_{\mathcal{X}}(x_t,\nabla F(x_t),\gamma)=\frac{1}{\gamma}(x_t-x^+_{t+1})$ with
\begin{align}
   x^+_{t+1} =  \arg\min_{x\in \mathcal{X}}\Big\{ \langle \nabla F(x_t), x\rangle + \frac{1}{\gamma}V_t(x,x_t) \Big\}.
\end{align}
According to the Proposition 1 of \cite{ghadimi2016mini}, we have $||\mathcal{G}_\mathcal{X}(x_t,w_t,\gamma)-\mathcal{G}_\mathcal{X} (x_t,\nabla F(x_t),\gamma)|| \leq \frac{1}{\rho} || \nabla F(x_t)-w_t ||$.
 Since $||\mathcal{G}_\mathcal{X}(x_t,\nabla F(x_t),\gamma)||  \leq  ||\mathcal{G}_\mathcal{X}(x_t,w_t,\gamma)|| + ||\mathcal{G}_\mathcal{X}(x_t,w_t,\gamma)-\mathcal{G}_\mathcal{X}(x_t,\nabla F(x_t),\gamma)||$,
 we have
\begin{align} \label{eq:19}
 &||\mathcal{G}_\mathcal{X}(x_t,\nabla F(x_t),\gamma)||
  \leq ||\mathcal{G}_\mathcal{X}(x_t,w_t,\gamma)|| + \frac{1}{\rho}||\nabla F(x_t)-w_t||  \nonumber \\
 & \leq \frac{1}{\gamma}\|\tilde{x}_{t+1} - x_t\| + \frac{1}{\rho}\Big(\sqrt{2}\|w_t - \bar{\nabla}f(x_t,y_t)\|    \nonumber \\
 &\quad + L_0\|y^*(x_t)-y_t\|\Big)=\mathcal{M}_t,
\end{align}
where the last inequality holds by the above Lemma \ref{lem:4}.
Thus, our new convergence measure $\mathbb{E} [\mathcal{M}_t]$ is tighter than the standard gradient mapping $\mathbb{E}||\mathcal{G}_\mathcal{X}(x_t,\nabla F(x_t),\gamma)||$
used in \cite{hong2020two}.
When $\mathcal{M}_t \rightarrow 0$, we have  $\|\mathcal{G}_{\mathcal{X}}(x_t,\nabla F(x_t),\gamma)\| \rightarrow 0$,
where $x_t$ is a stationary point or local minimum of the bilevel problem \eqref{eq:1} \cite{ghadimi2016mini,hong2020two}.

Next, we provide some useful lemmas.
\begin{lemma} \label{lem:6}
 Suppose that the sequence $\{x_t,y_t\}_{t=1}^T$ be generated from Algorithm \ref{alg:1} or \ref{alg:2}.
 Let $0<\eta_t \leq 1$ and $0< \gamma \leq \frac{\rho}{2L\eta_t}$,
 then we have
 \begin{align}
  F(x_{t+1}) \leq F(x_t) + \frac{\eta_t\gamma}{\rho}\|\nabla F(x_t)-w_t\|^2 -\frac{\rho\eta_t}{2\gamma}\|\tilde{x}_{t+1}-x_t\|^2. \nonumber
 \end{align}
\end{lemma}

\begin{lemma} \label{lem:7}
Suppose the sequence $\{x_t,y_t\}_{t=1}^T$ be generated from Algorithm \ref{alg:1} or \ref{alg:2}.
Under the above assumptions, given $0< \eta_t\leq 1$, $B_t=b_tI_p \ (b_u \geq b_t \geq b_l>0)$ for all $t\geq 1$,
and $0<\lambda\leq \frac{b_l}{6L_g}$, we have
\begin{align}
  & \|y_{t+1} - y^*(x_{t+1})\|^2 - \|y_t -y^*(x_t)\|^2 \nonumber \\
  & \leq -\frac{\eta_t\mu\lambda}{4b_t}\|y_t -y^*(x_t)\|^2 -\frac{3\eta_t}{4} \|\tilde{y}_{t+1}-y_t\|^2 \nonumber \\
  & \quad + \frac{25\eta_t\lambda}{6\mu b_t} \|\nabla_y g(x_t,y_t)-v_t\|^2 + \frac{25\kappa^2\eta_tb_t}{6\mu\lambda}\|\tilde{x}_{t+1} - x_t\|^2, \nonumber
\end{align}
where $\kappa = L_g/\mu$.
\end{lemma}

\subsection{ Convergence Analysis of BiAdam Algorithm }
In this subsection, we study the convergence properties of our BiAdam algorithm.
The detailed proofs are provided in the Appendix \ref{Appendix:A1}.

\begin{lemma} \label{lem:8}
 Assume that the stochastic partial derivatives $v_{t+1}$, and $w_{t+1}$ be generated from Algorithm \ref{alg:1}, we have
 \begin{align}
 & \mathbb{E}\| w_{t+1} - \bar{\nabla} f(x_{t+1},y_{t+1})-R_{t+1}\|^2  \\
 & \leq  (1-\beta_{t+1}) \mathbb{E}\|w_t - \bar{\nabla} f(x_t,y_t) -R_t\|^2  + \beta^2_{t+1}\sigma^2  \nonumber \\
 & \quad + \frac{3L^2_0\eta^2_t}{\beta_{t+1}}\big( \|\tilde{x}_{t+1}-x_t\|^2 + \|\tilde{y}_{t+1}-y_t\|^2 \big) \nonumber \\
 & \quad + \frac{3}{\beta_{t+1}}\big( \| R_t\|^2 + \|R_{t+1}\|^2 \big), \nonumber
 \end{align}
 \begin{align}
 & \mathbb{E}\|v_{t+1}- \nabla_y g(x_{t+1},y_{t+1})\|^2  \\
 & \leq (1-\alpha_{t+1}) \mathbb{E} \|v_t - \nabla_y g(x_t,y_t)\|^2 + \alpha_{t+1}^2\sigma^2 \nonumber \\
 & \quad + 2L_g^2\eta_t^2/\alpha_{t+1}\big(\mathbb{E}\|\tilde{x}_{t+1} - x_t\|^2 + \mathbb{E}\|\tilde{y}_{t+1} - y_t\|^2 \big), \nonumber
 \end{align}
where $R_t = \bar{\nabla}f(x_t,y_t) - \mathbb{E}_{\bar{\xi}}[\bar{\nabla}f(x_t,y_t;\bar{\xi})]$
 for all $t\geq 1$.
\end{lemma}

In the convergence analysis of BiAdam, we define a useful \emph{Lyapunov} function (i.e., potential function) $\Gamma_t$, for any $t\geq 1$,
 \begin{align}
 \Gamma_t & = \mathbb{E}\big [F(x_t) + \frac{5b_tL^2_0\gamma}{\lambda\mu\rho}\|y_t-y^*(x_t)\|^2   \nonumber \\
 & \quad + \frac{\gamma}{\rho} \big( \|v_t - \nabla_y g(x_t,y_t)\|^2+ \|w_t - \bar{\nabla} f(x_t,y_t) -R_t\|^2 \big) \big]. \nonumber
 \end{align}

\begin{theorem} \label{th:1}
 Under the above Assumptions (\ref{ass:1}, \ref{ass:2}, \ref{ass:4}, \ref{ass:6}, \ref{ass:7}), in the Algorithm \ref{alg:1}, given $\mathcal{X}\subset\mathbb{R}^{d}$, $\eta_t=\frac{k}{(m+t)^{1/2}}$ for all $t\geq 0$, $\alpha_{t+1}=c_1\eta_t$, $\beta_{t+1}=c_2\eta_t$, $m\geq \max\big(k^2, (c_1k)^2,(c_2k)^2\big)$, $k>0$, $\frac{125L^2_0}{6\mu^2} \leq c_1 \leq \frac{m^{1/2}}{k}$, $\frac{9}{2} \leq c_2 \leq \frac{m^{1/2}}{k}$, $0<\lambda\leq \min\big(\frac{15b_lL^2_0}{4L^2_1\mu}, \frac{b_l}{6L_g}\big)$, $0< \gamma \leq \min\big(\frac{\sqrt{6}\lambda\mu\rho}{\sqrt{6L^2_1\lambda^2\mu^2 + 125b^2_uL^2_0\kappa^2}}, \frac{m^{1/2}\rho}{4Lk}\big)$  and $K=\frac{L_g}{\mu}\log(C_{gxy}C_{fy}T/\mu)$, we have
\begin{align}
  & \frac{1}{T} \sum_{t=1}^T \mathbb{E}||\mathcal{G}_\mathcal{X}(x_t,\nabla F(x_t),\gamma)||  \leq \frac{1}{T} \sum_{t=1}^T \mathbb{E}[\mathcal{M}_t] \nonumber \\
  & \leq \frac{2\sqrt{3G}m^{1/4}}{T^{1/2}} + \frac{2\sqrt{3G}}{T^{1/4}} + \frac{\sqrt{2}}{T},
\end{align}
where $G = \frac{F(x_1) - F^*}{\rho k\gamma} + \frac{5b_1L^2_0\Delta_0}{\rho^2 k\lambda\mu} + \frac{2\sigma^2}{\rho^2 k} + \frac{2m\sigma^2}{\rho^2k}\ln(m+T) + \frac{4(m+T)}{9\rho^2 kT^2} + \frac{8k}{\rho^2T}$ and $\Delta_0=\|y_1-y^*(x_1)\|^2$, $L^2_1= \frac{12L^2_g\mu^2}{125L^2_0} + \frac{2L^2_0}{3}$.
\end{theorem}

\begin{remark}
Without loss of generality, let $k=O(1)$ and $m=O(1)$, we have $G=O(\ln(m+T))=\tilde{O}(1)$.
Thus our BiAdam algorithm has a convergence rate of $\tilde{O}(\frac{1}{T^{1/4}})$. Let $\mathbb{E}[\mathcal{M}_\zeta] = \frac{1}{T} \sum_{t=1}^T \mathbb{E}[\mathcal{M}_t] = \tilde{O}(\frac{1}{T^{1/4}}) \leq \epsilon$, we have $T=\tilde{O}(\epsilon^{-4})$. Since our BiAdam algorithm
only requires $K+2=\frac{L_g}{\mu}\log(C_{gxy}C_{fy}T/\mu)+2=\tilde{O}(1)$ samples to estimate stochastic partial derivatives in each iteration,
and needs $T$ iterations.
Thus our BiAdam algorithm requires sample complexity of $(K+2)T=\tilde{O}(\epsilon^{-4})$ for finding an $\epsilon$-stationary point of Problem (\ref{eq:1}).
\end{remark}

Next, we further give the convergence properties of our BiAdam algorithm for \emph{unconstrained} bilevel optimization.
\begin{theorem} \label{th:2}
 Under the above Assumptions (\ref{ass:1}, \ref{ass:2}, \ref{ass:4}, \ref{ass:6}, \ref{ass:7}), in the Algorithm \ref{alg:1}, given $\mathcal{X}=\mathbb{R}^{d}$, $\eta_t=\frac{k}{(m+t)^{1/2}}$ for all $t\geq 0$, $\alpha_{t+1}=c_1\eta_t$, $\beta_{t+1}=c_2\eta_t$, $m\geq \max\big(k^2, (c_1k)^2,(c_2k)^2\big)$, $k>0$, $\frac{125L^2_0}{6\mu^2} \leq c_1 \leq \frac{m^{1/2}}{k}$, $\frac{9}{2} \leq c_2 \leq \frac{m^{1/2}}{k}$, $0<\lambda\leq \min\big(\frac{15b_lL^2_0}{4L^2_1\mu}, \frac{b_l}{6L_g}\big)$, $0< \gamma \leq \min\big(\frac{\sqrt{6}\lambda\mu\rho}{\sqrt{6L^2_1\lambda^2\mu^2 + 125b^2_uL^2_0\kappa^2}}, \frac{m^{1/2}\rho}{4Lk}\big)$ and $K=\frac{L_g}{\mu}\log(C_{gxy}C_{fy}T/\mu)$, we have
\begin{align} \label{eq:t1}
 \frac{1}{T}\sum_{t=1}^T\mathbb{E}\|\nabla F(x_t)\| & \leq \frac{\sqrt{\frac{1}{T}\sum_{t=1}^T\mathbb{E}\|A_t\|^2}}{\rho}\Big( \frac{2\sqrt{6G'm}}{T^{1/2}} \nonumber \\
 & \quad + \frac{2\sqrt{6G'}}{T^{1/4}} + \frac{2\sqrt{3}}{T}\Big),
\end{align}
where $G' = \frac{\rho(F(x_1) - F^*)}{k\gamma} + \frac{5b_1L^2_0\Delta_0}{k\lambda\mu} + \frac{2\sigma^2}{k} + \frac{2m\sigma^2}{k}\ln(m+T) + \frac{4(m+T)}{9kT^2} + \frac{8k}{T}$.
\end{theorem}

\begin{remark}
Under the same conditions in Theorem \ref{th:1}, based on the metric $\mathbb{E}\|\nabla F(x)\|\leq \epsilon$,
our BiAdam algorithm still has a gradient complexity of $\tilde{O}(\epsilon^{-4})$ without relying on the large mini-batches.
Interestingly, the right hand side of the above inequality \eqref{eq:t1} includes a term $\frac{\sqrt{\frac{1}{T}\sum_{t=1}^T\mathbb{E}\|A_t\|^2}}{\rho}$ that
can be seen as an upper bound of the expected condition number of adaptive matrices $\{A_t\}_{t=1}^T$.
When $A_t$ given in Algorithm~\ref{alg:1}, we have $\frac{\sqrt{\frac{1}{T}\sum_{t=1}^T\mathbb{E}\|A_t\|^2}}{\rho} \leq \frac{G_1+\rho}{\rho}$
as in the existing adaptive gradient methods
assuming the bounded stochastic gradient $\|\nabla_x f(x,y;\xi)\| \leq G_1$. \textbf{NOTE THAT}: since our algorithms use the momentum techniques to update variables and estimate stochastic gradients simultaneously, our algorithms can use a relatively large $\rho$. So the term $\frac{G_1+\rho}{\rho}$ is not very large.
\end{remark}

\subsection{ Convergence Analysis of VR-BiAdam Algorithm }
In this subsection, we study convergence properties of our VR-BiAdam algorithm.
The detailed proofs are provided in the Appendix \ref{Appendix:A2}.

\begin{lemma} \label{lem:9}
Under the above Assumptions (\ref{ass:1}, \ref{ass:3}, \ref{ass:4}), assume the stochastic gradient estimators $v_t$ and $w_t$ be generated from Algorithm \ref{alg:2},
we have
 \begin{align}
& \mathbb{E}\|w_{t+1} - \bar{\nabla} f(x_{t+1},y_{t+1}) - R_{t+1}\|^2  \\
& \leq (1-\beta_{t+1}) \mathbb{E}\|w_t - \bar{\nabla} f(x_t,y_t) -R_t\|^2 + 2\beta^2_{t+1}\sigma^2 \nonumber \\
 &  \quad + 4L^2_K\eta^2_t\big( \|\tilde{x}_{t+1}-x_t\|^2 + \|\tilde{y}_{t+1}-y_t\|^2 \big), \nonumber
\end{align}
 \begin{align}
 & \mathbb{E}\|\nabla_y g(x_{t+1},y_{t+1}) - v_{t+1}\|^2  \\
 & \leq (1-\alpha_{t+1})\mathbb{E} \|\nabla_y g(x_t,y_t) - v_t\|^2 + 2\alpha_{t+1}^2\sigma^2 \nonumber \\
 & \quad + 4L_g^2\eta_t^2\big(\mathbb{E}\|\tilde{x}_{t+1} - x_t\|^2 + \mathbb{E}\|\tilde{y}_{t+1} - y_t\|^2 \big), \nonumber
 \end{align}
where $R_t = \bar{\nabla}f(x_t,y_t) - \mathbb{E}_{\bar{\xi}}[\bar{\nabla}f(x_t,y_t;\bar{\xi})]$
 for all $t\geq 1$.
\end{lemma}

In the convergence analysis of VR-BiAdam, we define a useful \emph{Lyapunov} function $\Theta_t$, for any $t\geq 1$
\begin{align}
 \Theta_t & = \mathbb{E}\big [F(x_t) + \frac{5b_t L^2_0\gamma}{\lambda\mu\rho}\|y_t-y^*(x_t)\|^2  \nonumber \\
 & + \frac{\gamma}{\rho\eta_{t-1}} \big(\|v_t - \nabla_y g(x_t,y_t)\|^2+ \|w_t-\bar{\nabla} f(x_t,y_t)-R_t\|^2 \big) \big]. \nonumber
\end{align}

\begin{theorem}  \label{th:3}
 Under the above Assumptions (\ref{ass:1}, \ref{ass:3}, \ref{ass:4}, \ref{ass:6}, \ref{ass:7}), in the Algorithm \ref{alg:2}, given $\mathcal{X}\subset\mathbb{R}^{d}$, $\eta_t=\frac{k}{(m+t)^{1/3}}$
 for all $t\geq 0$, $\alpha_{t+1}=c_1\eta_t^2$, $\beta_{t+1}=c_2\eta_t^2$, $m\geq \max\big(2,k^3, (c_1k)^3,(c_2k)^3\big)$, $k>0$, $c_1 \geq \frac{2}{3k^3} + \frac{125L^2_0}{6\mu^2}$,
 $c_2 \geq \frac{2}{3k^3} + \frac{9}{2}$, $0<\lambda\leq \min\big( \frac{15b_l L^2_0}{16L^2_2\mu}, \frac{b_l}{6L_g} \big)$,
 $0< \gamma \leq \min\big(\frac{\sqrt{6}\lambda\mu\rho}{2\sqrt{24L^2_2\lambda^2\mu^2 + 125b^2_uL^2_0\kappa^2}}, \frac{m^{1/3}\rho}{4Lk}\big)$ and $K=\frac{L_g}{\mu}\log(C_{gxy}C_{fy}T/\mu)$,
 we have
\begin{align}
  & \frac{1}{T} \sum_{t=1}^T \mathbb{E}||\mathcal{G}_\mathcal{X}(x_t,\nabla F(x_t),\gamma)|| \leq \frac{1}{T} \sum_{t=1}^T \mathbb{E}[\mathcal{M}_t]  \nonumber \\
  & \leq \frac{2\sqrt{3M}m^{1/6}}{T^{1/2}} + \frac{2\sqrt{3M}}{T^{1/3}}
 + \frac{\sqrt{2}}{T},
\end{align}
where $M =\frac{F(x_1) - F^*}{\rho k\gamma} + \frac{5b_1L^2_0\Delta_0}{ \rho^2 k\lambda\mu} + \frac{2m^{1/3}\sigma^2}{\rho^2 k^2} + \frac{2k^2(c^2_1+c^2_2)\sigma^2\ln(m+T)}{\rho^2} + \frac{6k(m+T)^{1/3}}{\rho^2 T}$, $\Delta_0 = \|y_1-y^*(x_1)\|^2$ and $L^2_2 = L^2_g+L^2_K$.
\end{theorem}

\begin{remark}
Without loss of generality, let $k=O(1)$ and $m=O(1)$, we have $M=O(\ln(m+T))=\tilde{O}(1)$.
Thus our VR-BiAdam algorithm has a convergence rate of $\tilde{O}(\frac{1}{T^{1/3}})$.
Let $\mathbb{E}[\mathcal{M}_\zeta] = \frac{1}{T} \sum_{t=1}^T \mathbb{E}[\mathcal{M}_t] = \tilde{O}(\frac{1}{T^{1/3}}) \leq \epsilon$, we have $T=\tilde{O}(\epsilon^{-3})$. Since our VR-BiAdam algorithm
 requires $K+2=\frac{L_g}{\mu}\log(C_{gxy}C_{fy}T/\mu)+2=\tilde{O}(1)$ samples to estimate stochastic partial derivatives in each iteration,
 and needs $T$ iterations.
Thus our VR-BiAdam algorithm requires sample complexity of $(K+2)T=\tilde{O}(\epsilon^{-3})$ for finding an $\epsilon$-stationary point of Problem (\ref{eq:1}).
\end{remark}

Next, we further give the convergence properties of our VR-BiAdam algorithm for \emph{unconstrained} bilevel optimization.
\begin{theorem} \label{th:4}
 Under the above Assumptions (\ref{ass:1}, \ref{ass:3}, \ref{ass:4}, \ref{ass:6}, \ref{ass:7}), in the Algorithm \ref{alg:2}, given $\mathcal{X}=\mathbb{R}^{d}$, $\eta_t=\frac{k}{(m+t)^{1/3}}$ for all $t\geq 0$, $\alpha_{t+1}=c_1\eta_t^2$, $\beta_{t+1}=c_2\eta_t^2$, $m\geq \max\big(2,k^3, (c_1k)^3,(c_2k)^3\big)$, $k>0$, $c_1 \geq \frac{2}{3k^3} + \frac{125L^2_0}{6\mu^2}$, $c_2 \geq \frac{2}{3k^3} + \frac{9}{2}$, $0<\lambda\leq \min\big( \frac{15b_lL^2_0}{16L^2_2\mu}, \frac{b_l}{6L_g} \big)$, $0< \gamma \leq \min\big(\frac{\sqrt{6}\lambda\mu\rho}{2\sqrt{24L^2_2\lambda^2\mu^2 + 125b^2_uL^2_0\kappa^2}}, \frac{m^{1/3}\rho}{4Lk}\big)$ and $K=\frac{L_g}{\mu}\log(C_{gxy}C_{fy}T/\mu)$, we have
\begin{align}
 \frac{1}{T}\sum_{t=1}^T\mathbb{E}\|\nabla F(x_t)\|  & \leq \frac{\sqrt{\frac{1}{T}\sum_{t=1}^T\mathbb{E}\|A_t\|^2}}{\rho}\Big( \frac{2\sqrt{6M'm}}{T^{1/2}}  \nonumber \\
 &\quad + \frac{2\sqrt{6M'}}{T^{1/3}} + \frac{2\sqrt{3}}{T}\Big),
\end{align}
where $M' = \frac{\rho(F(x_1) - F^*)}{k\gamma} + \frac{5b_1L^2_0\Delta_0}{k\lambda\mu} + \frac{2m^{1/3}\sigma^2}{k^2} + 2k^2(c^2_1+c^2_2)\sigma^2\ln(m+T) + \frac{6k(m+T)^{1/3}}{T}$.
\end{theorem}

\begin{remark}
Under the same conditions in Theorem \ref{th:3}, based on the metric $\mathbb{E}\|\nabla F(x)\|\leq \epsilon$,
our VR-BiAdam algorithm still has a gradient complexity of $\tilde{O}(\epsilon^{-3})$ in finding an $\epsilon$-stationary solution without relying on the large mini-batches.
\end{remark}

\begin{figure*}[ht]
	\begin{center}
		\includegraphics[width=0.63\columnwidth]{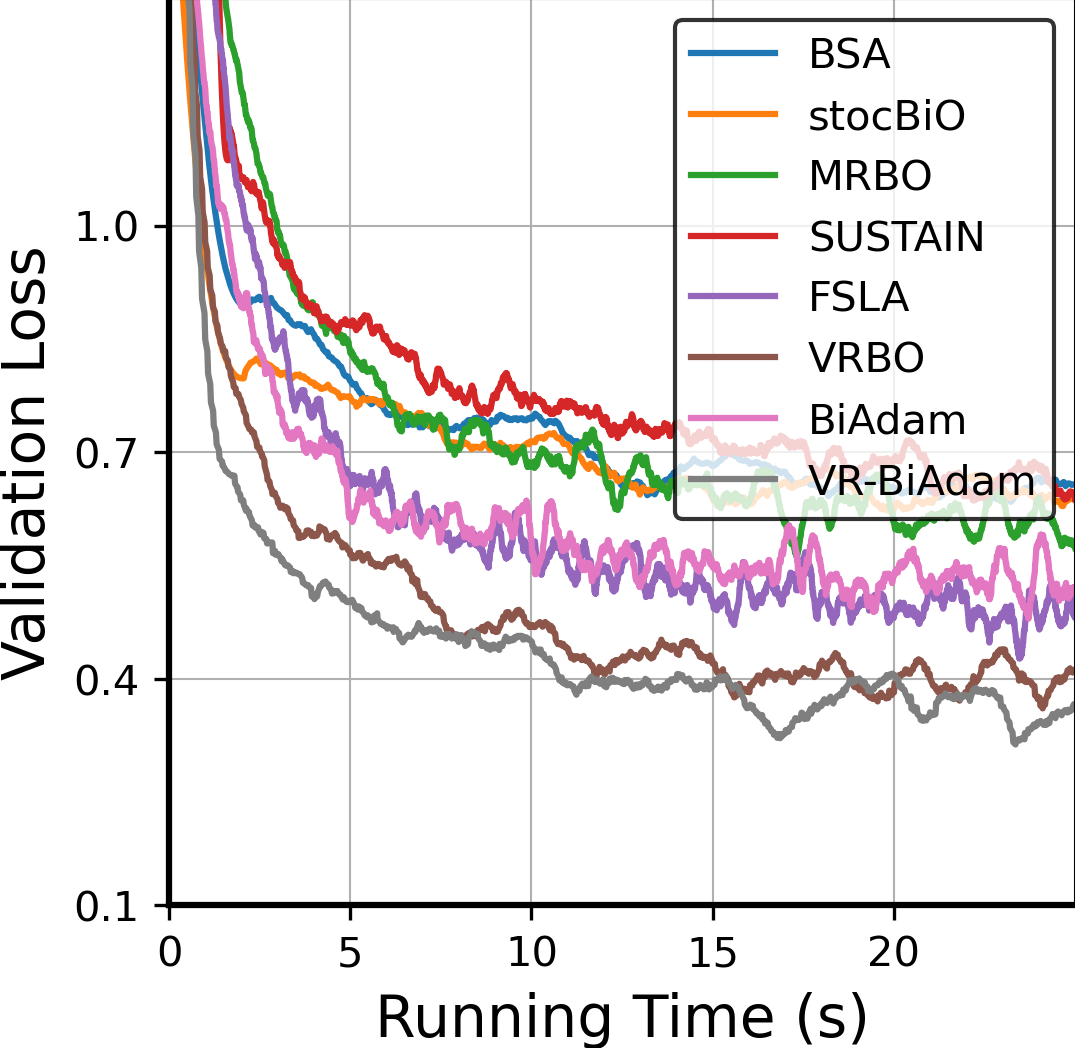}
		\includegraphics[width=0.63\columnwidth]{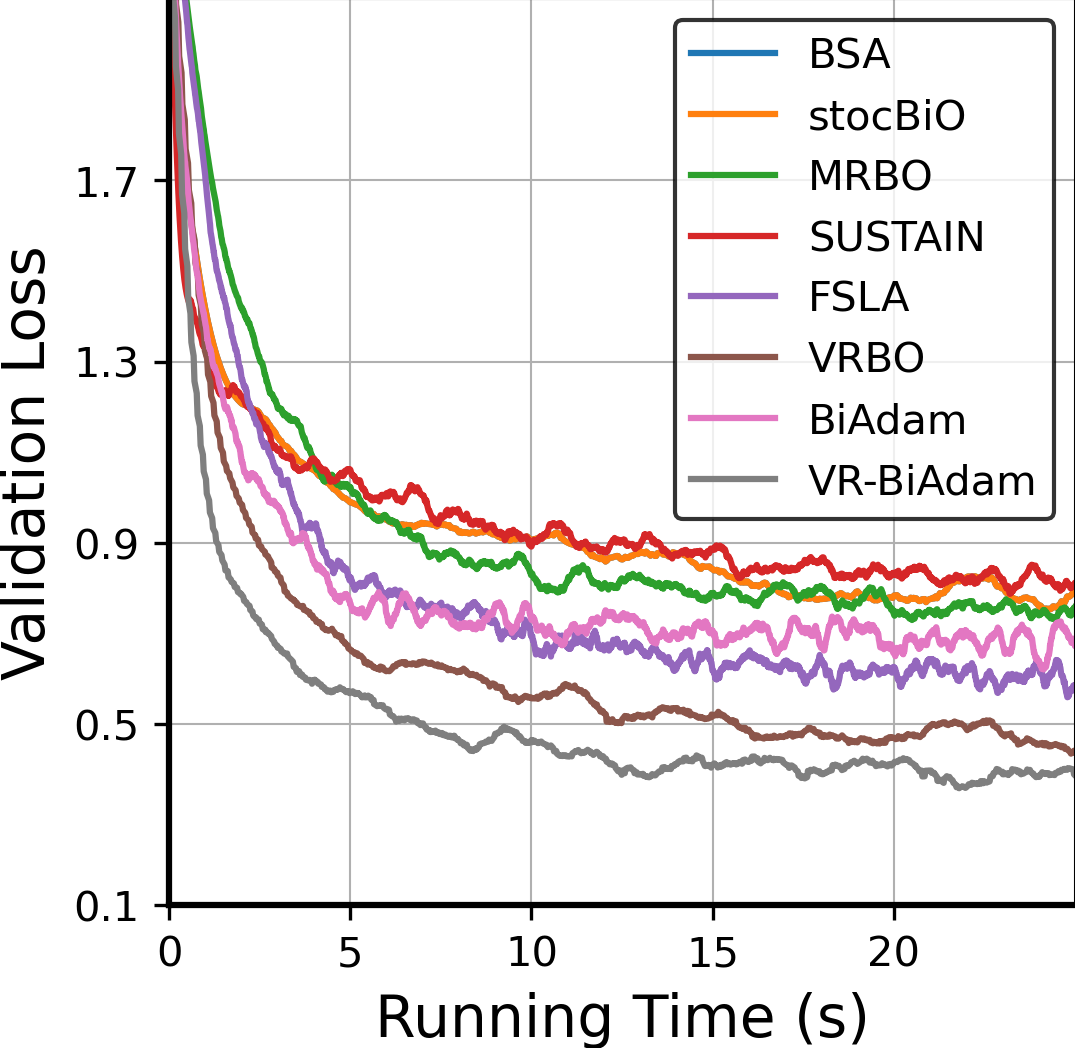}
		\includegraphics[width=0.63\columnwidth]{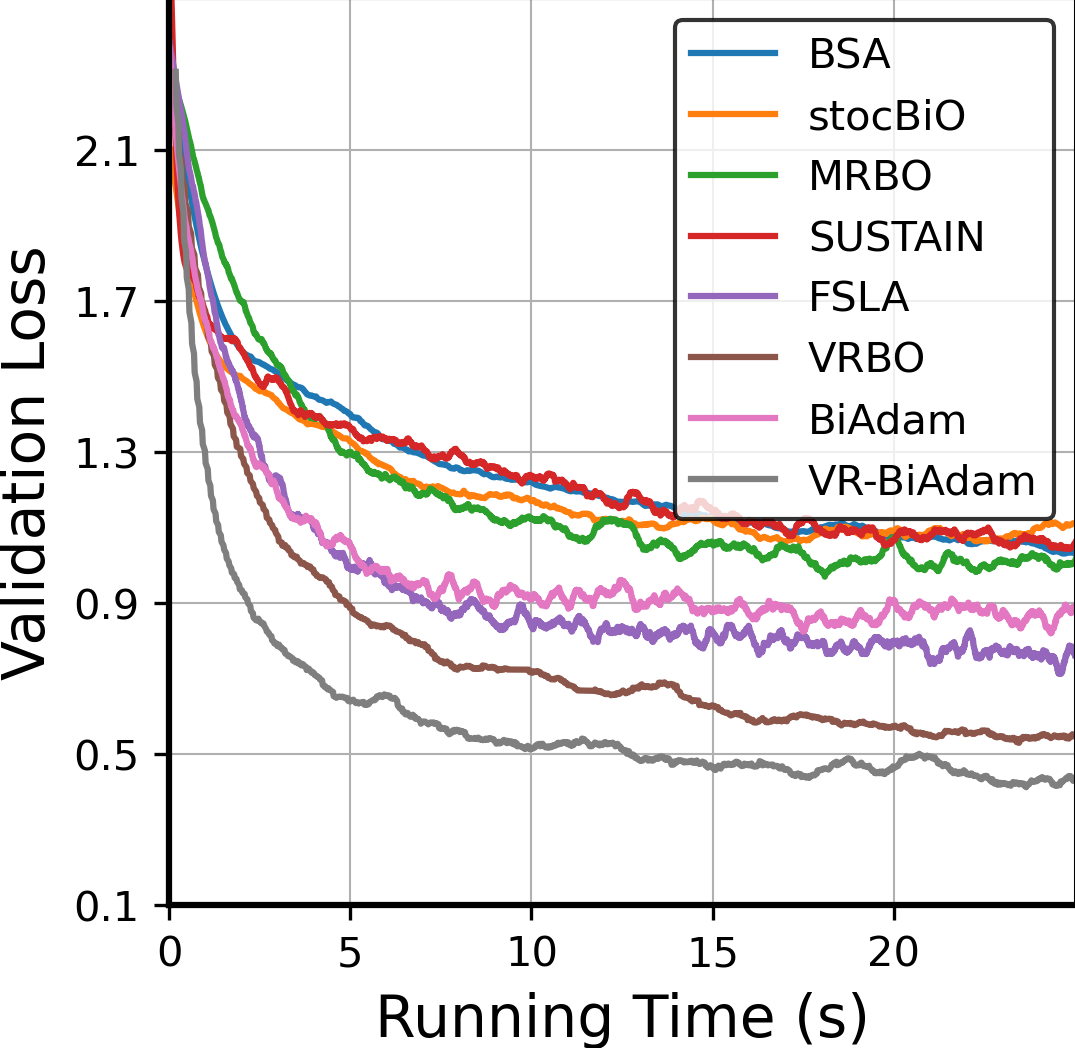}
            \caption{Validation Loss $vs.$ Running Time. We test three values of $\varrho$: 0.8, 0.6, 0.2 from left to right. Larger value of $\varrho$ represents a more noisy setting.}
		  \label{fig:loss}
	\end{center}
\end{figure*}

\section{ Numerical Experiments}
In this section, we perform two hyper-parameter optimization tasks to demonstrate the
efficiency of our algorithms: 1) data hyper-cleaning task
over the MNIST dataset; 2) hyper-representation learning
task over the Omniglot dataset.
In the following experiments, we compare our BiAdam and VR-BiAdam algorithms with the following bilevel optimization algorithms: reverse~\cite{franceschi2018bilevel}, AID-CG~\cite{grazzi2020iteration}, AID-FP~\cite{grazzi2020iteration}, stocBio~\cite{ji2021bilevel}, MRBO~\cite{yang2021provably}, VRBO~\cite{yang2021provably}, FSLA~\cite{li2021fully}, MSTSA/SUSTAIN~\cite{khanduri2021near}, SMB~\cite{guo2021stochastic}, SVRB~\cite{guo2021randomized}. For all methods, we perform grid search over hyper-parameters and choose the best setting. In all experiments,  we use a server with AMD EPYC 7763 64-Core CPU and 1 NVIDIA RTX A5000.

\subsection{ Data Hyper-Cleaning }
In the hyper-cleaning task, we clean a noisy dataset through a bilevel optimization formulation. In this task, we perform data hyper-cleaning over the MNIST dataset~\cite{lecun1998gradient}. Specifically, given the data $(a_i,b_i) \in D_{\mathcal{T}}\cup D_{\mathcal{V}}$,
the formulation of this problem is defined as:
\begin{align*}
  &\min_{z} \bigg\{ f_{val}\big(z,\theta^{\ast}(z)\big)  :=\frac{1}{|D_{\mathcal{V}}|} \sum_{(a_i,b_i)\in D_{\mathcal{V}}} f\big(a_i^T \theta^{\ast}(z), b_i\big) \bigg\} \\
  &\text{s.t. \ } \theta^{\ast}(z) =\arg\min_{\theta} \Bigg\{f_{tr}(z,\theta) \nonumber \\
  & \qquad :=\frac{1}{|D_{\mathcal{T}}|}\Big(\sum_{(a_i,b_i)\in D_{\mathcal{T}}}\sigma(z_i)f(a_i^T\theta,b_i)+C\|\theta\|^2\Big)\Bigg\},
\end{align*}
where $f_{val}(\cdot)$ and $f_{tr}(\cdot)$ denote the training and validation loss functions respectively, and $\theta$ denotes the parameter of model, and $z=\{z_i\}_{i\in \mathcal{D}_\mathcal{T}}$ are hyper-parameters. $D_{\mathcal{T}}$ and $D_{\mathcal{V}}$ are training and validation datasets, respectively. Here $C\geq0$ is a tuning parameter, and $\sigma(\cdot)$ denotes the sigmoid function. In the  experiment, we set $C=0.001$ and let $f_{val}(\cdot)$ and $f_{tr}(\cdot)$ be the cross entropy loss. In particular, we use a training set and a validation set, where each contains 5000 images in our experiments. A portion of the training data are corrupted by randomly changing their labels, and we denote the portion of corrupted images as $\varrho \in (0,1)$.

For training/validation batch-size, we use batch-size of 32, while for VRBO, we choose larger batch-size 5000 and sampling interval is set as 3. For AID-FP, AID-CG and reverse, we use the warm-start trick as our algorithms, \emph{i.e.} the inner variable starts from the state of last iteration. We fine tune the number of inner-loop iterations and set it to be 50 for these algorithms. For MRBO, VRBO, SUSTAIN and our BiAdam/VR-BiAdam, we set $K = 3$ to evaluate the hyper-gradient. For FSLA, $K=1$ as the hyper-gradient is evaluated recursively. As for learning rates, we set 1000 as the outer learning rate for all algorithms except our algorithms which use 0.5 as we change the learning rate adaptively. As for the inner learning rates, we set the stepsize as 0.05 for reverse, AID-CG, stocBiO/AID-FP, MRBO/SUSTAIN, FSLA; we set the stepsize as 0.2 for VRBO; we set the stepsize as 1 for SUSTAIN; we set the stepsize as $0.00025$ for BiAdam and $0.0005$ for VR-BiAdam.

The experimental results are summarized in Fig.~\ref{fig:loss}. As shown by the figure, our BiAdam algorithm outperforms its non-adaptive counterparts such as stocBiO, MRBO and SUSTAIN, furthermore, our VR-BiAdam gets the best performance, where it outperforms VRBO, which requires using large batch-sizes every a few iterations.

\subsection{Hyper-Representation Learning}
In the hyper-representation learning task, we learn a hyper-representation of the data such that a linear classifier can be learned quickly with a small number of data samples. In this task, we perform the hyper-representation learning task over the Omniglot dataset~\cite{lake2015human}.
Given the data $(a_i,b_i) \in D_{\mathcal{T},\xi}\cup D_{\mathcal{V},\xi}$ for all meta task $\xi$, the formulation of this problem is as follows:
\begin{align*}
  &\min_{z} \bigg\{ f_{val}\big(z,\theta^{\ast}(z)\big)\nonumber \\
  & \qquad :=\mathbb{E}_{\xi}\Big[\frac{1}{|D_{\mathcal{V}, \xi}|} \sum_{(a_i,b_i)\in D_{\mathcal{V}, \xi}} f\big(\theta^{*}(z)^{T}\phi(a_i; z),b_i\big)\Big] \bigg\} \\
  &\text{s.t. \ } \theta^{\ast}(z) =\arg\min_{\theta} \Bigg\{f_{tr}(z, \theta)\nonumber \\
  & \qquad :=\frac{1}{|D_{\mathcal{T},\xi}|}\sum_{(a_i,b_i)\in D_{\mathcal{T}}, \xi}f\big(\theta^{T}\phi(a_i; z),b_i\big)+C\|\theta\|^2 \Bigg\},
\end{align*}
where $f(\cdot)$ denotes the cross entropy loss, and $\theta$ denotes the parameter of model, and $z$ denotes the parameter vector of the representation mapping $\phi(\cdot,\cdot)$. $D_{\mathcal{T},\xi}$ and $D_{\mathcal{V},\xi}$ are training and validation datasets respectively for a randomly sampled meta task $\xi$. Here $C\geq0$ is a tuning parameter to guarantee the inner problem to be strongly convex. In experiment, we set $C=0.01$.

In every hyper-iteration, we choose 4 meta tasks, while for VRBO, we choose larger batch-size 16 and sampling interval is set as 3. For stocBiO/AID-FP, AID-CG and reverse, we use the warm-start trick as our algorithms, \emph{i.e.} the inner variable starts from the state of last iteration. We fine tune the number of inner-loop iterations and set it to be 16 for these algorithms. For MRBO, VRBO, SUSTAIN and our algorithms, we set $K = 5$ to evaluate the hyper-gradient. For FSLA, $K=1$ as the hyper-gradient is evaluated recursively. As for learning rates, we set 1000 as the outer learning rate for all algorithms except our algorithms which use 0.001 as we change the learning rate adaptively. As for the inner learning rates, we set the stepsize as 0.4 for all algorithms.

\begin{figure*}[ht]
	\begin{center}
		\includegraphics[width=0.76\columnwidth]{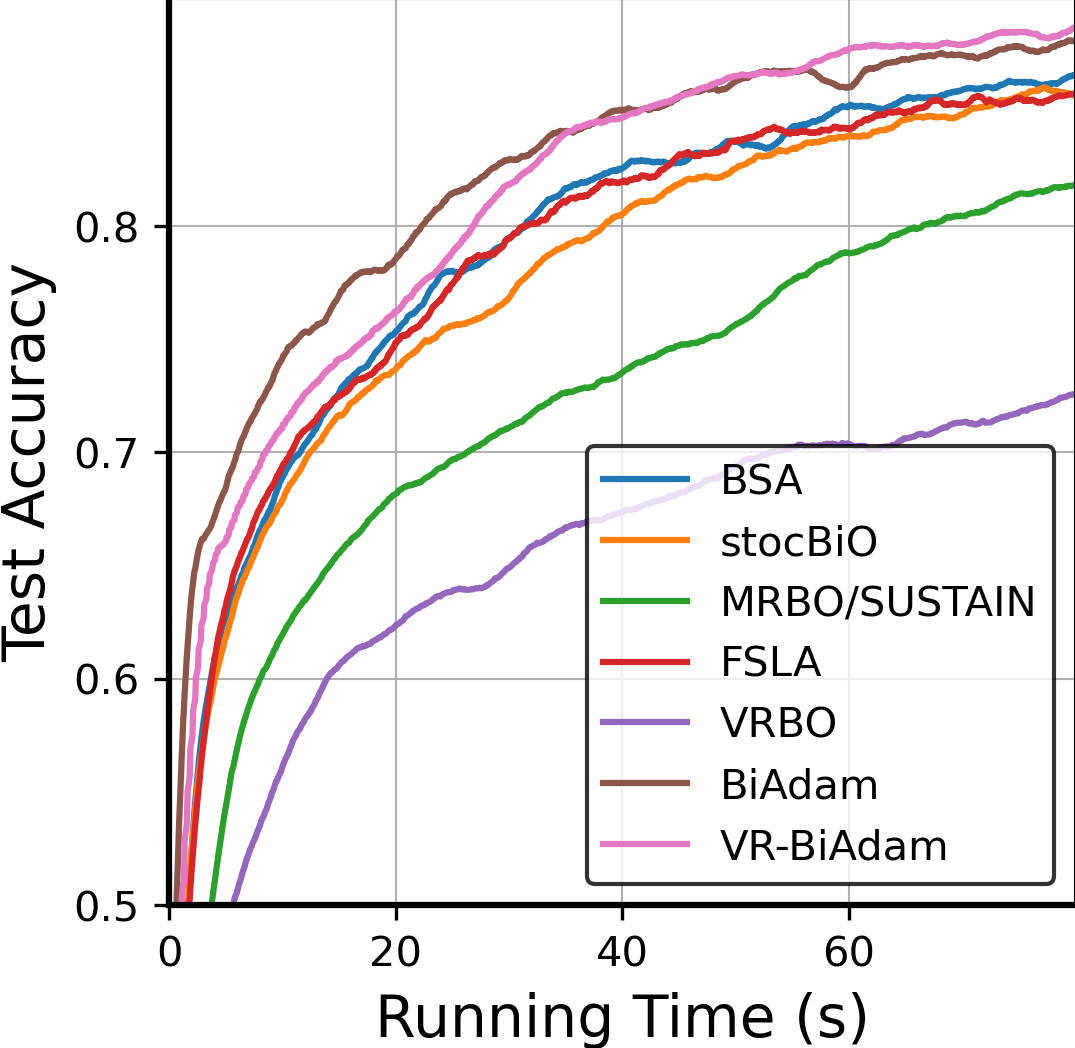} \qquad
		\includegraphics[width=0.76\columnwidth]{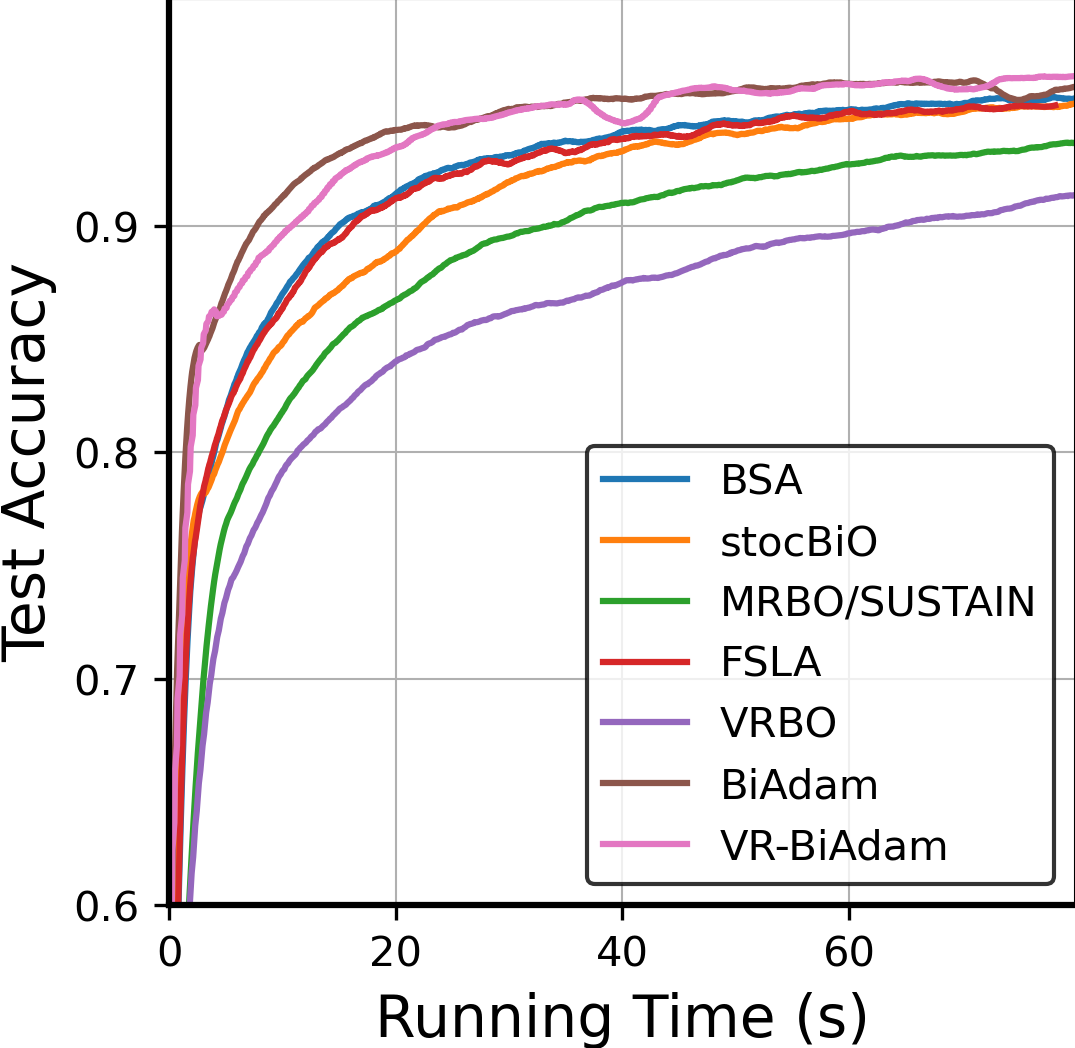}
		\includegraphics[width=0.76\columnwidth]{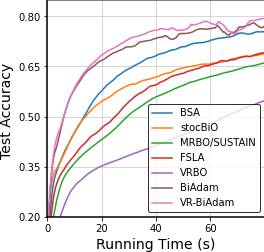} \qquad
		\includegraphics[width=0.76\columnwidth]{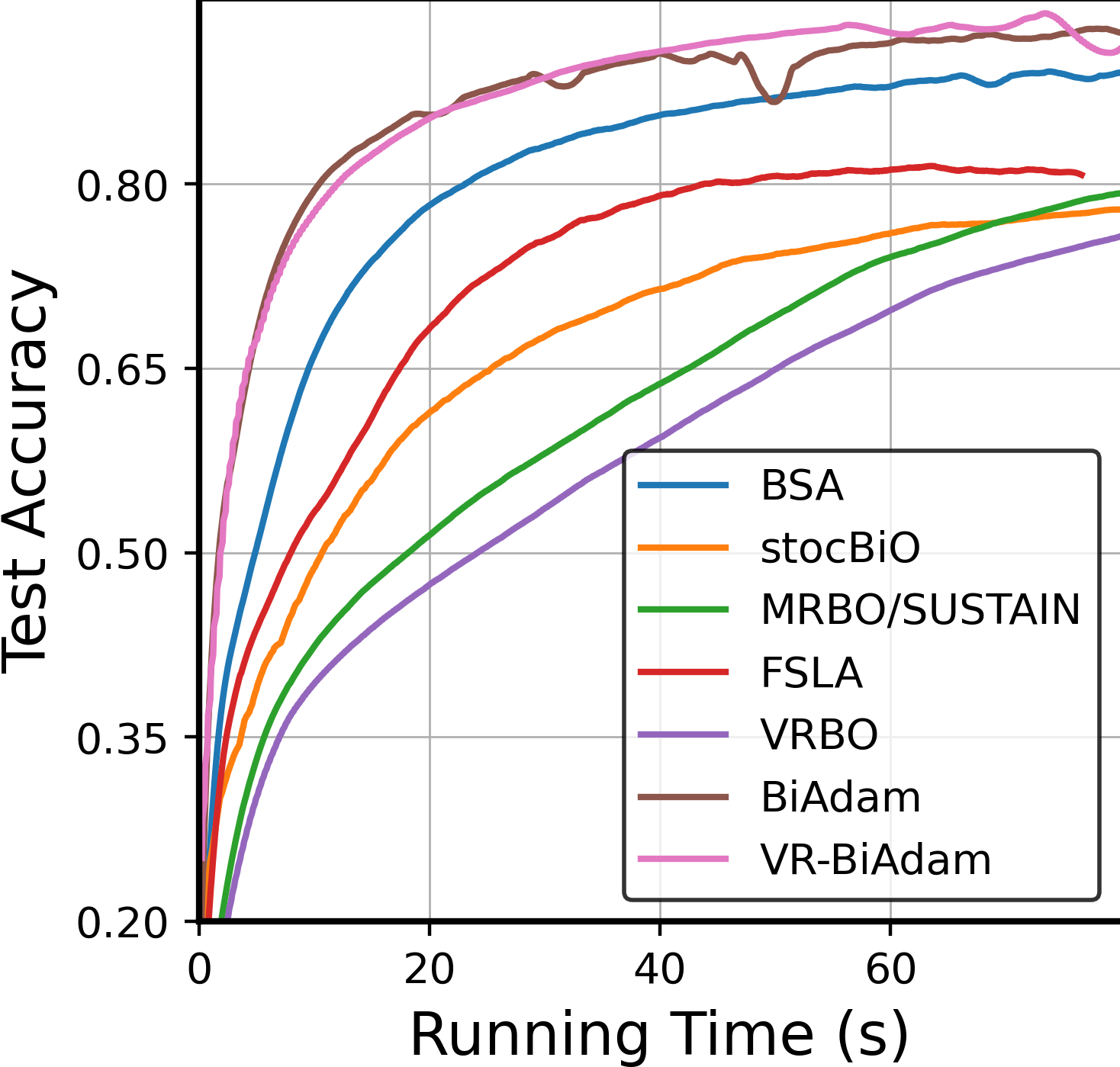}
            \caption{Test Accuracy $vs.$ Running Time. The plots corresponds to 5-way-1-shot, 5-way-5-shot, 20-way-1-shot and 20-way-5-shots from left to right.}
		\label{fig:acc-omniglot}
	\end{center}
\end{figure*}

The experimental results are summarized in Fig.~\ref{fig:acc-omniglot}. As shown by the figure, both our BiAdam and VR-BiAdam algorithms outperform other baselines.

\section{ Conclusions }
In this paper, we proposed a class of novel adaptive bilevel optimization methods for nonconvex-strongly-convex bilevel optimization problems.
Our methods use unified adaptive matrices including many types of adaptive learning rates, and
can flexibly use the momentum and variance reduced techniques.
Moreover, we provided a useful convergence analysis framework for both the constrained and unconstrained bilevel optimization. Meanwhile, our VR-BiAdam algorithm reaches the best known sample complexity of $\tilde{O}(\epsilon^{-3})$ for finding an $\epsilon$-stationary solution.
In particular, our algorithms use the the momentum techniques to update variables and estimate stochastic gradients simultaneously.


%

%

\ifCLASSOPTIONcaptionsoff
  \newpage
\fi

\bibliographystyle{IEEEtran}
\bibliography{IEEEabrv,Bilevel}

\begin{thebibliography}{10}
\providecommand{\url}[1]{#1}
\csname url@samestyle\endcsname
\providecommand{\newblock}{\relax}
\providecommand{\bibinfo}[2]{#2}
\providecommand{\BIBentrySTDinterwordspacing}{\spaceskip=0pt\relax}
\providecommand{\BIBentryALTinterwordstretchfactor}{4}
\providecommand{\BIBentryALTinterwordspacing}{\spaceskip=\fontdimen2\font plus
\BIBentryALTinterwordstretchfactor\fontdimen3\font minus
  \fontdimen4\font\relax}
\providecommand{\BIBforeignlanguage}[2]{{%
\expandafter\ifx\csname l@#1\endcsname\relax
\typeout{** WARNING: IEEEtran.bst: No hyphenation pattern has been}%
\typeout{** loaded for the language `#1'. Using the pattern for}%
\typeout{** the default language instead.}%
\else
\language=\csname l@#1\endcsname
\fi
#2}}
\providecommand{\BIBdecl}{\relax}
\BIBdecl

\bibitem{shaban2019truncated}
A.~Shaban, C.-A. Cheng, N.~Hatch, and B.~Boots, ``Truncated back-propagation
  for bilevel optimization,'' in \emph{The 22nd International Conference on
  Artificial Intelligence and Statistics}.\hskip 1em plus 0.5em minus
  0.4em\relax PMLR, 2019, pp. 1723--1732.

\bibitem{ji2021bilevel}
K.~Ji, J.~Yang, and Y.~Liang, ``Bilevel optimization: Convergence analysis and
  enhanced design,'' in \emph{International Conference on Machine
  Learning}.\hskip 1em plus 0.5em minus 0.4em\relax PMLR, 2021, pp. 4882--4892.

\bibitem{liu2021investigating}
R.~Liu, J.~Gao, J.~Zhang, D.~Meng, and Z.~Lin, ``Investigating bi-level
  optimization for learning and vision from a unified perspective: A survey and
  beyond,'' \emph{IEEE Transactions on Pattern Analysis and Machine
  Intelligence}, 2021.

\bibitem{hong2020two}
M.~Hong, H.-T. Wai, Z.~Wang, and Z.~Yang, ``A two-timescale framework for
  bilevel optimization: Complexity analysis and application to actor-critic,''
  \emph{arXiv preprint arXiv:2007.05170}, 2020.

\bibitem{franceschi2018bilevel}
L.~Franceschi, P.~Frasconi, S.~Salzo, R.~Grazzi, and M.~Pontil, ``Bilevel
  programming for hyperparameter optimization and meta-learning,'' in
  \emph{International Conference on Machine Learning}.\hskip 1em plus 0.5em
  minus 0.4em\relax PMLR, 2018, pp. 1568--1577.

\bibitem{liu2018darts}
H.~Liu, K.~Simonyan, and Y.~Yang, ``Darts: Differentiable architecture
  search,'' \emph{arXiv preprint arXiv:1806.09055}, 2018.

\bibitem{ghadimi2018approximation}
S.~Ghadimi and M.~Wang, ``Approximation methods for bilevel programming,''
  \emph{arXiv preprint arXiv:1802.02246}, 2018.

\bibitem{chen2022single}
T.~Chen, Y.~Sun, Q.~Xiao, and W.~Yin, ``A single-timescale method for
  stochastic bilevel optimization,'' in \emph{International Conference on
  Artificial Intelligence and Statistics}.\hskip 1em plus 0.5em minus
  0.4em\relax PMLR, 2022, pp. 2466--2488.

\bibitem{guo2021stochastic}
Z.~Guo, Y.~Xu, W.~Yin, R.~Jin, and T.~Yang, ``On stochastic moving-average
  estimators for non-convex optimization,'' \emph{arXiv preprint
  arXiv:2104.14840}, 2021.

\bibitem{khanduri2021near}
P.~Khanduri, S.~Zeng, M.~Hong, H.-T. Wai, Z.~Wang, and Z.~Yang, ``A
  near-optimal algorithm for stochastic bilevel optimization via
  double-momentum,'' \emph{Advances in Neural Information Processing Systems},
  vol.~34, pp. 30\,271--30\,283, 2021.

\bibitem{guo2021randomized}
Z.~Guo and T.~Yang, ``Randomized stochastic variance-reduced methods for
  stochastic bilevel optimization,'' \emph{arXiv preprint arXiv:2105.02266},
  2021.

\bibitem{yang2021provably}
J.~Yang, K.~Ji, and Y.~Liang, ``Provably faster algorithms for bilevel
  optimization,'' \emph{Advances in Neural Information Processing Systems},
  vol.~34, pp. 13\,670--13\,682, 2021.

\bibitem{cutkosky2019momentum}
A.~Cutkosky and F.~Orabona, ``Momentum-based variance reduction in non-convex
  sgd,'' \emph{Advances in neural information processing systems}, vol.~32,
  2019.

\bibitem{zhang2021revisiting}
Y.~Zhang, G.~Zhang, P.~Khanduri, M.~Hong, S.~Chang, and S.~Liu, ``Revisiting
  and advancing fast adversarial training through the lens of bi-level
  optimization,'' \emph{arXiv preprint arXiv:2112.12376}, 2021.

\bibitem{colson2007overview}
B.~Colson, P.~Marcotte, and G.~Savard, ``An overview of bilevel optimization,''
  \emph{Annals of operations research}, vol. 153, no.~1, pp. 235--256, 2007.

\bibitem{kunapuli2008classification}
G.~Kunapuli, K.~P. Bennett, J.~Hu, and J.-S. Pang, ``Classification model
  selection via bilevel programming,'' \emph{Optimization Methods \& Software},
  vol.~23, no.~4, pp. 475--489, 2008.

\bibitem{chen2021closing}
T.~Chen, Y.~Sun, and W.~Yin, ``Closing the gap: Tighter analysis of alternating
  stochastic gradient methods for bilevel problems,'' \emph{Advances in Neural
  Information Processing Systems}, vol.~34, pp. 25\,294--25\,307, 2021.

\bibitem{liu2021towards}
R.~Liu, Y.~Liu, S.~Zeng, and J.~Zhang, ``Towards gradient-based bilevel
  optimization with non-convex followers and beyond,'' \emph{Advances in Neural
  Information Processing Systems}, vol.~34, pp. 8662--8675, 2021.

\bibitem{liu2022general}
R.~Liu, P.~Mu, X.~Yuan, S.~Zeng, and J.~Zhang, ``A general descent aggregation
  framework for gradient-based bi-level optimization,'' \emph{IEEE Transactions
  on Pattern Analysis and Machine Intelligence}, 2022.

\bibitem{li2021fully}
J.~Li, B.~Gu, and H.~Huang, ``A fully single loop algorithm for bilevel
  optimization without hessian inverse,'' \emph{arXiv preprint
  arXiv:2112.04660}, 2021.

\bibitem{ji2022lower}
K.~Ji and Y.~Liang, ``Lower bounds and accelerated algorithms for bilevel
  optimization,'' \emph{Journal of Machine Learning Research}, vol.~23, pp.
  1--56, 2022.

\bibitem{chen2022fast}
Z.~Chen, B.~Kailkhura, and Y.~Zhou, ``A fast and convergent proximal algorithm
  for regularized nonconvex and nonsmooth bi-level optimization,'' \emph{arXiv
  preprint arXiv:2203.16615}, 2022.

\bibitem{duchi2011adaptive}
J.~Duchi, E.~Hazan, and Y.~Singer, ``Adaptive subgradient methods for online
  learning and stochastic optimization.'' \emph{Journal of machine learning
  research}, vol.~12, no.~7, 2011.

\bibitem{kingma2014adam}
D.~P. Kingma and J.~Ba, ``Adam: A method for stochastic optimization,''
  \emph{arXiv preprint arXiv:1412.6980}, 2014.

\bibitem{loshchilov2018decoupled}
I.~Loshchilov and F.~Hutter, ``Decoupled weight decay regularization,'' in
  \emph{International Conference on Learning Representations}, 2018.

\bibitem{zhuang2020adabelief}
J.~Zhuang, T.~Tang, Y.~Ding, S.~C. Tatikonda, N.~Dvornek, X.~Papademetris, and
  J.~Duncan, ``Adabelief optimizer: Adapting stepsizes by the belief in
  observed gradients,'' \emph{Advances in neural information processing
  systems}, vol.~33, pp. 18\,795--18\,806, 2020.

\bibitem{reddi2019convergence}
S.~J. Reddi, S.~Kale, and S.~Kumar, ``On the convergence of adam and beyond,''
  \emph{arXiv preprint arXiv:1904.09237}, 2019.

\bibitem{chen2019convergence}
X.~Chen, S.~Liu, R.~Sun, and M.~Hong, ``On the convergence of a class of
  adam-type algorithms for non-convex optimization,'' in \emph{7th
  International Conference on Learning Representations (ICLR)}, 2019.

\bibitem{chen2018closing}
J.~Chen, D.~Zhou, Y.~Tang, Z.~Yang, Y.~Cao, and Q.~Gu, ``Closing the
  generalization gap of adaptive gradient methods in training deep neural
  networks,'' \emph{arXiv preprint arXiv:1806.06763}, 2018.

\bibitem{ward2019adagrad}
R.~Ward, X.~Wu, and L.~Bottou, ``Adagrad stepsizes: Sharp convergence over
  nonconvex landscapes,'' in \emph{International Conference on Machine
  Learning}.\hskip 1em plus 0.5em minus 0.4em\relax PMLR, 2019, pp. 6677--6686.

\bibitem{huang2021super}
F.~Huang, J.~Li, and H.~Huang, ``Super-adam: faster and universal framework of
  adaptive gradients,'' \emph{Advances in Neural Information Processing
  Systems}, vol.~34, pp. 9074--9085, 2021.

\bibitem{guo2021novel}
Z.~Guo, Y.~Xu, W.~Yin, R.~Jin, and T.~Yang, ``A novel convergence analysis for
  algorithms of the adam family,'' \emph{arXiv preprint arXiv:2112.03459},
  2021.

\bibitem{censor1992proximal}
Y.~Censor and S.~A. Zenios, ``Proximal minimization algorithm
  withd-functions,'' \emph{Journal of Optimization Theory and Applications},
  vol.~73, no.~3, pp. 451--464, 1992.

\bibitem{li2019convergence}
X.~Li and F.~Orabona, ``On the convergence of stochastic gradient descent with
  adaptive stepsizes,'' in \emph{The 22nd International Conference on
  Artificial Intelligence and Statistics}.\hskip 1em plus 0.5em minus
  0.4em\relax PMLR, 2019, pp. 983--992.

\bibitem{censor1981iterative}
Y.~Censor and A.~Lent, ``An iterative row-action method for interval convex
  programming,'' \emph{Journal of Optimization theory and Applications},
  vol.~34, no.~3, pp. 321--353, 1981.

\bibitem{ghadimi2016mini}
S.~Ghadimi, G.~Lan, and H.~Zhang, ``Mini-batch stochastic approximation methods
  for nonconvex stochastic composite optimization,'' \emph{Mathematical
  Programming}, vol. 155, no. 1-2, pp. 267--305, 2016.

\bibitem{grazzi2020iteration}
R.~Grazzi, L.~Franceschi, M.~Pontil, and S.~Salzo, ``On the iteration
  complexity of hypergradient computation,'' in \emph{International Conference
  on Machine Learning}.\hskip 1em plus 0.5em minus 0.4em\relax PMLR, 2020, pp.
  3748--3758.

\bibitem{lecun1998gradient}
Y.~LeCun, L.~Bottou, Y.~Bengio, and P.~Haffner, ``Gradient-based learning
  applied to document recognition,'' \emph{Proceedings of the IEEE}, vol.~86,
  no.~11, pp. 2278--2324, 1998.

\bibitem{lake2015human}
B.~M. Lake, R.~Salakhutdinov, and J.~B. Tenenbaum, ``Human-level concept
  learning through probabilistic program induction,'' \emph{Science}, vol. 350,
  no. 6266, pp. 1332--1338, 2015.

\end{thebibliography}

\vfill


\newpage

\begin{onecolumn}

\appendices

\begin{appendices}
\section{Detailed Proofs in Convergence Analysis}

In this section, we provide the detailed convergence analysis of our algorithms.
We first review and provide some useful lemmas.

Given a $\rho$-strongly convex function $\phi(x)$, we define a prox-function (Bregman distance) \cite{censor1981iterative,censor1992proximal}
associated with $\phi(x)$ as follows:
\begin{align}
 V(z,x) = \phi(z) - \big[ \phi(x) + \langle\nabla \phi(x),z-x\rangle\big].
\end{align}
Then we define a generalized projection problem as in \cite{ghadimi2016mini}:
\begin{align} \label{eq:A1}
 x^* = \arg\min_{z\in \mathcal{X}} \big\{\langle z,w\rangle + \frac{1}{\gamma}V(z,x) + h(z)\big\},
\end{align}
where $ \mathcal{X} \subseteq \mathbb{R}^d$, $w \in \mathbb{R}^d$ and $\gamma>0$.
Here $h(x)$ is convex and possibly nonsmooth function.
At the same time, we define a generalized gradient as follows:
\begin{align}
 \mathcal{G}_{\mathcal{X}}(x,w,\gamma) = \frac{1}{\gamma}(x-x^*).
\end{align}

\begin{lemma} \label{lem:A1}
(Lemma 1 in \cite{ghadimi2016mini})
Let $x^*$ be given in (\ref{eq:A1}). Then, for any $x\in \mathcal{X}$, $w\in \mathbb{R}^d$ and $\gamma >0$,
we have
\begin{align}
 \langle w,  \mathcal{G}_{\mathcal{X}}(x,w,\gamma)\rangle \geq \rho \|\mathcal{G}_{\mathcal{X}}(x,w,\gamma)\|^2
 + \frac{1}{\gamma}\big[h(x^*)-h(x)\big],
\end{align}
where $\rho>0$ depends on $\rho$-strongly convex function $\phi(x)$.
\end{lemma}
When $h(x)=0$, in the above lemma \ref{lem:A1}, we have
\begin{align}
 \langle w, \mathcal{G}_{\mathcal{X}}(x,w,\gamma)\rangle \geq \rho \|\mathcal{G}_{\mathcal{X}}(x,w,\gamma)\|^2.
\end{align}

\begin{lemma} \label{lem:A2}
 (Restatement of Lemma 5)
When the gradient estimator $w_t$ generated from Algorithm \ref{alg:1} or \ref{alg:2}, for all $t\geq 1$,
 we have
 \begin{align}
 \|w_t-\nabla F(x_t)\|^2 \leq L^2_0\|y^*(x_t)-y_t\|^2 + 2\|w_t-\bar{\nabla} f(x_t,y_t)\|^2,
\end{align}
where $L^2_0 = 8\big(L^2_f+ \frac{L^2_{gxy}C^2_{fy}}{\mu^2} + \frac{L^2_{gyy} C^2_{gxy}C^2_{fy}}{\mu^4} +
 \frac{L^2_fC^2_{gxy}}{\mu^2}\big)$.
\end{lemma}

\begin{proof}
 We first consider the term $\|\nabla F(x_t)-\bar{\nabla} f(x_t,y_t)\|^2$. Since $\nabla f(x_t,y^*(x_t))=\nabla F(x_t)$, we have
\begin{align}
 & \|\nabla f(x_t,y^*(x_t))-\bar{\nabla} f(x_t,y_t)\|^2 \nonumber \\
 & = \|\nabla_xf(x_t,y^*(x_t)) - \nabla^2_{xy}g(x_t,y^*(x_t)) \big(\nabla^2_{yy}g(x_t,y^*(x_t))\big)^{-1}\nabla_yf(x_t,y^*(x)) \nonumber \\
 & \quad  - \nabla_xf(x_t,y_t) + \nabla^2_{xy}g(x_t,y_t) \big(\nabla^2_{yy}g(x_t,y_t)\big)^{-1}\nabla_yf(x_t,y_t)\|^2 \nonumber \\
 & = \|\nabla_xf(x_t,y^*(x_t)) - \nabla_xf(x_t,y_t) - \nabla^2_{xy}g(x_t,y^*(x_t)) \big(\nabla^2_{yy}g(x_t,y^*(x_t))\big)^{-1}\nabla_yf(x_t,y^*(x_t)) \nonumber \\
 & \quad  + \nabla^2_{xy}g(x_t,y_t) \big(\nabla^2_{yy}g(x_t,y^*(x_t))\big)^{-1}\nabla_yf(x_t,y^*(x_t)) - \nabla^2_{xy}g(x_t,y_t) \big(\nabla^2_{yy}g(x_t,y^*(x_t))\big)^{-1}\nabla_yf(x_t,y^*(x_t)) \nonumber \\
 & \quad + \nabla^2_{xy}g(x_t,y_t) \big(\nabla^2_{yy}g(x_t,y_t)\big)^{-1}\nabla_yf(x_t,y^*(x_t)) -  \nabla^2_{xy}g(x_t,y_t) \big(\nabla^2_{yy}g(x_t,y_t)\big)^{-1}\nabla_yf(x_t,y^*(x_t)) \nonumber \\
 & \quad + \nabla^2_{xy}g(x_t,y_t) \big(\nabla^2_{yy}g(x_t,y_t)\big)^{-1}\nabla_yf(x_t,y_t)\|^2 \nonumber \\
 & \leq 4\|\nabla_xf(x_t,y^*(x_t)) - \nabla_xf(x_t,y_t)\|^2 + \frac{4C^2_{fy}}{\mu^2}\|\nabla^2_{xy}g(x_t,y^*(x_t))-\nabla^2_{xy}g(x_t,y_t)\|^2
 \nonumber \\
 & \quad + \frac{4C^2_{gxy}C^2_{fy}}{\mu^4}\|\nabla^2_{yy}g(x_t,y^*(x_t))-\nabla^2_{yy}g(x_t,y_t)\|^2 +
 \frac{4C^2_{gxy}}{\mu^2}\|\nabla_yf(x_t,y^*(x_t))-\nabla_yf(x_t,y_t)\|^2 \nonumber \\
 & \leq 4\big(L^2_f+ \frac{L^2_{gxy}C^2_{fy}}{\mu^2} + \frac{L^2_{gyy} C^2_{gxy}C^2_{fy}}{\mu^4} +
 \frac{L^2_fC^2_{gxy}}{\mu^2}\big)\|y^*(x_t)-y_t\|^2 \nonumber \\
 & = 4\bar{L}^2\|y^*(x_t)-y_t\|^2,
\end{align}
where the second last inequality is due to Assumptions \ref{ass:1}, \ref{ass:2} and \ref{ass:4};
the last equality holds by $ \bar{L}^2=L^2_f+ \frac{L^2_{gxy}C^2_{fy}}{\mu^2} + \frac{L^2_{gyy} C^2_{gxy}C^2_{fy}}{\mu^4} +
 \frac{L^2_fC^2_{gxy}}{\mu^2}$.

Then we have
\begin{align}
 \|w_t-\nabla F(x_t)\|^2 & = \|w_t-\bar{\nabla} f(x_t,y_t) + \bar{\nabla} f(x_t,y_t)-\nabla F(x_t)\|^2 \nonumber \\
 & \leq 2\|w_t-\bar{\nabla} f(x_t,y_t)\|^2 + 2\|\bar{\nabla} f(x_t,y_t)-\nabla F(x_t)\|^2 \nonumber \\
 & \leq 2\|w_t-\bar{\nabla} f(x_t,y_t)\|^2 + 8\bar{L}^2\|y^*(x_t)-y_t\|^2.
\end{align}

\end{proof}

\begin{lemma} \label{lem:A3}
Under the Assumptions \ref{ass:1}, \ref{ass:2}, \ref{ass:4}, we have
\begin{align}
  \|\bar{\nabla} f(x_{t+1},y_{t+1})-\bar{\nabla} f(x_t,y_t)\|^2 \leq L^2_0\big(\|x_{t+1}-x_t\|^2 + \|y_{t+1}-y_t\|^2\big),
\end{align}
where $L^2_0 = 8\big( L^2_f+ \frac{L^2_{gxy}C^2_{fy}}{\mu^2} + \frac{L^2_{gyy} C^2_{gxy}C^2_{fy}}{\mu^4} +
 \frac{L^2_fC^2_{gxy}}{\mu^2} \big)$.
\end{lemma}

\begin{proof}
\begin{align}
 & \|\bar{\nabla} f(x_{t+1},y_{t+1})-\bar{\nabla} f(x_t,y_t)\|^2 \nonumber \\
 & = \|\nabla_xf(x_{t+1},y_{t+1}) - \nabla^2_{xy}g(x_{t+1},y_{t+1}) \big(\nabla^2_{yy}g(x_{t+1},y_{t+1})\big)^{-1}\nabla_yf(x_{t+1},y_{t+1}) \nonumber \\
 & \quad  - \nabla_xf(x_t,y_t) + \nabla^2_{xy}g(x_t,y_t) \big(\nabla^2_{yy}g(x_t,y_t)\big)^{-1}\nabla_yf(x_t,y_t)\|^2 \nonumber \\
 & = \|\nabla_xf(x_{t+1},y_{t+1}) - \nabla_xf(x_t,y_t) - \nabla^2_{xy}g(x_{t+1},y_{t+1}) \big(\nabla^2_{yy}g(x_{t+1},y_{t+1})\big)^{-1}\nabla_yf(x_{t+1},y_{t+1}) \nonumber \\
 & \quad  + \nabla^2_{xy}g(x_t,y_t) \big(\nabla^2_{yy}g(x_{t+1},y_{t+1})\big)^{-1}\nabla_yf(x_{t+1},y_{t+1}) - \nabla^2_{xy}g(x_t,y_t) \big(\nabla^2_{yy}g(x_{t+1},y_{t+1})\big)^{-1}\nabla_yf(x_{t+1},y_{t+1}) \nonumber \\
 & \quad + \nabla^2_{xy}g(x_t,y_t) \big(\nabla^2_{yy}g(x_t,y_t)\big)^{-1}\nabla_yf(x_{t+1},y_{t+1}) -  \nabla^2_{xy}g(x_t,y_t) \big(\nabla^2_{yy}g(x_t,y_t)\big)^{-1}\nabla_yf(x_{t+1},y_{t+1}) \nonumber \\
 & \quad + \nabla^2_{xy}g(x_t,y_t) \big(\nabla^2_{yy}g(x_t,y_t)\big)^{-1}\nabla_yf(x_t,y_t)\|^2 \nonumber \\
 & \leq 4\|\nabla_xf(x_{t+1},y_{t+1}) - \nabla_xf(x_t,y_t)\|^2 + \frac{4C^2_{fy}}{\mu^2}\|\nabla^2_{xy}g(x_{t+1},y_{t+1})-\nabla^2_{xy}g(x_t,y_t)\|^2
 \nonumber \\
 & \quad + \frac{4C^2_{gxy}C^2_{fy}}{\mu^4}\|\nabla^2_{yy}g(x_{t+1},y_{t+1})-\nabla^2_{yy}g(x_t,y_t)\|^2 +
 \frac{4C^2_{gxy}}{\mu^2}\|\nabla_yf(x_{t+1},y_{t+1})-\nabla_yf(x_t,y_t)\|^2 \nonumber \\
 & \leq 8L^2_f\big(\|x_{t+1}-x_t\|^2 + \|y_{t+1}-y_t\|^2\big) + \frac{8L^2_{gxy}C^2_{fy}}{\mu^2}\big(\|x_{t+1}-x_t\|^2 + \|y_{t+1}-y_t\|^2\big)
 \nonumber \\
 & \quad + \frac{8L^2_{gyy} C^2_{gxy}C^2_{fy}}{\mu^4}\big(\|x_{t+1}-x_t\|^2 + \|y_{t+1}-y_t\|^2\big) +
 \frac{8L^2_fC^2_{gxy}}{\mu^2}\big(\|x_{t+1}-x_t\|^2 + \|y_{t+1}-y_t\|^2\big) \nonumber \\
 & = L^2_0\big(\|x_{t+1}-x_t\|^2 + \|y_{t+1}-y_t\|^2\big),
\end{align}
where the first inequality holds by the Assumptions \ref{ass:1} and \ref{ass:4}, and the second inequality holds by the Assumption \ref{ass:2}.

\end{proof}

\begin{lemma} \label{lem:A4}
 Suppose that the sequence $\{x_t,y_t\}_{t=1}^T$ be generated from Algorithm \ref{alg:1} or \ref{alg:2}.
 Let $0<\eta_t \leq 1$ and $0< \gamma \leq \frac{\rho}{2L\eta_t}$,
 then we have
 \begin{align}
  F(x_{t+1}) \leq F(x_t) + \frac{\eta_t\gamma}{\rho}\|\nabla F(x_t)-w_t\|^2 -\frac{\rho\eta_t}{2\gamma}\|\tilde{x}_{t+1}-x_t\|^2.
 \end{align}
\end{lemma}

\begin{proof}
According to Lemma \ref{lem:1}, the function $F(x)$ is $L$-smooth. Thus we have
 \begin{align} \label{eq:B1}
  F(x_{t+1}) & \leq F(x_t) + \langle\nabla F(x_t), x_{t+1}-x_t\rangle + \frac{L}{2}\|x_{t+1}-x_t\|^2 \\
  & = F(x_t)+ \langle \nabla F(x_t),\eta_t(\tilde{x}_{t+1}-x_t)\rangle + \frac{L}{2}\|\eta_t(\tilde{x}_{t+1}-x_t)\|^2 \nonumber \\
  & = F(x_t) + \eta_t\underbrace{\langle w_t,\tilde{x}_{t+1}-x_t\rangle}_{=T_1} + \eta_t\underbrace{\langle \nabla F(x_t)-w_t,\tilde{x}_{t+1}-x_t\rangle}_{=T_2} + \frac{L\eta_t^2}{2}\|\tilde{x}_{t+1}-x_t\|^2, \nonumber
 \end{align}
where the second equality is due to $x_{t+1}=x_t + \eta_t(\tilde{x}_{t+1}-x_t)$.

According to Assumption \ref{ass:7}, i.e., $A_t\succ \rho I_d$ for any $t\geq 1$,
the function $\phi_t(x)=\frac{1}{2}x^TA_t x$ is $\rho$-strongly convex, then we define a prox-function (a.k.a. Bregman distance) associated with $\phi_t(x)$ as in \cite{censor1992proximal,ghadimi2016mini},
\begin{align}
 V_t(x,x_t) = \phi_t(x) - \big[ \phi_t(x_t) + \langle\nabla \phi_t(x_t), x-x_t\rangle\big] = \frac{1}{2}(x-x_t)^TA_t(x-x_t).
\end{align}
By using the above Lemma \ref{lem:A1} to the problem $\tilde{x}_{t+1} = \arg\min_{x \in \mathcal{X}} \big\{\langle w_t, x\rangle + \frac{1}{2\gamma}(x-x_t)^T A_t (x-x_t)\big\}$ at
the line 5 of Algorithm \ref{alg:1} or \ref{alg:2}, we have
\begin{align}
  \langle w_t, \frac{1}{\gamma}(x_t - \tilde{x}_{t+1})\rangle \geq \rho\|\frac{1}{\gamma}(x_t - \tilde{x}_{t+1})\|^2.
\end{align}
Then we obtain
\begin{align} \label{eq:B3}
  T_1=\langle w_t, \tilde{x}_{t+1}-x_t\rangle \leq -\frac{\rho }{\gamma }\|\tilde{x}_{t+1}-x_t\|^2.
\end{align}

Next, consider the bound of the term $T_2$, we have
\begin{align} \label{eq:B4}
  T_2 & = \langle \nabla F(x_t)-w_t,\tilde{x}_{t+1}-x_t\rangle \nonumber \\
  & \leq \|\nabla F(x_t)-w_t\|\cdot\|\tilde{x}_{t+1}-x_t\| \nonumber \\
  & \leq \frac{\gamma}{\rho}\|\nabla F(x_t)-w_t\|^2+\frac{\rho}{4\gamma}\|\tilde{x}_{t+1}-x_t\|^2,
\end{align}
where the first inequality is due to the Cauchy-Schwarz inequality and the last is due to Young's inequality.
By combining the above inequalities (\ref{eq:B1}), (\ref{eq:B3}) with (\ref{eq:B4}),
we obtain
\begin{align} \label{eq:B5}
  F(x_{t+1}) &\leq F(x_t) + \eta_t\langle \nabla F(x_t)-w_t,\tilde{x}_{t+1}-x_t\rangle + \eta_t\langle w_t,\tilde{x}_{t+1}-x_t\rangle + \frac{L\eta_t^2}{2}\|\tilde{x}_{t+1}-x_t\|^2  \nonumber \\
  & \leq F(x_t) + \frac{\eta_t\gamma}{\rho}\|\nabla F(x_t)-w_t\|^2 + \frac{\rho\eta_t}{4\gamma}\|\tilde{x}_{t+1}-x_t\|^2  -\frac{\rho\eta_t}{\gamma}\|\tilde{x}_{t+1}-x_t\|^2 + \frac{L\eta_t^2}{2}\|\tilde{x}_{t+1}-x_t\|^2 \nonumber \\
  & = F(x_t) + \frac{\eta_t\gamma}{\rho}\|\nabla F(x_t)-w_t\|^2 -\frac{\rho\eta_t}{2\gamma}\|\tilde{x}_{t+1}-x_t\|^2  -\big(\frac{\rho\eta_t}{4\gamma}-\frac{L\eta_t^2}{2}\big)\|\tilde{x}_{t+1}-x_t\|^2 \nonumber \\
  & \leq F(x_t) + \frac{\eta_t\gamma}{\rho}\|\nabla F(x_t)-w_t\|^2 -\frac{\rho\eta_t}{2\gamma}\|\tilde{x}_{t+1}-x_t\|^2,
\end{align}
where the last inequality is due to $0< \gamma \leq \frac{\rho}{2L\eta_t}$.

\end{proof}

\begin{lemma} \label{lem:A5}
Suppose the sequence $\{x_t,y_t\}_{t=1}^T$ be generated from Algorithm \ref{alg:1} or \ref{alg:2}.
Under the above assumptions, given $0< \eta_t\leq 1$, $B_t=b_tI_p \ (b_u \geq b_t \geq b_l>0)$ for all $t\geq 1$,
and $0<\lambda\leq \frac{b_l}{6L_g}$, we have
\begin{align}
  \|y_{t+1} - y^*(x_{t+1})\|^2 &\leq (1-\frac{\eta_t\mu\lambda}{4b_t})\|y_t -y^*(x_t)\|^2 -\frac{3\eta_t}{4} \|\tilde{y}_{t+1}-y_t\|^2 \nonumber \\
     & \quad + \frac{25\eta_t\lambda}{6\mu b_t} \|\nabla_y g(x_t,y_t)-v_t\|^2 + \frac{25\kappa^2\eta_tb_t}{6\mu\lambda}\|\tilde{x}_{t+1} - x_t\|^2,
\end{align}
where $\kappa = L_g/\mu$.
\end{lemma}

\begin{proof}
According to Assumption \ref{ass:1}, i.e., the function $g(x,y)$ is $\mu$-strongly convex w.r.t $y$,
we have
\begin{align} \label{eq:F1}
 g(x_t,y) & \geq g(x_t,y_t) + \langle\nabla_y g(x_t,y_t), y-y_t\rangle + \frac{\mu}{2}\|y-y_t\|^2 \nonumber \\
 & = g(x_t,y_t) + \langle v_t, y-\tilde{y}_{t+1}\rangle + \langle\nabla_y g(x_t,y_t)-v_t, y-\tilde{y}_{t+1}\rangle \nonumber \\
 & \quad +\langle\nabla_y g(x_t,y_t), \tilde{y}_{t+1}-y_t\rangle + \frac{\mu}{2}\|y-y_t\|^2.
\end{align}
According to the Assumption \ref{ass:2}, i.e., the function $g(x,y)$ is $L_g$-smooth, we have
\begin{align} \label{eq:F2}
 g(x_t,\tilde{y}_{t+1}) \leq g(x_t,y_{t}) + \langle\nabla_y g(x_t,y_t), \tilde{y}_{t+1}-y_t\rangle + \frac{L_g}{2}\|\tilde{y}_{t+1}-y_t\|^2 .
\end{align}
Combining the about inequalities (\ref{eq:F1}) with (\ref{eq:F2}), we have
\begin{align} \label{eq:F3}
 g(x_t,y) & \geq g(x_t,\tilde{y}_{t+1}) + \langle v_t, y-\tilde{y}_{t+1}\rangle + \langle\nabla_y g(x_t,y_t)-v_t, y-\tilde{y}_{t+1}\rangle \nonumber \\
 & \quad + \frac{\mu}{2}\|y-y_t\|^2 - \frac{L_g}{2}\|\tilde{y}_{t+1}-y_t\|^2.
\end{align}

By the optimality condition of the problem $\tilde{y}_{t+1} = \arg\min_{y\in \mathcal{Y}}\big\{ \langle v_t, y\rangle + \frac{1}{2\lambda}(y-y_t)^TB_t(y-y_t) \big\}$
at the line 6 of Algorithm \ref{alg:1} or \ref{alg:2},
given $B_t=b_tI_p$, we have
\begin{align}
  \langle v_t + \frac{b_t}{\lambda}(\tilde{y}_{t+1} - y_t), y-\tilde{y}_{t+1} \rangle \geq 0, \quad \forall y\in \mathcal{Y}.
\end{align}
Then we obtain
\begin{align} \label{eq:F4}
  \langle v_t, y-\tilde{y}_{t+1}\rangle & \geq \frac{b_t}{\lambda}\langle \tilde{y}_{t+1}- y_t, \tilde{y}_{t+1}-y\rangle  \nonumber \\
  & = \frac{b_t}{\lambda}\|\tilde{y}_{t+1}- y_t\|^2 + \frac{b_t}{\lambda}\langle \tilde{y}_{t+1}- y_t, y_t- y\rangle.
\end{align}

By pugging the inequalities (\ref{eq:F4}) into (\ref{eq:F3}), we have
\begin{align}
 g(x_t,y) & \geq g(x_t,\tilde{y}_{t+1}) + \frac{b_t}{\lambda}\langle \tilde{y}_{t+1}- y_t, y_t- y\rangle + \frac{b_t}{\lambda}\|\tilde{y}_{t+1}- y_t\|^2 \nonumber \\
 & \quad + \langle\nabla_y g(x_t,y_t)-v_t, y-\tilde{y}_{t+1}\rangle +\frac{\mu}{2}\|y-y_t\|^2 -\frac{L_g}{2}\|\tilde{y}_{t+1}-y_t\|^2.
\end{align}
Let $y=y^*(x_t)$, then we have
\begin{align}
 g(x_t,y^*(x_t)) & \geq g(x_t,\tilde{y}_{t+1}) + \frac{b_t}{\lambda}\langle \tilde{y}_{t+1}- y_t, y_t - y^*(x_t)\rangle
  + (\frac{b_t}{\lambda}-\frac{L_g}{2})\|\tilde{y}_{t+1}- y_t\|^2 \nonumber \\
 & \quad + \langle\nabla_y g(x_t,y_t)-v_t, y^*(x_t)-\tilde{y}_{t+1}\rangle + \frac{\mu}{2}\|y^*(x_t)-y_t\|^2.
\end{align}
Due to the strongly-convexity of $g(\cdot,y)$ and $y^*(x_t) =\arg\min_{y\in \mathcal{Y}}g(x_t,y)$, we have $g(x_t,y^*(x_t)) \leq g(x_t,\tilde{y}_{t+1})$.
Thus, we obtain
\begin{align} \label{eq:F5}
 0 & \geq  \frac{b_t}{\lambda}\langle \tilde{y}_{t+1}- y_t, y_t - y^*(x_t)\rangle
 + \langle\nabla_y g(x_t,y_t)-v_t, y^*(x_t)-\tilde{y}_{t+1}\rangle \nonumber \\
 & \quad + (\frac{b_t}{\lambda}-\frac{L_g}{2})\|\tilde{y}_{t+1}- y_t\|^2 + \frac{\mu}{2}\|y^*(x_t)-y_t\|^2.
\end{align}

By $y_{t+1} = y_t + \eta_t(\tilde{y}_{t+1}-y_t) $, we have
\begin{align}
\|y_{t+1}-y^*(x_t)\|^2
& = \|y_t + \eta_t(\tilde{y}_{t+1}-y_t) -y^*(x_t)\|^2 \nonumber \\
& =  \|y_t -y^*(x_t)\|^2 + 2\eta_t\langle \tilde{y}_{t+1}-y_t, y_t -y^*(x_t)\rangle + \eta_t^2 \|\tilde{y}_{t+1}-y_t\|^2.
\end{align}
Then we obtain
\begin{align} \label{eq:F6}
 \langle \tilde{y}_{t+1}-y_t, y_t - y^*(x_t)\rangle = \frac{1}{2\eta_t}\|y_{t+1}-y^*(x_t)\|^2 - \frac{1}{2\eta_t}\|y_t -y^*(x_t)\|^2 - \frac{\eta_t}{2}\|\tilde{y}_{t+1}-y_t\|^2.
\end{align}
Consider the upper bound of the term $\langle\nabla_y g(x_t,y_t)-v_t, y^*(x_t)-\tilde{y}_{t+1}\rangle$, we have
\begin{align} \label{eq:F7}
 &\langle\nabla_y g(x_t,y_t)-v_t, y^*(x_t)-\tilde{y}_{t+1}\rangle \nonumber \\
 & = \langle\nabla_y g(x_t,y_t)-v_t, y^*(x_t)-y_t\rangle + \langle\nabla_y g(x_t,y_t)-v_t, y_t-\tilde{y}_{t+1}\rangle \nonumber \\
 & \geq -\frac{1}{\mu} \|\nabla_y g(x_t,y_t)-v_t\|^2 - \frac{\mu}{4}\|y^*(x_t)-y_t\|^2 - \frac{1}{\mu} \|\nabla_y g(x_t,y_t)-v_t\|^2 - \frac{\mu}{4}\|y_t-\tilde{y}_{t+1}\|^2 \nonumber \\
 & = -\frac{2}{\mu} \|\nabla_y g(x_t,y_t)-v_t\|^2 - \frac{\mu}{4}\|y^*(x_t)-y_t\|^2 - \frac{\mu}{4}\|y_t-\tilde{y}_{t+1}\|^2.
\end{align}

By plugging the inequalities (\ref{eq:F6}) and (\ref{eq:F7}) into (\ref{eq:F5}),
we obtain
\begin{align}
 & \frac{b_t}{2\eta_t\lambda}\|y_{t+1}-y^*(x_t)\|^2 \nonumber \\
 & \leq (\frac{b_t}{2\eta_t\lambda}-\frac{\mu}{4})\|y_t -y^*(x_t)\|^2 + ( \frac{b_t\eta_t}{2\lambda} + \frac{\mu}{4} + \frac{L_g}{2}-\frac{b_t}{\lambda}) \|\tilde{y}_{t+1}-y_t\|^2 + \frac{2}{\mu} \|\nabla_y g(x_t,y_t)-v_t\|^2 \nonumber \\
 & \leq ( \frac{b_t}{2\eta_t\lambda}-\frac{\mu}{4})\|y_t -y^*(x_t)\|^2 + (\frac{3L_g}{4} -\frac{b_t}{2\lambda}) \|\tilde{y}_{t+1}-y_t\|^2 + \frac{2}{\mu} \|\nabla_y g(x_t,y_t)-v_t\|^2 \nonumber \\
 & = ( \frac{b_t}{2\eta_t\lambda}-\frac{\mu}{4})\|y_t -y^*(x_t)\|^2 - \big( \frac{3b_t}{8\lambda} + \frac{b_t}{8\lambda} -\frac{3L_g}{4}\big) \|\tilde{y}_{t+1}-y_t\|^2 + \frac{2}{\mu} \|\nabla_y g(x_t,y_t)-v_t\|^2 \nonumber \\
 & \leq  (\frac{b_t}{2\eta_t\lambda}-\frac{\mu}{4})\|y_t -y^*(x_t)\|^2 - \frac{3b_t}{8\lambda} \|\tilde{y}_{t+1}-y_t\|^2 + \frac{2}{\mu}\|\nabla_y g(x_t,y_t)-v_t\|^2,
\end{align}
where the second inequality holds by $L_g \geq \mu$ and $0< \eta_t \leq 1$, and the last inequality is due to
$0< \lambda \leq \frac{b_l}{6L_g} \leq \frac{b_t}{6L_g}$.
It implies that
\begin{align} \label{eq:F8}
\|y_{t+1}-y^*(x_t)\|^2 \leq ( 1 -\frac{\eta_t\mu\lambda}{2b_t})\|y_t -y^*(x_t)\|^2 - \frac{3\eta_t}{4} \|\tilde{y}_{t+1}-y_t\|^2
+ \frac{4\eta_t\lambda}{\mu b_t}\|\nabla_y g(x_t,y_t)-v_t\|^2.
\end{align}

Next, we decompose the term $\|y_{t+1} - y^*(x_{t+1})\|^2$ as follows:
\begin{align} \label{eq:F9}
  \|y_{t+1} - y^*(x_{t+1})\|^2 & = \|y_{t+1} - y^*(x_t) + y^*(x_t) - y^*(x_{t+1})\|^2    \nonumber \\
  & = \|y_{t+1} - y^*(x_t)\|^2 + 2\langle y_{t+1} - y^*(x_t), y^*(x_t) - y^*(x_{t+1})\rangle  + \|y^*(x_t) - y^*(x_{t+1})\|^2  \nonumber \\
  & \leq (1+\frac{\eta_t\mu\lambda}{4b_t})\|y_{t+1} - y^*(x_t)\|^2  + (1+\frac{4b_t}{\eta_t\mu\lambda})\|y^*(x_t) - y^*(x_{t+1})\|^2 \nonumber \\
  & \leq (1+\frac{\eta_t\mu\lambda}{4b_t})\|y_{t+1} - y^*(x_t)\|^2  + (1+\frac{4b_t}{\eta_t\mu\lambda})\kappa^2\|x_t - x_{t+1}\|^2,
\end{align}
where the first inequality holds by Cauchy-Schwarz inequality and Young's inequality, and the second inequality is due to
Lemma \ref{lem:2}, and the last equality holds by $x_{t+1}=x_t + \eta_t(\tilde{x}_{t+1}-x_t)$.

By combining the above inequalities (\ref{eq:F8}) and (\ref{eq:F9}), we have
\begin{align}
 \|y_{t+1} - y^*(x_{t+1})\|^2 & \leq (1+\frac{\eta_t\mu\lambda}{4b_t})( 1-\frac{\eta_t\mu\lambda}{2b_t})\|y_t -y^*(x_t)\|^2
 - (1+\frac{\eta_t\mu\lambda}{4b_t})\frac{3\eta_t}{4} \|\tilde{y}_{t+1}-y_t\|^2     \nonumber \\
 & \quad + (1+\frac{\eta_t\mu\lambda}{4b_t})\frac{4\eta_t\lambda}{\mu b_t}\|\nabla_y g(x_t,y_t)-v_t\|^2
  + (1+\frac{4b_t}{\eta_t\mu\lambda})\kappa^2\|x_t - x_{t+1}\|^2.   \nonumber
\end{align}
Since $0 < \eta_t \leq 1$, $0< \lambda \leq \frac{b_l}{6L_g} \leq \frac{b_t}{6L_g}$ and $L_g\geq \mu$, we have $\lambda \leq \frac{b_t}{6L_g} \leq \frac{b_t}{6\mu}$
and $\eta_t\leq 1\leq \frac{b_t}{6\mu \lambda}$. Then we have
\begin{align}
  (1+\frac{\eta_t\mu\lambda}{4b_t})(1-\frac{\eta_t\mu\lambda}{2b_t})&= 1-\frac{\eta_t\mu\lambda}{2b_t} +\frac{\eta_t\mu\lambda}{4b_t}
  - \frac{\eta_t^2\mu^2\lambda^2}{8b^2_t} \leq 1-\frac{\eta_t\mu\lambda}{4b_t}, \nonumber  \\
 - (1+\frac{\eta_t\mu\lambda}{4b_t})\frac{3\eta_t}{4} &\leq -\frac{3\eta_t}{4}, \nonumber  \\
 (1+\frac{\eta_t\mu\lambda}{4b_t})\frac{4\eta_t\lambda}{\mu b_t} & \leq (1+\frac{1}{24})\frac{4\eta_t\lambda}{\mu b_t}
 = \frac{25\eta_t\lambda}{6\mu b_t}, \nonumber  \\
  (1+\frac{4b_t}{\eta_t\mu\lambda})\kappa^2 &
  \leq \frac{b_t\kappa^2}{6\eta_t\mu\lambda} +\frac{4b_t\kappa^2}{\eta_t\mu\lambda} = \frac{25b_t\kappa^2}{6\eta_t\mu\lambda}, \nonumber
\end{align}
where the second last inequality is due to $\frac{\eta_t\mu\lambda}{b_t}\leq \frac{1}{6}$ and the last inequality holds by
$\frac{b_t}{6\mu\lambda\eta_t} \geq 1$.
By using $x_{t+1}=x_t + \eta_t(\tilde{x}_{t+1}-x_t)$, then we have
\begin{align}
     \|y_{t+1} - y^*(x_{t+1})\|^2 &\leq (1-\frac{\eta_t\mu\lambda}{4b_t})\|y_t -y^*(x_t)\|^2 -\frac{3\eta_t}{4} \|\tilde{y}_{t+1}-y_t\|^2 \nonumber \\
     & \quad + \frac{25\eta_t\lambda}{6\mu b_t} \|\nabla_y g(x_t,y_t)-v_t\|^2 + \frac{25\kappa^2\eta_tb_t}{6\mu\lambda}\|\tilde{x}_{t+1} - x_t\|^2.
\end{align}

\end{proof}

\subsection{ Convergence Analysis of BiAdam Algorithm }
\label{Appendix:A1}
In the subsection, we provide the detail convergence analysis of BiAdam algorithm.
For notational simplicity, let $R_t = R(x_t,y_t)$ for all $t\geq 1$.

\begin{lemma} \label{lem:C3}
 Assume that the stochastic partial derivatives $v_{t+1}$, and $w_{t+1}$ be generated from Algorithm \ref{alg:1}, we have
 \begin{align}
 \mathbb{E}\| w_{t+1} - \bar{\nabla} f(x_{t+1},y_{t+1})-R_{t+1}\|^2 & \leq  (1-\beta_{t+1}) \mathbb{E}\|w_t - \bar{\nabla} f(x_t,y_t) -R_t\|^2  + \beta^2_{t+1}\sigma^2   \\
 & \quad + \frac{3L^2_0\eta^2_t}{\beta_{t+1}}\big( \|\tilde{x}_{t+1}-x_t\|^2 + \|\tilde{y}_{t+1}-y_t\|^2 \big) + \frac{3}{\beta_{t+1}}\big( \| R_t\|^2 + \|R_{t+1}\|^2 \big), \nonumber
 \end{align}
  \begin{align}
 \mathbb{E}\|v_{t+1}- \nabla_y g(x_{t+1},y_{t+1})\|^2
 & \leq (1-\alpha_{t+1}) \mathbb{E} \|v_t - \nabla_y g(x_t,y_t)\|^2 + \alpha_{t+1}^2\sigma^2  \\
 & \quad + 2L_g^2\eta_t^2/\alpha_{t+1}\big(\mathbb{E}\|\tilde{x}_{t+1} - x_t\|^2 + \mathbb{E}\|\tilde{y}_{t+1} - y_t\|^2 \big), \nonumber
 \end{align}
where $L^2_0 = 8\big(L^2_f+ \frac{L^2_{gxy}C^2_{fy}}{\mu^2} + \frac{L^2_{gyy} C^2_{gxy}C^2_{fy}}{\mu^4} +
 \frac{L^2_fC^2_{gxy}}{\mu^2}\big)$ and $R_t = \bar{\nabla}f(x_t,y_t) - \mathbb{E}_{\bar{\xi}}[\bar{\nabla}f(x_t,y_t;\bar{\xi})]$
 for all $t\geq 1$.
\end{lemma}

\begin{proof}
Without loss of generality, we only prove the term $\mathbb{E}\|w_{t+1} - \bar{\nabla} f(x_{t+1},y_{t+1})-R_{t+1}\|^2$.
The other term is similar for this term.
Since $w_{t+1} = \beta_{t+1}\bar{\nabla} f(x_{t+1},y_{t+1};\bar{\xi}_{t+1}) + (1-\beta_{t+1})w_t$, we have
\begin{align}
 &\mathbb{E}\|w_{t+1} - \bar{\nabla} f(x_{t+1},y_{t+1}) - R_{t+1}\|^2 \\
 & = \mathbb{E}\|\beta_{t+1}\bar{\nabla} f(x_{t+1},y_{t+1};\bar{\xi}_{t+1}) + (1-\beta_{t+1})w_t  - \bar{\nabla} f(x_{t+1},y_{t+1}) - R_{t+1}\|^2 \nonumber \\
 & = \mathbb{E}\|(1-\beta_{t+1})(w_t - \bar{\nabla} f(x_t,y_t) -R_t) + \beta_{t+1}\big( \bar{\nabla} f(x_{t+1},y_{t+1};\bar{\xi}_{t+1}) - \bar{\nabla} f(x_{t+1},y_{t+1}) - R_{t+1} \big) \nonumber \\
 & \quad + (1-\beta_{t+1})\big( \bar{\nabla} f(x_t,y_t) + R_t - (\bar{\nabla} f(x_{t+1},y_{t+1}) + R_{t+1}) \big)\|^2 \nonumber \\
 & = \beta_{t+1}^2 \mathbb{E}\| \bar{\nabla} f(x_{t+1},y_{t+1};\bar{\xi}_{t+1}) - \bar{\nabla} f(x_{t+1},y_{t+1}) - R_{t+1}\|^2 \nonumber \\
 & \quad + \mathbb{E}\|(1-\beta_{t+1})\big( w_t - \bar{\nabla} f(x_t,y_t) - R_t + \bar{\nabla} f(x_t,y_t) + R_t - (\bar{\nabla} f(x_{t+1},y_{t+1}) + R_{t+1}) \big)\|^2 \nonumber \\
 & \leq \beta_{t+1}^2 \mathbb{E}\| \bar{\nabla} f(x_{t+1},y_{t+1};\bar{\xi}_{t+1}) - \bar{\nabla} f(x_{t+1},y_{t+1}) - R_{t+1}\|^2 + (1-\beta_{t+1})^2(1+\beta_{t+1}) \mathbb{E}\|w_t - \bar{\nabla} f(x_t,y_t) -R_t\|^2  \nonumber \\
 & \quad + (1-\beta_{t+1})^2(1+\frac{1}{\beta_{t+1}})\mathbb{E}\| \bar{\nabla} f(x_t,y_t) + R_t - (\bar{\nabla} f(x_{t+1},y_{t+1}) + R_{t+1}) \|^2  \nonumber \\
 & \leq (1-\beta_{t+1}) \mathbb{E}\|w_t - \bar{\nabla} f(x_t,y_t) -R_t\|^2  + \beta^2_{t+1}\sigma^2
  + \frac{1}{\beta_{t+1}}\| \bar{\nabla} f(x_t,y_t) + R_t - (\bar{\nabla} f(x_{t+1},y_{t+1}) + R_{t+1})\|^2 \nonumber \\
 & \leq (1-\beta_{t+1}) \mathbb{E}\|w_t - \bar{\nabla} f(x_t,y_t) -R_t\|^2  + \beta^2_{t+1}\sigma^2 \nonumber \\
 & \quad + \frac{3}{\beta_{t+1}}\| \bar{\nabla} f(x_{t+1},y_{t+1})-\bar{\nabla} f(x_t,y_t) \|^2 + \frac{3}{\beta_{t+1}}\big(\| R_t\|^2 + \|R_{t+1})\|^2\big) \nonumber \\
 & \leq (1-\beta_{t+1}) \mathbb{E}\|w_t - \bar{\nabla} f(x_t,y_t) -R_t\|^2  + \beta^2_{t+1}\sigma^2 + \frac{3}{\beta_{t+1}}\big( \| R_t\|^2 + \|R_{t+1}\|^2 \big) \nonumber \\
 & \quad + \frac{3L^2_0\eta^2_t}{\beta_{t+1}}\big( \|\tilde{x}_{t+1}-x_t\|^2 + \|\tilde{y}_{t+1}-y_t\|^2 \big),
\end{align}
where the third equality is due to $\mathbb{E}_{\bar{\xi}_{t+1}}[\bar{\nabla} f(x_{t+1},y_{t+1};\bar{\xi}_{t+1})] = \bar{\nabla} f(x_{t+1},y_{t+1})+R_{t+1}$; the second last inequality holds by $0\leq \beta_{t+1} \leq 1$ such that  $(1-\beta_{t+1})^2(1+\beta_{t+1})=1-\beta_{t+1}-\beta_{t+1}^2+
  \beta_{t+1}^3\leq 1-\beta_{t+1}$ and $(1-\beta_{t+1})^2(1+\frac{1}{\beta_{t+1}}) \leq (1-\beta_{t+1})(1+\frac{1}{\beta_{t+1}})  =-\beta_{t+1}+\frac{1}{\beta_{t+1}}\leq \frac{1}{\beta_{t+1}}$,
  and the last inequality holds by the above Lemma \ref{lem:A3} and $x_{t+1}=x_t+\eta_t(\tilde{x}_{t+1}-x_t)$,  $y_{t+1}=y_t+\eta_t(\tilde{y}_{t+1}-y_t)$.
\end{proof}

\begin{theorem} \label{th:A1}
 (Restatement of Theorem 1)
 Under the above Assumptions (\ref{ass:1}, \ref{ass:2}, \ref{ass:4}, \ref{ass:6}, \ref{ass:7}), in the Algorithm \ref{alg:1}, given $\mathcal{X}\subset\mathbb{R}^{d}$, $\eta_t=\frac{k}{(m+t)^{1/2}}$ for all $t\geq 0$, $\alpha_{t+1}=c_1\eta_t$, $\beta_{t+1}=c_2\eta_t$, $m\geq \max\big(k^2, (c_1k)^2,(c_2k)^2\big)$, $k>0$, $\frac{125L^2_0}{6\mu^2} \leq c_1 \leq \frac{m^{1/2}}{k}$, $\frac{9}{2} \leq c_2 \leq \frac{m^{1/2}}{k}$, $0<\lambda\leq \min\big(\frac{15b_lL^2_0}{4L^2_1\mu}, \frac{b_l}{6L_g}\big)$, $0< \gamma \leq \min\big(\frac{\sqrt{6}\lambda\mu\rho}{\sqrt{6L^2_1\lambda^2\mu^2 + 125b^2_uL^2_0\kappa^2}}, \frac{m^{1/2}\rho}{4Lk}\big)$ and $K=\frac{L_g}{\mu}\log(C_{gxy}C_{fy}T/\mu)$, we have
\begin{align}
 \frac{1}{T} \sum_{t=1}^T \mathbb{E}||\mathcal{G}_\mathcal{X}(x_t,\nabla F(x_t),\gamma)|| \leq \frac{1}{T} \sum_{t=1}^T \mathbb{E}[\mathcal{M}_t] \leq \frac{2\sqrt{3G}m^{1/4}}{T^{1/2}} + \frac{2\sqrt{3G}}{T^{1/4}} + \frac{\sqrt{2}}{T},
\end{align}
where $G = \frac{F(x_1) - F^*}{\rho k\gamma} + \frac{5b_1L^2_0\Delta_0}{\rho^2 k\lambda\mu} + \frac{2\sigma^2}{\rho^2 k} + \frac{2m\sigma^2}{\rho^2 k}\ln(m+T) + \frac{4(m+T)}{9\rho^2 kT^2} + \frac{8k}{\rho^2T}$, $\Delta_0=\|y_1-y^*(x_1)\|^2$ and $L^2_1= \frac{12L^2_g\mu^2}{125L^2_0} + \frac{2L^2_0}{3}$.
\end{theorem}

\begin{proof}
Since $\eta_t=\frac{k}{(m+t)^{1/2}}$ on $t$ is decreasing and $m\geq k^2$, we have $\eta_t \leq \eta_0 = \frac{k}{m^{1/2}} \leq 1$ and $\gamma \leq \frac{m^{1/2}\rho}{4Lk}\leq \frac{\rho}{2L\eta_0} \leq \frac{\rho}{2L\eta_t}$ for any $t\geq 0$.
 Due to $0 < \eta_t \leq 1$ and $m\geq (c_1k)^2$, we have $\alpha_{t+1} = c_1\eta_t \leq \frac{c_1k}{m^{1/2}}\leq 1$.
 Similarly, due to $m\geq (c_2k)^2$, we have $\beta_{t+1}\leq 1$. At the same time, we have $c_1,c_2\leq \frac{m^{1/2}}{k}$.
According to Lemma \ref{lem:C3}, we have
 \begin{align} \label{eq:S1}
  & \mathbb{E} \|v_{t+1} - \nabla_y g(x_{t+1},y_{t+1})\|^2 - \mathbb{E} \|v_t - \nabla_y g(x_t,y_t)\|^2  \\
  & \leq -c_1 \eta_t\mathbb{E} \|\nabla_y g(x_t,y_t) - v_t\|^2 + 2L^2_g\eta_t/c_1\big(\|\tilde{x}_{t+1}-x_t\|^2 + \|\tilde{y}_{t+1}-y_t\|^2\big) + c_1^2\eta_t^2\sigma^2 \nonumber \\
  & \leq -\frac{125L^2_0}{6\mu^2} \mathbb{E} \|\nabla_y g(x_t,y_t) -v_t\|^2 + \frac{12L^2_g\mu^2\eta_t}{125L^2_0}\big(\|\tilde{x}_{t+1}-x_t\|^2 + \|\tilde{y}_{t+1}-y_t\|^2\big) + \frac{m\eta^2_t\sigma^2}{k^2}, \nonumber
 \end{align}
 where the above equality holds by $\alpha_{t+1}=c_1\eta_t$, and the last inequality is due to $\frac{125L^2_0}{6\mu^2} \leq c_1 \leq \frac{m^{1/2}}{k}$.
 Similarly, we have
 \begin{align} \label{eq:S2}
 & \mathbb{E}\| w_{t+1} - \bar{\nabla} f(x_{t+1},y_{t+1})-R_{t+1}\|^2 - \mathbb{E}\|w_t - \bar{\nabla} f(x_t,y_t) -R_t\|^2 \\
 & \leq -\beta_{t+1}\mathbb{E}\|w_t - \bar{\nabla} f(x_t,y_t) -R_t\|^2 + \frac{3L^2_0\eta^2_t}{\beta_{t+1}}\big( \|\tilde{x}_{t+1}-x_t\|^2 + \|\tilde{y}_{t+1}-y_t\|^2 \big) \nonumber \\
 & \quad  + \frac{3}{\beta_{t+1}}\big( \| R_t\|^2 + \|R_{t+1}\|^2 \big) + \beta^2_{t+1}\sigma^2 \nonumber \\
 & \leq -\frac{9\eta_t}{2}\mathbb{E}\|w_t - \bar{\nabla} f(x_t,y_t) -R_t\|^2 + \frac{2L^2_0\eta_t}{3}\big( \|\tilde{x}_{t+1}-x_t\|^2 + \|\tilde{y}_{t+1}-y_t\|^2 \big) \nonumber \\
 & \quad  + \frac{2}{3\eta_t}\big( \| R_t\|^2 + \|R_{t+1}\|^2 \big) + \frac{m\eta^2_t\sigma^2}{k^2}, \nonumber
 \end{align}
where the last inequality holds by $\beta_{t+1}=c_2\eta_t$ and $ \frac{9}{2}\leq c_2 \leq \frac{m^{1/2}}{k}$.

According to Lemmas \ref{lem:A2} and \ref{lem:A4}, we have
\begin{align} \label{eq:S3}
  & F(x_{t+1}) - F(x_t) \\
  & \leq \frac{2\eta_t\gamma}{\rho}\|w_t-\bar{\nabla} f(x_t,y_t)\|^2 + \frac{L^2_0\eta_t\gamma}{\rho}\|y^*(x_t)-y_t\|^2 -\frac{\rho\eta_t}{2\gamma}\|\tilde{x}_{t+1}-x_t\|^2 \nonumber \\
  & \leq \frac{4\eta_t\gamma}{\rho}\|w_t-\bar{\nabla} f(x_t,y_t)-R_t\|^2 + \frac{4\eta_t\gamma}{\rho}\|R_t\|^2 + \frac{L^2_0\eta_t\gamma}{\rho}\|y^*(x_t)-y_t\|^2 -\frac{\rho\eta_t}{2\gamma}\|\tilde{x}_{t+1}-x_t\|^2. \nonumber
\end{align}
According to Lemma \ref{lem:A5}, we have
\begin{align} \label{eq:S4}
  & \|y_{t+1}-y^*(x_{t+1})\|^2 - \|y_t-y^*(x_t)\|^2 \\
  & \leq -\frac{\eta_t\mu\lambda}{4b_t}\|y_t -y^*(x_t)\|^2 -\frac{3\eta_t}{4} \|\tilde{y}_{t+1}-y_t\|^2
  + \frac{25\eta_t\lambda}{6\mu b_t} \|\nabla_y g(x_t,y_t)-v_t\|^2 + \frac{25\kappa^2\eta_tb_t}{6\mu\lambda}\|\tilde{x}_{t+1} - x_t\|^2. \nonumber
\end{align}

 Next, we define a \emph{Lyapunov} function (i.e., potential function), for any $t\geq 1$,
 \begin{align}
 \Gamma_t & = \mathbb{E}\big [F(x_t) + \frac{5b_tL^2_0\gamma}{\lambda\mu\rho}\|y_t-y^*(x_t)\|^2 + \frac{\gamma}{\rho} \big( \|v_t - \nabla_y g(x_t,y_t)\|^2  + \|w_t - \bar{\nabla} f(x_t,y_t) -R_t\|^2 \big) \big]. \nonumber
 \end{align}
 For notational simplicity, let $L^2_1= \frac{12L^2_g\mu^2}{125L^2_0} + \frac{2L^2_0}{3}$.
 Then we have
 \begin{align}
 & \Gamma_{t+1} - \Gamma_t \nonumber \\
 & = F(x_{t+1}) - F(x_t) + \frac{5b_tL^2_0\gamma}{\lambda\mu\rho}
 \big( \|y_{t+1}-y^*(x_{t+1})\|^2 - \|y_t-y^*(x_t)\|^2 \big)
 + \frac{\gamma}{\rho} \big(\|v_{t+1} - \nabla_y g(x_{t+1},y_{t+1})\|^2 \nonumber \\
 & \quad - \|v_t - \nabla_y g(x_t,y_t)\|^2 + \|w_{t+1} - \bar{\nabla} f(x_{t+1},y_{t+1}) -R_{t+1}\|^2 -\|w_t - \bar{\nabla} f(x_t,y_t) -R_t\|^2 \big) \nonumber \\
 & \leq \frac{L^2_0\eta_t\gamma}{\rho}\|y^*(x_t)-y_t\|^2 + \frac{4\eta_t\gamma}{\rho}\|w_t-\bar{\nabla} f(x_t,y_t)-R_t\|^2 + \frac{4\eta_t\gamma}{\rho}\|R_t\|^2-\frac{\rho\eta_t}{2\gamma}\|\tilde{x}_{t+1}-x_t\|^2 \nonumber \\
 & \quad + \frac{5b_t L^2_0\gamma}{\lambda\mu\rho} \bigg( -\frac{\eta_t\mu\lambda}{4b_t}\|y_t -y^*(x_t)\|^2 -\frac{3\eta_t}{4} \|\tilde{y}_{t+1}-y_t\|^2
  + \frac{25\eta_t\lambda}{6\mu b_t} \|\nabla_y g(x_t,y_t)-v_t\|^2 + \frac{25\kappa^2\eta_tb_t}{6\mu\lambda}\|\tilde{x}_{t+1} - x_t\|^2 \bigg)  \nonumber \\
 & \quad + \frac{\gamma}{\rho} \bigg( -\frac{125L^2_0}{6\mu^2} \mathbb{E} \|\nabla_y g(x_t,y_t) -v_t\|^2 + \frac{12L^2_g\mu^2\eta_t}{125L^2_0}\big(\|\tilde{x}_{t+1}-x_t\|^2 + \|\tilde{y}_{t+1}-y_t\|^2\big) + \frac{m\eta^2_t\sigma^2}{k^2} \nonumber \\
 & \quad  -\frac{9\eta_t}{2}\mathbb{E}\|w_t - \bar{\nabla} f(x_t,y_t) -R_t\|^2 + \frac{2L^2_0\eta_t}{3}\big( \|\tilde{x}_{t+1}-x_t\|^2 + \|\tilde{y}_{t+1}-y_t\|^2 \big) + \frac{2}{3\eta_t}\big( \| R_t\|^2 + \|R_{t+1}\|^2 \big) + \frac{m\eta^2_t\sigma^2}{k^2} \bigg)\nonumber \\
 & = - \frac{\gamma\eta_t}{4\rho} \bigg( L^2_0 \|y_t - y^*(x_t)\|^2 + 2\mathbb{E}\|w_t - \bar{\nabla} f(x_t,y_t) -R_t\|^2 \bigg) - \big( \frac{\rho}{2\gamma} - \frac{L^2_1\gamma}{\rho} - \frac{125b^2_t L^2_0\kappa^2\gamma}{6\mu^2\lambda^2\rho}\big)\eta_t\|\tilde{x}_{t+1}-x_t\|^2
 \nonumber \\
 & \quad - \big( \frac{15b_t L^2_0\gamma}{4\lambda\mu\rho} -
 \frac{L^2_1\gamma}{\rho}\big)\eta_t\|\tilde{y}_{t+1}-y_t\|^2 + + \frac{2m\gamma\sigma^2}{k^2\rho}\eta^2_t
  + \frac{2\gamma}{3\rho\eta_t}\big( \| R_t\|^2 + \|R_{t+1}\|^2 \big) + \frac{4\eta_t\gamma}{\rho}\|R_t\|^2 \nonumber \\
 & \leq - \frac{\gamma\eta_t}{4\rho} \bigg( L^2_0 \|y_t - y^*(x_t)\|^2 + 2\mathbb{E}\|w_t - \bar{\nabla} f(x_t,y_t) -R_t\|^2 \bigg) - \frac{\rho\eta_t}{4\gamma}\|\tilde{x}_{t+1}-x_t\|^2 + \frac{2m\gamma\sigma^2}{k^2\rho}\eta^2_t \nonumber \\
 & \quad  + \frac{2\gamma}{3\rho\eta_t}\big( \| R_t\|^2 + \|R_{t+1}\|^2 \big) + \frac{4\eta_t\gamma}{\rho}\|R_t\|^2,
 \end{align}
where the first inequality holds by the above inequalities (\ref{eq:S1}), (\ref{eq:S2}), (\ref{eq:S3}) and (\ref{eq:S4});
the last inequality is due to $0< \gamma \leq \frac{\sqrt{6}\lambda\mu\rho}{\sqrt{6L^2_1\lambda^2\mu^2 + 125b^2_uL^2_0\kappa^2}} \leq \frac{\sqrt{6}\lambda\mu\rho}{\sqrt{6L^2_1\lambda^2\mu^2 + 125b^2_tL^2_0\kappa^2}}$, $0< \lambda \leq \frac{15b_lL^2_0}{4L^2_1\mu} \leq \frac{15b_tL^2_0}{4L^2_1\mu}$ for all $t\geq 1$.

Let $\Phi_t = L^2_0 \|y_t - y^*(x_t)\|^2 + 2\|w_t - \bar{\nabla} f(x_t,y_t) -R_t\|^2$, we have
\begin{align} \label{eq:S5}
 \frac{\gamma\eta_t}{4\rho} \Phi_t + \frac{\rho\eta_t}{4\gamma}\|\tilde{x}_{t+1}-x_t\|^2 & \leq \Gamma_t - \Gamma_{t+1}+ \frac{2m\gamma\sigma^2}{k^2\rho}\eta^2_t + \frac{2\gamma}{3\rho\eta_t}\big(\| R_t\|^2 + \|R_{t+1}\|^2 \big) + \frac{4\gamma\eta_t}{\rho}\|R_t\|^2.
\end{align}
Taking average over $t=1,2,\cdots,T$ on both sides of (\ref{eq:S5}), we have
\begin{align}
  \frac{1}{T} \sum_{t=1}^T \mathbb{E} \big[ \frac{\eta_t}{4}\Phi_t + \frac{\rho^2\eta_t}{4\gamma^2}\|\tilde{x}_{t+1}-x_t\|^2 \big]
  & \leq  \sum_{t=1}^T \frac{\rho(\Gamma_t - \Gamma_{t+1})}{T\gamma} + \frac{1}{T}\sum_{t=1}^T\frac{2m\sigma^2}{k^2}\eta^2_t \nonumber \\
  & \quad + \frac{1}{T}\sum_{t=1}^T \bigg( \frac{2}{3\eta_t}\big(\|R_t\|^2 + \|R_{t+1}\|^2 \big) + 4\eta_t\|R_t\|^2 \bigg). \nonumber
\end{align}
Given $x_1\in \mathcal{X}$ and $y_1\in \mathcal{Y}$, let $\Delta_0 = \|y_1-y^*(x_1)\|^2$, we have
\begin{align} \label{eq:S6}
 \Gamma_1 &= \mathbb{E}\big [F(x_t) + \frac{5b_1L^2_0\gamma}{\lambda\mu\rho}\|y_1-y^*(x_1)\|^2 + \frac{\gamma}{\rho} \big( \|v_1 - \nabla_y g(x_1,y_1)\|^2  + \|w_1 - \bar{\nabla} f(x_1,y_1) -R_1\|^2 \big) \big] \nonumber \\
 & \leq F(x_1) + \frac{5b_1L^2_0\gamma \Delta_0}{\lambda\mu\rho} + \frac{2\gamma\sigma^2}{\rho},
\end{align}
where the last inequality holds by Assumption \ref{ass:2}.
Since $\eta_t$ is decreasing on $t$, i.e., $\eta_T^{-1} \geq \eta_t^{-1}$ for any $0\leq t\leq T$, we have
 \begin{align}
 & \frac{1}{T} \sum_{t=1}^T \mathbb{E}\Big( \frac{\Phi_t}{4} + \frac{\rho^2}{4\gamma^2}\|\tilde{x}_{t+1}-x_t\|^2 \Big) \nonumber \\
 & \leq  \frac{\rho}{T\gamma\eta_T} \sum_{t=1}^T\big(\Gamma_t - \Gamma_{t+1}\big) + \frac{1}{T\eta_T}\sum_{t=1}^T\frac{2m\sigma^2}{k^2}\eta^2_t + \frac{1}{T}\sum_{t=1}^T \bigg( \frac{2}{3\eta_t}\big(\|R_t\|^2 + \|R_{t+1}\|^2 \big) + 4\eta_t\|R_t\|^2 \bigg) \nonumber \\
 & \leq \frac{\rho}{T\gamma\eta_T} \Big( F(x_1) + \frac{5b_1L^2_0\gamma \Delta_0}{\lambda\mu\rho} + \frac{2\gamma\sigma^2}{\rho} - F^* \Big) + \frac{1}{T\eta_T}\sum_{t=1}^T\frac{2m\sigma^2}{k^2}\eta^2_t + \frac{2}{3T^3}\sum_{t=1}^T\frac{1}{\eta_t} + \frac{4}{T^2}\sum_{t=1}^T\eta_t \nonumber \\
 & \leq \frac{\rho(F(x_1) - F^*)}{T\gamma\eta_T} +  \frac{5b_1L^2_0\Delta_0}{T\eta_T\lambda\mu} + \frac{2\sigma^2}{T\eta_T} + \frac{2m\sigma^2}{T\eta_T k^2}\int^T_1\frac{k^2}{m+t}dt + \frac{2}{3T^3}\int_{1}^T\frac{(m+t)^{1/2}}{k}dt \nonumber \\
 & \quad + \frac{4}{T^2}\int_{1}^T\frac{k}{(m+t)^{1/2}} dt \nonumber \\
 & \leq \frac{\rho(F(x_1) - F^*)}{T\gamma\eta_T} + \frac{5b_1L^2_0\Delta_0}{T\eta_T\lambda\mu} + \frac{2\sigma^2}{T\eta_T} + \frac{2m\sigma^2}{T\eta_T }\ln(m+T) + \frac{4}{9kT^3}(m+T)^{3/2} + \frac{8k}{T^2}(m+T)^{1/2} \nonumber \\
 & = \bigg( \frac{\rho(F(x_1) - F^*)}{k\gamma} + \frac{5b_1L^2_0\Delta_0}{k\lambda\mu} + \frac{2\sigma^2}{k} + \frac{2m\sigma^2}{k}\ln(m+T) + \frac{4(m+T)}{9kT^2} + \frac{8k}{T}\bigg)\frac{(m+T)^{1/2}}{T}, \nonumber
\end{align}
where the second inequality holds by the above inequality (\ref{eq:S6}) and $\|R_t\|\leq \frac{1}{T}$ for all $t\geq 1$ by choosing $K=\frac{L_g}{\mu}\log(C_{gxy}C_{fy}T/\mu)$ in Algorithm \ref{alg:1}.
Let $G = \frac{F(x_1) - F^*}{\rho k\gamma} + \frac{5b_1L^2_0\Delta_0}{\rho^2 k\lambda\mu} + \frac{2\sigma^2}{\rho^2k} + \frac{2m\sigma^2}{\rho^2 k}\ln(m+T) + \frac{4(m+T)}{9\rho^2kT^2} + \frac{8k}{\rho^2 T}$,
we have
\begin{align} \label{eq:S7}
 \frac{1}{T} \sum_{t=1}^T \mathbb{E}\big[ \frac{\Phi_t}{4\rho^2} + \frac{1}{4\gamma^2}\|\tilde{x}_{t+1}-x_t\|^2 \big] \leq \frac{G}{T}(m+T)^{1/2}.
\end{align}

According to the Jensen's inequality, we have
\begin{align}
& \frac{1}{T} \sum_{t=1}^T \mathbb{E}\Big[ \frac{1}{2\gamma}\|\tilde{x}_{t+1}-x_t\| + \frac{1}{2\rho} \big( L_0 \|y_t - y^*(x_t)\| + \sqrt{2}\|w_t - \bar{\nabla} f(x_t,y_t) -R_t\|  \big) \Big] \nonumber \\
& \leq \bigg( \frac{3}{T} \sum_{t=1}^T \big( \frac{1}{4\gamma^2}\|\tilde{x}_{t+1}-x_t\|^2  + \frac{L^2_0}{4\rho^2}\|y_t - y^*(x_t)\|^2 + \frac{2}{4\rho^2}\mathbb{E}\|w_t - \bar{\nabla} f(x_t,y_t) -R_t\|^2 \big)\bigg)^{1/2} \nonumber \\
& = \bigg( \frac{3}{T} \sum_{t=1}^T \big( \frac{\Phi_t}{4\rho^2} + \frac{1}{4\gamma^2}\|\tilde{x}_{t+1}-x_t\|^2 \big)\bigg)^{1/2} \nonumber \\
& \leq \frac{\sqrt{3G}}{T^{1/2}}(m+T)^{1/4} \leq \frac{\sqrt{3G}m^{1/4}}{T^{1/2}} + \frac{\sqrt{3G}}{T^{1/4}},
\end{align}
where the last inequality is due to $(a+b)^{1/4} \leq a^{1/4} + b^{1/4}$ for all $a,b>0$.
Thus we have
\begin{align} \label{eq:S8}
& \frac{1}{T} \sum_{t=1}^T \mathbb{E}\big[ \frac{1}{\gamma}\|\tilde{x}_{t+1}-x_t\| + \frac{1}{\rho} \big(L_0 \|y_t - y^*(x_t)\| + \sqrt{2}\|w_t - \bar{\nabla} f(x_t,y_t) -R_t\|\big) \big] \nonumber \\
& \leq \frac{2\sqrt{3G}m^{1/4}}{T^{1/2}} + \frac{2\sqrt{3G}}{T^{1/4}}.
\end{align}

Since $\|w_t - \bar{\nabla} f(x_t,y_t) -R_t\|\geq \|w_t - \bar{\nabla} f(x_t,y_t)\| - \|R_t\|$, by the above inequality (\ref{eq:S8}),
we can obtain
\begin{align}
\frac{1}{T} \sum_{t=1}^T \mathbb{E}[\mathcal{M}_t] & = \frac{1}{T} \sum_{t=1}^T \mathbb{E}\big[ \frac{1}{\gamma}\|\tilde{x}_{t+1}-x_t\|
+ \frac{1}{\rho} \big(L_0 \|y_t - y^*(x_t)\| + \sqrt{2}\|w_t - \bar{\nabla} f(x_t,y_t)\|\big) \big] \nonumber \\
& \leq \frac{2\sqrt{3G}m^{1/4}}{T^{1/2}} + \frac{2\sqrt{3G}}{T^{1/4}} + \frac{\sqrt{2}}{T} \sum_{t=1}^T \mathbb{E}\|R_t\| \nonumber \\
& = \frac{2\sqrt{3G}m^{1/4}}{T^{1/2}} + \frac{2\sqrt{3G}}{T^{1/4}} + \frac{\sqrt{2}}{T},
\end{align}
where the last inequality is due to $\|R_t\|\leq \frac{1}{T}$ for all $t\geq 1$ by choosing $K=\frac{L_g}{\mu}\log(C_{gxy}C_{fy}T/\mu)$
in Algorithm \ref{alg:1}.
According to the above inequality (\ref{eq:19}), we have
\begin{align}
\frac{1}{T} \sum_{t=1}^T \mathbb{E}||\mathcal{G}_\mathcal{X}(x_t,\nabla F(x_t),\gamma)|| \leq \frac{1}{T} \sum_{t=1}^T \mathbb{E}[\mathcal{M}_t] \leq \frac{2\sqrt{3G}m^{1/4}}{T^{1/2}} + \frac{2\sqrt{3G}}{T^{1/4}} + \frac{\sqrt{2}}{T}.
\end{align}

\end{proof}

\begin{theorem} \label{th:A2}
 (Restatement of Theorem 2)
 Under the above Assumptions (\ref{ass:1}, \ref{ass:2}, \ref{ass:4}, \ref{ass:6}, \ref{ass:7}), in the Algorithm \ref{alg:1}, given $\mathcal{X}=\mathbb{R}^{d}$, $\eta_t=\frac{k}{(m+t)^{1/2}}$ for all $t\geq 0$, $\alpha_{t+1}=c_1\eta_t$, $\beta_{t+1}=c_2\eta_t$, $m\geq \max\big(k^2, (c_1k)^2,(c_2k)^2\big)$, $k>0$, $\frac{125L^2_0}{6\mu^2} \leq c_1 \leq \frac{m^{1/2}}{k}$, $\frac{9}{2} \leq c_2 \leq \frac{m^{1/2}}{k}$, $0<\lambda\leq \min\big(\frac{15b_l L^2_0}{4L^2_1\mu}, \frac{b_l}{6L_g}\big)$, $0< \gamma \leq \min\big(\frac{\sqrt{6}\lambda\mu\rho}{\sqrt{6L^2_1\lambda^2\mu^2 + 125b^2_u L^2_0\kappa^2}}, \frac{m^{1/2}\rho}{4Lk}\big)$ and $K=\frac{L_g}{\mu}\log(C_{gxy}C_{fy}T/\mu)$, we have
\begin{align}
 \frac{1}{T}\sum_{t=1}^T\mathbb{E}\|\nabla F(x_t)\| \leq \frac{\sqrt{\frac{1}{T}\sum_{t=1}^T\mathbb{E}\|A_t\|^2}}{\rho}\Big( \frac{2\sqrt{6G'm}}{T^{1/2}} + \frac{2\sqrt{6G'}}{T^{1/4}} + \frac{2\sqrt{3}}{T}\Big),
\end{align}
where $G' = \frac{\rho(F(x_1) - F^*)}{k\gamma} + \frac{5b_u L^2_0\Delta_0}{k\lambda\mu} + \frac{2\sigma^2}{k} + \frac{2m\sigma^2}{k}\ln(m+T) + \frac{4(m+T)}{9kT^2} + \frac{8k}{T}$.
\end{theorem}

\begin{proof}
According to Lemma \ref{lem:A2}, we have
\begin{align}
\mathcal{M}_t & =\frac{1}{\gamma}\|x_t - \tilde{x}_{t+1}\|
+ \frac{1}{\rho}\big( L_0 \|y^*(x_t) - y_t\|  + \sqrt{2}\|\bar{\nabla} f(x_t,y_t)- w_t\|\big) \nonumber \\
& \geq \frac{1}{\gamma}\|x_t - \tilde{x}_{t+1}\|
+ \frac{1}{\rho}\|\nabla  F(x_t) - w_t\| \nonumber \\
& \mathop{=}^{(i)} \|A_t^{-1}w_t\|
+ \frac{1}{\rho}\|\nabla  F(x_t)-w_t\| \nonumber \\
& = \frac{1}{\|A_t\|}\|A_t\|\|A_t^{-1}w_t\|
+ \frac{1}{\rho}\|\nabla  F(x_t)-w_t\| \nonumber \\
& \geq \frac{1}{\|A_t\|}\|w_t\|
+ \frac{1}{\rho}\|\nabla  F(x_t)-w_t\| \nonumber \\
& \mathop{\geq}^{(ii)} \frac{1}{\|A_t\|}\|w_t\|
+ \frac{1}{\|A_t\|}\|\nabla  F(x_t)-w_t\| \nonumber \\
& \geq \frac{1}{\|A_t\|}\|\nabla  F(x_t)\|,
\end{align}
where the equality $(i)$ holds by $\tilde{x}_{t+1} = x_t - \gamma A^{-1}_tw_t$ that can be easily obtained from the step 5 of Algorithm \ref{alg:1} when $\mathcal{X}=\mathbb{R}^{d}$, and the inequality $(ii)$ holds by $\|A_t\| \geq \rho$ for all $t\geq1$ due to Assumption \ref{ass:7}. Then we have
\begin{align}
 \|\nabla  F(x_t)\| \leq \|A_t\|\mathcal{M}_t.
\end{align}
According to Cauchy-Schwarz inequality, we have
\begin{align} \label{eq:A77}
\frac{1}{T}\sum_{t=1}^T\mathbb{E}\|\nabla F(x_t)\| \leq \frac{1}{T}\sum_{t=1}^T\mathbb{E}\big[\mathcal{M}_t \|A_t\|\big] \leq \sqrt{\frac{1}{T}\sum_{t=1}^T\mathbb{E}[\mathcal{M}_t^2]} \sqrt{\frac{1}{T}\sum_{t=1}^T\mathbb{E}\|A_t\|^2}.
\end{align}
Then we have
\begin{align} \label{eq:A78}
 \frac{1}{T} \sum_{t=1}^T \mathbb{E}[\mathcal{M}_t^2] & \leq \frac{1}{T} \sum_{t=1}^T \mathbb{E}\big[ \frac{3L^2_0 \|y_t - y^*(x_t)\|^2}{\rho^2} + \frac{6\|w_t - \bar{\nabla} f(x_t,y_t)\|^2}{\rho^2} + \frac{3}{\gamma^2}\|\tilde{x}_{t+1}-x_t\|^2 \big] \nonumber \\
 & \leq \frac{1}{T} \sum_{t=1}^T \mathbb{E}\big[ \frac{3L^2_0 \|y_t - y^*(x_t)\|^2}{\rho^2} + \frac{12\|w_t - \bar{\nabla} f(x_t,y_t) -R_t\|^2}{\rho^2} + \frac{12\|R_t\|^2}{\rho^2}+ \frac{3}{\gamma^2}\|\tilde{x}_{t+1}-x_t\|^2 \big] \nonumber \\
 & \leq \frac{24G}{T}(m+T)^{1/2} + \frac{1}{T} \sum_{t=1}^T\frac{12\|R_t\|^2}{\rho^2} \nonumber \\
 & \leq \frac{24G}{T}(m+T)^{1/2} + \frac{12}{\rho^2T^2},
\end{align}
where the third inequality holds by the above inequality (\ref{eq:S7}),
and the last inequality holds by $\|R_t\|\leq \frac{1}{T}$ for all $t\geq 1$ by choosing $K=\frac{L_g}{\mu}\log(C_{gxy}C_{fy}T/\mu)$.

By combining the inequalities (\ref{eq:A77}) and (\ref{eq:A78}), we have
\begin{align}
\frac{1}{T}\sum_{t=1}^T\mathbb{E}\|\nabla F(x_t)\| & \leq  \sqrt{\frac{1}{T}\sum_{t=1}^T\mathbb{E}[\mathcal{M}_t^2]} \sqrt{\frac{1}{T}\sum_{t=1}^T\mathbb{E}\|A_t\|^2} \nonumber \\
& \leq \frac{\sqrt{\frac{1}{T}\sum_{t=1}^T\mathbb{E}\|A_t\|^2}}{\rho}\Big( \frac{2\sqrt{6G'm}}{T^{1/2}} + \frac{2\sqrt{6G'}}{T^{1/4}} + \frac{2\sqrt{3}}{T} \Big),
\end{align}
where $G'=\rho^2G$.

\end{proof}

\subsection{ Convergence Analysis of VR-BiAdam Algorithm }
\label{Appendix:A2}
In the subsection, we detail convergence analysis of VR-BiAdam algorithm.

\begin{lemma} \label{lem:C4}
Under the above Assumptions (\ref{ass:1}, \ref{ass:3}, \ref{ass:4}), assume the stochastic gradient estimators $v_t$ and $w_t$ be generated from Algorithm \ref{alg:2},
we have
 \begin{align}
 \mathbb{E}\|\nabla_y g(x_{t+1},y_{t+1}) - v_{t+1}\|^2
 & \leq (1-\alpha_{t+1})\mathbb{E} \|\nabla_y g(x_t,y_t) - v_t\|^2 + 2\alpha_{t+1}^2\sigma^2 \nonumber \\
 & \quad + 4L_g^2\eta_t^2\big(\mathbb{E}\|\tilde{x}_{t+1} - x_t\|^2 + \mathbb{E}\|\tilde{y}_{t+1} - y_t\|^2 \big),
 \end{align}
 \begin{align}
\mathbb{E}\|w_{t+1} - \bar{\nabla} f(x_{t+1},y_{t+1}) - R_{t+1}\|^2 & \leq (1-\beta_{t+1}) \mathbb{E}\|w_t - \bar{\nabla} f(x_t,y_t) -R_t\|^2 + 2\beta^2_{t+1}\sigma^2 \nonumber \\
 &  \quad + 4L^2_K\eta^2_t\big( \|\tilde{x}_{t+1}-x_t\|^2 + \|\tilde{y}_{t+1}-y_t\|^2 \big),
\end{align}
where $L_K^2 = 2L^2_f + 6C^2_{gxy}L^2_f\frac{K}{2\mu L_g - \mu^2} + 6C^2_{fy}L^2_{gxy}\frac{K}{2\mu L_g - \mu^2} + 6C^2_{gxy}L^2_f\frac{K^3L^2_{gyy}}{(L_g-\mu)^2(2\mu L_g - \mu^2)}$.
\end{lemma}

\begin{proof}
Without loss of generality, we only prove the term $\mathbb{E}\|w_{t+1} - \bar{\nabla} f(x_{t+1},y_{t+1})-R_{t+1}\|^2$.
The other term is similar for this term.
Since $w_{t+1} = \bar{\nabla} f(x_{t+1},y_{t+1};\bar{\xi}_{t+1}) + (1-\beta_{t+1})\big(w_t - \bar{\nabla} f(x_t,y_t;\bar{\xi}_{t+1})\big)$, we have
\begin{align}
 &\mathbb{E}\|w_{t+1} - \bar{\nabla} f(x_{t+1},y_{t+1}) - R_{t+1}\|^2 \\
 & = \mathbb{E}\|\bar{\nabla} f(x_{t+1},y_{t+1};\bar{\xi}_{t+1}) + (1-\beta_{t+1})\big(w_t - \bar{\nabla} f(x_t,y_t;\bar{\xi}_{t+1})\big)  - \bar{\nabla} f(x_{t+1},y_{t+1}) - R_{t+1}\|^2 \nonumber \\
 & = \mathbb{E}\|(1-\beta_{t+1})(w_t - \bar{\nabla} f(x_t,y_t) -R_t) + \beta_{t+1}\big( \bar{\nabla} f(x_{t+1},y_{t+1};\bar{\xi}_{t+1}) - \bar{\nabla} f(x_{t+1},y_{t+1}) - R_{t+1} \big) \nonumber \\
 & \quad + (1-\beta_{t+1})\big( \bar{\nabla} f(x_{t+1},y_{t+1};\bar{\xi}_{t+1}) - \bar{\nabla} f(x_{t+1},y_{t+1}) - R_{t+1} -(\bar{\nabla} f(x_t,y_t;\bar{\xi}_t)) - \bar{\nabla} f(x_t,y_t) - R_t) \big)\|^2 \nonumber \\
 & =(1-\beta_{t+1})^2 \mathbb{E}\|w_t - \bar{\nabla} f(x_t,y_t) -R_t\|^2 + \mathbb{E}\|\beta_{t+1}\big( \bar{\nabla} f(x_{t+1},y_{t+1};\bar{\xi}_{t+1}) - \bar{\nabla} f(x_{t+1},y_{t+1}) - R_{t+1} \big) \nonumber \\
 & \quad + (1-\beta_{t+1})\big( \bar{\nabla} f(x_{t+1},y_{t+1};\bar{\xi}_{t+1}) - \bar{\nabla} f(x_{t+1},y_{t+1}) - R_{t+1} -(\bar{\nabla} f(x_t,y_t;\bar{\xi}_t) - \bar{\nabla} f(x_t,y_t) - R_t) \big)\|^2 \nonumber \\
 & \leq (1-\beta_{t+1})^2 \mathbb{E}\|w_t - \bar{\nabla} f(x_t,y_t) -R_t\|^2 + 2\beta^2_{t+1}\mathbb{E}\| \bar{\nabla} f(x_{t+1},y_{t+1};\bar{\xi}_{t+1}) - \bar{\nabla} f(x_{t+1},y_{t+1}) - R_{t+1} \|^2 \nonumber \\
 & \quad + 2(1-\beta_{t+1})^2\|\bar{\nabla} f(x_{t+1},y_{t+1};\bar{\xi}_{t+1}) - \bar{\nabla} f(x_{t+1},y_{t+1}) - R_{t+1} -(\bar{\nabla} f(x_t,y_t;\bar{\xi}_t)) - \bar{\nabla} f(x_t,y_t) - R_t) \|^2 \nonumber \\
 & \leq (1-\beta_{t+1})^2 \mathbb{E}\|w_t - \bar{\nabla} f(x_t,y_t) -R_t\|^2 + 2\beta^2_{t+1}\sigma^2 +  2(1-\beta_{t+1})^2\|\bar{\nabla} f(x_{t+1},y_{t+1};\bar{\xi}_{t+1})) - \bar{\nabla} f(x_t,y_t;\bar{\xi}_t) \|^2 \nonumber \\
 & \leq (1-\beta_{t+1})^2 \mathbb{E}\|w_t - \bar{\nabla} f(x_t,y_t) -R_t\|^2 + 2\beta^2_{t+1}\sigma^2 +  4(1-\beta_{t+1})^2L^2_K\big(
  \|x_{t+1}-x_t\|^2 + \|y_{t+1}-y_t\|^2 \big) \nonumber \\
 & \leq (1-\beta_{t+1}) \mathbb{E}\|w_t - \bar{\nabla} f(x_t,y_t) -R_t\|^2 + 2\beta^2_{t+1}\sigma^2 + 4L^2_K\eta^2_t\big(
  \|\tilde{x}_{t+1}-x_t\|^2 + \|\tilde{y}_{t+1}-y_t\|^2 \big) \nonumber,
\end{align}
where the third equality holds by $\mathbb{E}_{\bar{\xi}} \big[\bar{\nabla} f(x_{t+1},y_{t+1};\bar{\xi}_{t+1})\big]=\bar{\nabla} f(x_{t+1},y_{t+1})+ R_{t+1}$ and $\mathbb{E}_{\bar{\xi}} \big[\bar{\nabla} f(x_t,y_t;\bar{\xi}_t))\big]= \bar{\nabla} f(x_t,y_t) + R_t$;
the third last inequality holds by the inequality $\mathbb{E}\|\zeta-\mathbb{E}[\zeta]\|^2 \leq \mathbb{E}\|\zeta\|^2$;
the second last inequality is due to Lemma \ref{lem:3}; the last inequality holds by $0<\beta_{t+1} \leq 1$ and $x_{t+1}=x_t+\eta_t(\tilde{x}_{t+1}-x_t)$,
$y_{t+1}=y_t+\eta_t(\tilde{y}_{t+1}-y_t)$.

\end{proof}

\begin{theorem}  \label{th:A3}
(Restatement of Theorem 3)
 Under the above Assumptions (\ref{ass:1}, \ref{ass:3}, \ref{ass:4}, \ref{ass:6}, \ref{ass:7}), in the Algorithm \ref{alg:2}, given $\mathcal{X}\subset\mathbb{R}^{d}$, $\eta_t=\frac{k}{(m+t)^{1/3}}$ for all $t\geq 0$, $\alpha_{t+1}=c_1\eta_t^2$, $\beta_{t+1}=c_2\eta_t^2$, $m\geq \max\big(2,k^3, (c_1k)^3,(c_2k)^3\big)$, $k>0$, $c_1 \geq \frac{2}{3k^3} + \frac{125L^2_0}{6\mu^2}$, $c_2 \geq \frac{2}{3k^3} + \frac{9}{2}$, $0<\lambda\leq \min\big( \frac{15b_l L^2_0}{16L^2_2\mu}, \frac{b_l}{6L_g} \big)$, $0< \gamma \leq \min\big(\frac{\sqrt{6}\lambda\mu\rho}{2\sqrt{24L^2_2\lambda^2\mu^2 + 125b^2_uL^2_0\kappa^2}}, \frac{m^{1/3}\rho}{4Lk}\big)$ and $K=\frac{L_g}{\mu}\log(C_{gxy}C_{fy}T/\mu)$, we have
\begin{align}
 \frac{1}{T} \sum_{t=1}^T \mathbb{E}||\mathcal{G}_\mathcal{X}(x_t,\nabla F(x_t),\gamma)|| \leq \frac{1}{T} \sum_{t=1}^T \mathbb{E}[\mathcal{M}_t]
 \leq \frac{2\sqrt{3M}m^{1/6}}{T^{1/2}} + \frac{2\sqrt{3M}}{T^{1/3}} + \frac{\sqrt{2}}{T},
\end{align}
where $M = \frac{F(x_1) - F^*}{\rho k\gamma} + \frac{5b_1L^2_0\Delta_0}{ \rho^2 k\lambda\mu} + \frac{2m^{1/3}\sigma^2}{\rho^2 k^2} + \frac{2k^2(c^2_1+c^2_2)\sigma^2\ln(m+T)}{\rho^2} + \frac{6k(m+T)^{1/3}}{\rho^2 T}$, $\Delta_0 = \|y_1-y^*(x_1)\|^2$ and $L^2_2 = L^2_g+L^2_K$.
\end{theorem}

\begin{proof}
Since $\eta_t=\frac{k}{(m+t)^{1/3}}$ on $t$ is decreasing and $m\geq k^3$, we have $\eta_t \leq \eta_0 = \frac{k}{m^{1/3}} \leq 1$ and $\gamma \leq \frac{m^{1/3}\rho}{4Lk}\leq \frac{\rho}{2L\eta_0} \leq \frac{\rho}{2L\eta_t}$ for any $t\geq 0$.
 Due to $0 < \eta_t \leq 1$ and $m\geq (c_1k)^3$, we have $\alpha_{t+1} = c_1\eta_t^2 \leq c_1\eta_t \leq \frac{c_1k}{m^{1/3}}\leq 1$.
 Similarly, due to $m\geq (c_2k)^3$, we have $\beta_{t+1}\leq 1$.
 According to Lemma \ref{lem:C4}, we have
 \begin{align}
  & \frac{1}{\eta_t}\mathbb{E} \|\nabla_y g(x_{t+1},y_{t+1}) - v_{t+1}\|^2 - \frac{1}{\eta_{t-1}}\mathbb{E} \|\nabla_y g(x_t,y_t) - v_t\|^2  \\
  & \leq \big(\frac{1-\alpha_{t+1}}{\eta_t} - \frac{1}{\eta_{t-1}}\big)\mathbb{E} \|\nabla_y g(x_t,y_t) - v_t\|^2 + 4L^2_g\eta_t\big(\|\tilde{x}_{t+1}-x_t\|^2 + \|\tilde{y}_{t+1}-y_t\|^2\big) + \frac{2\alpha_{t+1}^2\sigma^2}{\eta_t}\nonumber \\
  & = \big(\frac{1}{\eta_t} - \frac{1}{\eta_{t-1}} - c_1\eta_t\big)\mathbb{E} \|\nabla_y g(x_t,y_t) - v_t\|^2  + 4L^2_g\eta_t\big(\|\tilde{x}_{t+1}-x_t\|^2 +\|\tilde{y}_{t+1}-y_t\|^2\big) + 2c_1^2\eta^3_t\sigma^2, \nonumber
 \end{align}
 where the second equality is due to $\alpha_{t+1}=c_1\eta^2_t$.
 Similarly, we have
 \begin{align}
 & \frac{1}{\eta_t}\mathbb{E}\|w_{t+1} - \bar{\nabla} f(x_{t+1},y_{t+1}) - R_{t+1}\|^2 - \frac{1}{\eta_{t-1}}\mathbb{E}\|w_t - \bar{\nabla} f(x_t,y_t) - R_t\|^2 \\
 & \leq \big(\frac{1-\beta_{t+1}}{\eta_t}  - \frac{1}{\eta_{t-1}}\big) \mathbb{E}\|w_t - \bar{\nabla} f(x_t,y_t) -R_t\|^2 + 4L^2_K\eta_t\big( \|\tilde{x}_{t+1}-x_t\|^2 + \|\tilde{y}_{t+1}-y_t\|^2 \big) + \frac{2\beta^2_{t+1}\sigma^2}{\eta_t} \nonumber \\
 & =  \big(\frac{1}{\eta_t}  - \frac{1}{\eta_{t-1}} -c_2\eta_t\big) \mathbb{E}\|w_t - \bar{\nabla} f(x_t,y_t) -R_t\|^2 + 4L^2_K\eta_t\big( \|\tilde{x}_{t+1}-x_t\|^2 + \|\tilde{y}_{t+1}-y_t\|^2 \big) + 2c^2_2\eta^3_t\sigma^2. \nonumber
 \end{align}
By $\eta_t = \frac{k}{(m+t)^{1/3}}$, we have
 \begin{align}
  \frac{1}{\eta_t} - \frac{1}{\eta_{t-1}} &= \frac{1}{k}\big( (m+t)^{\frac{1}{3}} - (m+t-1)^{\frac{1}{3}}\big) \leq \frac{1}{3k(m+t-1)^{2/3}} \leq \frac{1}{3k\big(m/2+t\big)^{2/3}} \nonumber \\
  & \leq \frac{2^{2/3}}{3k(m+t)^{2/3}} = \frac{2^{2/3}}{3k^3}\frac{k^2}{(m+t)^{2/3}} = \frac{2^{2/3}}{3k^3}\eta_t^2 \leq \frac{2}{3k^3}\eta_t,
 \end{align}
 where the first inequality holds by the concavity of function $f(x)=x^{1/3}$, \emph{i.e.}, $(x+y)^{1/3}\leq x^{1/3} + \frac{y}{3x^{2/3}}$; the second inequality is due to $m\geq 2$,  and
 the last inequality is due to $0<\eta_t\leq 1$.

Let $c_1 \geq \frac{2}{3k^3} + \frac{125L^2_0}{6\mu^2}$, we have
 \begin{align} \label{eq:W1}
  & \frac{1}{\eta_t}\mathbb{E} \|\nabla_y g(x_{t+1},y_{t+1}) - v_{t+1}\|^2 - \frac{1}{\eta_{t-1}}\mathbb{E} \|\nabla_y g(x_t,y_t) - v_t\|^2  \\
  & \leq -\frac{125L^2_0\eta_t}{6\mu^2} \mathbb{E} \|\nabla_y g(x_t,y_t) - v_t\|^2 + 4L^2_g\eta_t\big(\|\tilde{x}_{t+1}-x_t\|^2 +\|\tilde{y}_{t+1}-y_t\|^2\big) + 2c_1^2\eta^3_t\sigma^2. \nonumber
 \end{align}
Let $c_2 \geq \frac{2}{3k^3} + \frac{9}{2}$, we have
 \begin{align} \label{eq:W2}
  & \frac{1}{\eta_t}\mathbb{E}\|w_{t+1} - \bar{\nabla} f(x_{t+1},y_{t+1}) - R_{t+1}\|^2 - \frac{1}{\eta_{t-1}}\mathbb{E}\|w_t - \bar{\nabla} f(x_t,y_t) - R_t\|^2  \\
  & \leq -\frac{9\eta_t}{2} \mathbb{E}\|w_t - \bar{\nabla} f(x_t,y_t) - R_t\|^2 +
  4L^2_K\eta_t\big(\|\tilde{x}_{t+1}-x_t\|^2 +\|\tilde{y}_{t+1}-y_t\|^2\big)  + 2c_2^2\eta_t^3\sigma^2. \nonumber
 \end{align}

According to Lemmas \ref{lem:A2} and \ref{lem:A4}, we have
\begin{align} \label{eq:W3}
  & F(x_{t+1}) - F(x_t) \\
  & \leq \frac{2\eta_t\gamma}{\rho}\|w_t-\bar{\nabla} f(x_t,y_t)\|^2 + \frac{L^2_0\eta_t\gamma}{\rho}\|y^*(x_t)-y_t\|^2 -\frac{\rho\eta_t}{2\gamma}\|\tilde{x}_{t+1}-x_t\|^2 \nonumber \\
  & \leq \frac{4\eta_t\gamma}{\rho}\|w_t-\bar{\nabla} f(x_t,y_t)-R_t\|^2 + \frac{4\eta_t\gamma}{\rho}\|R_t\|^2 + \frac{L^2_0\eta_t\gamma}{\rho}\|y^*(x_t)-y_t\|^2 -\frac{\rho\eta_t}{2\gamma}\|\tilde{x}_{t+1}-x_t\|^2. \nonumber
\end{align}
According to Lemma \ref{lem:A5}, we have
\begin{align} \label{eq:W4}
  & \|y_{t+1}-y^*(x_{t+1})\|^2 - \|y_t-y^*(x_t)\|^2 \\
  & \leq -\frac{\eta_t\mu\lambda}{4b_t}\|y_t -y^*(x_t)\|^2 -\frac{3\eta_t}{4} \|\tilde{y}_{t+1}-y_t\|^2
  + \frac{25\eta_t\lambda}{6\mu b_t} \|\nabla_y g(x_t,y_t)-v_t\|^2 + \frac{25\kappa^2\eta_t b_t}{6\mu\lambda}\|\tilde{x}_{t+1} - x_t\|^2. \nonumber
\end{align}

Next, we define a \emph{Lyapunov} function, for any $t\geq 1$
\begin{align}
 \Theta_t & = \mathbb{E}\big [F(x_t) + \frac{5b_t L^2_0\gamma}{\lambda\mu\rho}\|y_t-y^*(x_t)\|^2 + \frac{\gamma}{\rho\eta_{t-1}} \big(\|v_t - \nabla_y g(x_t,y_t)\|^2 + \|w_t-\bar{\nabla} f(x_t,y_t)-R_t\|^2 \big) \big]. \nonumber
\end{align}
For notational simplicity, let $L^2_2 = L^2_g+L^2_K$.
Then we have
 \begin{align}
 & \Theta_{t+1} - \Theta_t \nonumber \\
 & = F(x_{t+1}) - F(x_t) + \frac{5b_t L^2_0\gamma}{\lambda\mu\rho}
 \big(\|y_{t+1}-y^*(x_{t+1})\|^2 - \|y_t-y^*(x_t)\|^2 \big)
 + \frac{\gamma}{\rho} \bigg( \frac{1}{\eta_t}\mathbb{E}\|v_{t+1} - \nabla_y g(x_{t+1},y_{t+1})\|^2 \nonumber \\
 & \quad - \frac{1}{\eta_{t-1}}\mathbb{E}\|v_t - \nabla_y g(x_t,y_t)\|^2 + \frac{1}{\eta_t}\mathbb{E}\|w_{t+1}-\bar{\nabla} f(x_{t+1},y_{t+1})-R_{t+1}\|^2
 - \frac{1}{\eta_{t-1}}\mathbb{E}\|w_t-\bar{\nabla} f(x_t,y_t)-R_t\|^2 \bigg) \nonumber \\
 & \leq \frac{L^2_0\eta_t\gamma}{\rho}\|y^*(x_t)-y_t\|^2 + \frac{4\eta_t\gamma}{\rho}\|w_t-\bar{\nabla} f(x_t,y_t)-R_t\|^2 + \frac{4\eta_t\gamma}{\rho}\|R_t\|^2 -\frac{\rho\eta_t}{2\gamma}\|\tilde{x}_{t+1}-x_t\|^2  \nonumber \\
 & \quad + \frac{5b_t L^2_0\gamma}{\lambda\mu\rho} \bigg( -\frac{\eta_t\mu\lambda}{4b_t}\|y_t -y^*(x_t)\|^2 -\frac{3\eta_t}{4} \|\tilde{y}_{t+1}-y_t\|^2 + \frac{25\eta_t\lambda}{6\mu b_t} \|\nabla_y g(x_t,y_t)-v_t\|^2
 + \frac{25\kappa^2\eta_tb_t}{6\mu\lambda}\|\tilde{x}_{t+1} - x_t\|^2 \bigg)  \nonumber \\
 & \quad + \frac{\gamma}{\rho} \bigg( -\frac{125L^2_0\eta_t}{6\mu^2} \mathbb{E} \|\nabla_y g(x_t,y_t) - v_t\|^2 + 4L^2_g\eta_t\big(\|\tilde{x}_{t+1}-x_t\|^2 +\|\tilde{y}_{t+1}-y_t\|^2\big) + 2c_1^2\eta^3_t\sigma^2 \nonumber \\
 & \quad  -\frac{9\eta_t}{2} \mathbb{E}\|w_t - \bar{\nabla} f(x_t,y_t) - R_t\|^2 +
  4L^2_K\eta_t\big(\|\tilde{x}_{t+1}-x_t\|^2 +\|\tilde{y}_{t+1}-y_t\|^2\big)  + 2c_2^2\eta_t^3\sigma^2 \bigg)\nonumber \\
 & = - \frac{\gamma\eta_t}{4\rho} \bigg( L^2_0 \|y_t - y^*(x_t)\|^2 + 2\mathbb{E}\|w_t - \bar{\nabla} f(x_t,y_t) - R_t\|^2 \bigg)
 - \big( \frac{\rho}{2\gamma} - \frac{4L^2_2\gamma}{\rho} - \frac{125b^2_t L^2_0\kappa^2\gamma}{6\mu^2\lambda^2\rho}\big)\eta_t\|\tilde{x}_{t+1}-x_t\|^2  \nonumber \\
 & \quad - \big( \frac{15b_t L^2_0\gamma}{4\lambda\mu\rho} -
 \frac{4L^2_2\gamma}{\rho}\big)\eta_t\|\tilde{y}_{t+1}-y_t\|^2 +\frac{4\eta_t\gamma}{\rho}\|R_t\|^2
 + \frac{2(c^2_1+c^2_2)\gamma\sigma^2}{\rho}\eta^3_t \nonumber \\
 & \leq - \frac{\gamma\eta_t}{4\rho} \bigg( L^2_0 \|y_t - y^*(x_t)\|^2 + 2\mathbb{E}\|w_t - \bar{\nabla} f(x_t,y_t) - R_t\|^2 \bigg) - \frac{\rho\eta_t}{4\gamma}\|\tilde{x}_{t+1}-x_t\|^2 \nonumber \\
 & \quad +\frac{4\eta_t\gamma}{\rho}\|R_t\|^2 + \frac{2(c^2_1+c^2_2)\gamma\sigma^2}{\rho}\eta^3_t,
 \end{align}
where the first inequality holds by the above inequalities (\ref{eq:W1}), (\ref{eq:W2}), (\ref{eq:W3}) and (\ref{eq:W4});
the last inequality is due to $0< \gamma \leq \frac{\sqrt{6}\lambda\mu\rho}{2\sqrt{24L^2_2\lambda^2\mu^2 + 125b^2_u L^2_0\kappa^2}}
\leq \frac{\sqrt{6}\lambda\mu\rho}{2\sqrt{24L^2_2\lambda^2\mu^2 + 125b^2_t L^2_0\kappa^2}}$, $0 < \lambda \leq \frac{15b_l L^2_0}{16L^2_2\mu} \leq \frac{15b_t L^2_0}{16L^2_2\mu}$ for all $t\geq 1$.

Let $\Phi_t = L^2_0\|y_t-y^*(x_t)\|^2 + 2\|w_t-\bar{\nabla}f(x_t,y_t)-R_t\|^2$, then we have
\begin{align} \label{eq:W5}
 \frac{\gamma\eta_t}{4\rho}\mathbb{E}\Big[ \Phi_t + \frac{\rho\eta_t}{4\gamma}\|\tilde{x}_{t+1}-x_t\|^2 \Big] \leq \Theta_t - \Theta_{t+1} + \frac{4\eta_t\gamma}{\rho}\|R_t\|^2 + \frac{2(c^2_1+c^2_2)\gamma\sigma^2}{\rho}\eta^3_t.
\end{align}
Taking average over $t=1,2,\cdots,T$ on both sides of (\ref{eq:W5}), we have
\begin{align}
  \frac{1}{T} \sum_{t=1}^T \mathbb{E} \Big[ \frac{\eta_t}{4}\Phi_t + \frac{\rho^2\eta_t}{4\gamma^2}\|\tilde{x}_{t+1}-x_t\|^2 \Big]
  \leq  \sum_{t=1}^T \frac{\rho(\Theta_t - \Theta_{t+1})}{T\gamma} + \frac{4}{T}\sum_{t=1}^T\eta_t\|R_t\|^2 + \frac{2(c^2_1+c^2_2)\sigma^2}{T}\sum_{t=1}^T \eta^3_t. \nonumber
\end{align}

Given $x_1\in \mathcal{X}$ and $y_1 \in \mathcal{Y}$, let $\Delta_0 = \|y_1-y^*(x_1)\|^2$, we have
\begin{align} \label{eq:W6}
 \Theta_1 &= \mathbb{E}\big [F(x_1) + \frac{5b_1 L^2_0\gamma}{\lambda\mu\rho}\|y_1-y^*(x_1)\|^2 + \frac{\gamma}{\rho\eta_{0}} \big(\|v_1 - \nabla_y g(x_1,y_1)\|^2 + \|w_1-\bar{\nabla} f(x_1,y_1)-R_1\|^2 \big) \big] \nonumber \\
 & \leq F(x_1) + \frac{5b_1L^2_0\gamma\Delta_0}{\lambda\mu\rho} + \frac{2\gamma\sigma^2}{\rho\eta_0},
\end{align}
where the last inequality holds by Assumption \ref{ass:2}.
Since $\eta_t$ is decreasing, i.e., $\eta_T^{-1} \geq \eta_t^{-1}$ for any $0\leq t\leq T$, we have
 \begin{align}
 & \frac{1}{T} \sum_{t=1}^T \mathbb{E}\Big[ \frac{\Phi_t}{4} + \frac{\rho^2}{4\gamma^2}\|\tilde{x}_{t+1}-x_t\|^2 \Big]  \\
 & \leq  \frac{\rho}{T\gamma\eta_T} \sum_{t=1}^T\big(\Theta_t - \Theta_{t+1}\big) + \frac{2(c^2_1+c^2_2)\sigma^2}{T\eta_T}\sum_{t=1}^T\eta^3_t + \frac{4}{T}\sum_{t=1}^T\eta_t\|R_t\|^2 \nonumber \\
 & \leq \frac{\rho}{T\gamma\eta_T} \Big( F(x_1) + \frac{5b_1 L^2_0\gamma\Delta_0}{\lambda\mu\rho} + \frac{2\gamma\sigma^2}{\rho\eta_0} - F^* \Big) + \frac{2(c^2_1+c^2_2)\sigma^2}{T\eta_T}\sum_{t=1}^T\eta^3_t + \frac{4}{T^2}\sum_{t=1}^T\eta_t  \nonumber \\
 & \leq \frac{\rho(F(x_1) - F^*)}{T\gamma\eta_T} + \frac{5b_1 L^2_0\Delta_0}{\lambda\mu\eta_TT} + \frac{2\sigma^2}{\eta_0\eta_TT} + \frac{2(c^2_1+c^2_2)\sigma^2}{T\eta_T}\int^T_1\frac{k^3}{m+t} dt + \frac{4}{T^2}\int_{1}^T\frac{k}{(m+t)^{1/3}} dt\nonumber \\
 & \leq  \frac{\rho(F(x_1) - F^*)}{T\gamma\eta_T} + \frac{5b_1 L^2_0\Delta_0}{\lambda\mu\eta_TT} + \frac{2\sigma^2}{\eta_0\eta_TT} + \frac{2k^3(c^2_1+c^2_2)\sigma^2}{T\eta_T}\ln(m+T) + \frac{6k}{T^2}(m+T)^{2/3} \nonumber \\
 & = \bigg( \frac{\rho(F(x_1) - F^*)}{k\gamma} + \frac{5b_1 L^2_0\Delta_0}{ k\lambda\mu} + \frac{2m^{1/3}\sigma^2}{k^2} + 2k^2(c^2_1+c^2_2)\sigma^2\ln(m+T) + \frac{6k(m+T)^{1/3}}{T} \bigg)\frac{(m+T)^{1/3}}{T}, \nonumber
\end{align}
where the second inequality holds by the above inequality (\ref{eq:W6}).
Let $M = \frac{F(x_1) - F^*}{\rho k\gamma} + \frac{5b_1 L^2_0\Delta_0}{ \rho^2 k\lambda\mu} + \frac{2m^{1/3}\sigma^2}{\rho^2 k^2} + \frac{2k^2(c^2_1+c^2_2)\sigma^2\ln(m+T)}{\rho^2} + \frac{6k(m+T)^{1/3}}{\rho^2 T}$,
we have
\begin{align} \label{eq:W7}
 \frac{1}{T} \sum_{t=1}^T \mathbb{E}\Big[ \frac{\Phi_t}{4\rho^2} + \frac{1}{4\gamma^2}\|\tilde{x}_{t+1}-x_t\|^2 \Big]  \leq \frac{M}{T}(m+T)^{1/3}.
\end{align}
According to Jensen's inequality, we have
\begin{align}
& \frac{1}{T} \sum_{t=1}^T \mathbb{E}\Big[ \frac{1}{2\gamma}\|\tilde{x}_{t+1}-x_t\|  + \frac{1}{2\rho} \big( L_0 \|y_t - y^*(x_t)\| + \sqrt{2}\|w_t - \bar{\nabla} f(x_t,y_t) -R_t\|\big) \Big] \nonumber \\
& \leq \bigg( \frac{3}{T} \sum_{t=1}^T \big(  \frac{1}{4\gamma^2}\|\tilde{x}_{t+1}-x_t\|^2 + \frac{L^2_0}{4\rho^2}\|y_t - y^*(x_t)\|^2 + \frac{2}{4\rho^2}\mathbb{E}\|w_t - \bar{\nabla} f(x_t,y_t) -R_t\|^2  \big)\bigg)^{1/2} \nonumber \\
& = \bigg( \frac{3}{T} \sum_{t=1}^T \big( \frac{\Phi_t}{4\rho^2} + \frac{1}{4\gamma^2}\|\tilde{x}_{t+1}-x_t\|^2 \big)\bigg)^{1/2} \nonumber \\
& \leq \frac{\sqrt{3M}}{T^{1/2}}(m+T)^{1/6} \leq \frac{\sqrt{3M}m^{1/6}}{T^{1/2}} + \frac{\sqrt{3M}}{T^{1/3}},
\end{align}
where the last inequality is due to $(a+b)^{1/6} \leq a^{1/6} + b^{1/6}$ for all $a,b>0$.
Thus we have
\begin{align} \label{eq:W8}
& \frac{1}{T} \sum_{t=1}^T \mathbb{E}\big[  \frac{1}{\gamma}\|\tilde{x}_{t+1}-x_t\| + \frac{1}{\rho}\big(L_0\|y_t - y^*(x_t)\| + \sqrt{2}\|w_t - \bar{\nabla} f(x_t,y_t) -R_t\|\big) \big] \nonumber \\
& \leq \frac{2\sqrt{3M}m^{1/6}}{T^{1/2}} + \frac{2\sqrt{3M}}{T^{1/3}}.
\end{align}

Since $\|w_t - \bar{\nabla} f(x_t,y_t) -R_t\|\geq \|w_t - \bar{\nabla} f(x_t,y_t)\| - \|R_t\|$,
by the above inequality (\ref{eq:W8}), we can obtain
\begin{align}
 & \frac{1}{T} \sum_{t=1}^T \mathbb{E}\big[ \frac{1}{\gamma}\|\tilde{x}_{t+1}-x_t\| + \frac{1}{\rho}\big(L_0\|y_t - y^*(x_t)\| + \sqrt{2}\|w_t - \bar{\nabla} f(x_t,y_t) -R_t\|\big) \big] \nonumber \\
& \leq \frac{2\sqrt{3M}m^{1/6}}{T^{1/2}} + \frac{2\sqrt{3M}}{T^{1/3}} + \frac{\sqrt{2}}{T} \sum_{t=1}^T \|R_t\| \nonumber \\
& = \frac{2\sqrt{3M}m^{1/6}}{T^{1/2}} + \frac{2\sqrt{3M}}{T^{1/3}} + \frac{\sqrt{2}}{T},
\end{align}
where the last inequality is due to $\|R_t\|\leq \frac{1}{T}$ for all $t\geq 1$ by choosing $K=\frac{L_g}{\mu}\log(C_{gxy}C_{fy}T/\mu)$
in Algorithm \ref{alg:2}.

According to the above inequality (\ref{eq:19}), we have
\begin{align}
\frac{1}{T} \sum_{t=1}^T \mathbb{E}||\mathcal{G}_\mathcal{X}(x_t,\nabla F(x_t),\gamma)|| \leq \frac{1}{T} \sum_{t=1}^T \mathbb{E}[\mathcal{M}_t] \leq \frac{2\sqrt{3M}m^{1/6}}{T^{1/2}} + \frac{2\sqrt{3M}}{T^{1/3}} + \frac{\sqrt{2}}{T}.
\end{align}

\end{proof}

\begin{theorem}  \label{th:A4}
(Restatement of Theorem 4)
 Under the above Assumptions (\ref{ass:1}, \ref{ass:3}, \ref{ass:4}, \ref{ass:6}, \ref{ass:7}), in the Algorithm \ref{alg:2}, given $\mathcal{X}=\mathbb{R}^{d}$, $\eta_t=\frac{k}{(m+t)^{1/3}}$ for all $t\geq 0$, $\alpha_{t+1}=c_1\eta_t^2$, $\beta_{t+1}=c_2\eta_t^2$, $m\geq \max\big(2,k^3, (c_1k)^3,(c_2k)^3\big)$, $k>0$, $c_1 \geq \frac{2}{3k^3} + \frac{125L^2_0}{6\mu^2}$, $c_2 \geq \frac{2}{3k^3} + \frac{9}{2}$, $0<\lambda\leq \min\big( \frac{15b_l L^2_0}{16L^2_2\mu}, \frac{b_l}{6L_g} \big)$, $0< \gamma \leq \min\Big(\frac{\sqrt{6}\lambda\mu\rho}{2\sqrt{24L^2_2\lambda^2\mu^2 + 125b^2_u L^2_0\kappa^2}}, \frac{m^{1/3}\rho}{4Lk}\Big)$ and $K=\frac{L_g}{\mu}\log(C_{gxy}C_{fy}T/\mu)$, we have
\begin{align}
 \frac{1}{T}\sum_{t=1}^T\mathbb{E}\|\nabla F(x_t)\| \leq \frac{\sqrt{\frac{1}{T}\sum_{t=1}^T\mathbb{E}\|A_t\|^2}}{\rho}\Big( \frac{2\sqrt{6M'm}}{T^{1/2}} + \frac{2\sqrt{6M'}}{T^{1/3}} + \frac{2\sqrt{3}}{T}\Big),
\end{align}
where $M' = \frac{\rho(F(x_1) - F^*)}{k\gamma} + \frac{5b_1L^2_0\Delta_0}{k\lambda\mu} + \frac{2m^{1/3}\sigma^2}{k^2} + 2k^2(c^2_1+c^2_2)\sigma^2\ln(m+T) + \frac{6k(m+T)^{1/3}}{T}$.
\end{theorem}

\begin{proof}
According to Lemma \ref{lem:A2}, we have
\begin{align}
\mathcal{M}_t & =\frac{1}{\gamma}\|x_t - \tilde{x}_{t+1}\|
+ \frac{1}{\rho}\big( L_0 \|y^*(x_t) - y_t\|  + \sqrt{2}\|\bar{\nabla} f(x_t,y_t)- w_t\|\big) \nonumber \\
& \geq \frac{1}{\gamma}\|x_t - \tilde{x}_{t+1}\|
+ \frac{1}{\rho}\|\nabla  F(x_t) - w_t\| \nonumber \\
& \mathop{=}^{(i)} \|A_t^{-1}w_t\|
+ \frac{1}{\rho}\|\nabla  F(x_t)-w_t\| \nonumber \\
& = \frac{1}{\|A_t\|}\|A_t\|\|A_t^{-1}w_t\|
+ \frac{1}{\rho}\|\nabla  F(x_t)-w_t\| \nonumber \\
& \geq \frac{1}{\|A_t\|}\|w_t\|
+ \frac{1}{\rho}\|\nabla  F(x_t)-w_t\| \nonumber \\
& \mathop{\geq}^{(ii)} \frac{1}{\|A_t\|}\|w_t\|
+ \frac{1}{\|A_t\|}\|\nabla  F(x_t)-w_t\| \nonumber \\
& \geq \frac{1}{\|A_t\|}\|\nabla  F(x_t)\|,
\end{align}
where the equality $(i)$ holds by $\tilde{x}_{t+1} = x_t - \gamma A^{-1}_tw_t$ that can be easily obtained from the step 5 of Algorithm \ref{alg:2} when $\mathcal{X}=\mathbb{R}^{d}$, and the inequality $(ii)$ holds by $\|A_t\| \geq \rho$ for all $t\geq1$ due to Assumption \ref{ass:7}. Then we have
\begin{align}
 \|\nabla  F(x_t)\| \leq \|A_t\|\mathcal{M}_t.
\end{align}
According to Cauchy-Schwarz inequality, we have
\begin{align} \label{eq:A102}
\frac{1}{T}\sum_{t=1}^T\mathbb{E}\|\nabla F(x_t)\| \leq \frac{1}{T}\sum_{t=1}^T\mathbb{E}\big[\mathcal{M}_t \|A_t\|\big] \leq \sqrt{\frac{1}{T}\sum_{t=1}^T\mathbb{E}[\mathcal{M}_t^2]} \sqrt{\frac{1}{T}\sum_{t=1}^T\mathbb{E}\|A_t\|^2}.
\end{align}
Then we have
\begin{align} \label{eq:A103}
 \frac{1}{T} \sum_{t=1}^T \mathbb{E}[\mathcal{M}_t^2] & \leq \frac{1}{T} \sum_{t=1}^T \mathbb{E}\Big[ \frac{3L^2_0 \|y_t - y^*(x_t)\|^2}{\rho^2} + \frac{6\|w_t - \bar{\nabla} f(x_t,y_t)\|^2}{\rho^2} + \frac{3}{\gamma^2}\|\tilde{x}_{t+1}-x_t\|^2 \Big] \nonumber \\
 & \leq \frac{1}{T} \sum_{t=1}^T \mathbb{E}\Big[ \frac{3L^2_0 \|y_t - y^*(x_t)\|^2}{\rho^2} + \frac{12\|w_t - \bar{\nabla} f(x_t,y_t) -R_t\|^2}{\rho^2} + \frac{12\|R_t\|^2}{\rho^2} + \frac{3}{\gamma^2}\|\tilde{x}_{t+1}-x_t\|^2 \Big] \nonumber \\
 & \leq \frac{24M}{T}(m+T)^{1/3} + \frac{1}{T} \sum_{t=1}^T\frac{12\|R_t\|^2}{\rho^2} \nonumber \\
 & \leq \frac{24M}{T}(m+T)^{1/3} + \frac{12}{\rho^2T^2},
\end{align}
where the third inequality holds by the above inequality (\ref{eq:W7}),
and the last inequality holds by $\|R_t\|\leq \frac{1}{T}$ for all $t\geq 1$ by choosing $K=\frac{L_g}{\mu}\log(C_{gxy}C_{fy}T/\mu)$.

By combining the above inequalities (\ref{eq:A102}) and (\ref{eq:A103}), we have
\begin{align}
\frac{1}{T}\sum_{t=1}^T\mathbb{E}\|\nabla F(x_t)\| & \leq  \sqrt{\frac{1}{T}\sum_{t=1}^T\mathbb{E}[\mathcal{M}_t^2]} \sqrt{\frac{1}{T}\sum_{t=1}^T\mathbb{E}\|A_t\|^2} \nonumber \\
& \leq \frac{\sqrt{\frac{1}{T}\sum_{t=1}^T\mathbb{E}\|A_t\|^2}}{\rho}\Big( \frac{2\sqrt{6M'm}}{T^{1/2}} + \frac{2\sqrt{6M'}}{T^{1/3}} + \frac{2\sqrt{3}}{T} \Big),
\end{align}
where $M'=\rho^2M$.

\end{proof}

\end{appendices}

\end{onecolumn}

\end{document}